\documentclass[12pt,a4paper]{article}

\usepackage{amsmath,amsthm,amsfonts,amssymb,oldgerm}
\usepackage[frenchb,english]{babel}
\usepackage[utf8x]{inputenc}
\usepackage{indentfirst}
\usepackage{url}
\numberwithin{equation}{section}
\usepackage[citecolor=blue,colorlinks=true]{hyperref}
\usepackage{color}
\allowdisplaybreaks[1]
\usepackage[top=3cm, bottom=3cm, right=2.2cm, left=3cm]{geometry}

\newcommand{\disp}{\displaystyle}

\newcommand{\E}{\mathds{E}}
\renewcommand{\H}{\mathds{H}}
\newcommand{\N}{\mathds{N}}
\renewcommand{\P}{\mathds{P}}

\newcommand{\R}{\mathds{R}}

\newcommand{\1}{\mathds{1}}

\newcommand\cA{{\mathcal A}}

\newcommand\cF{{\mathcal F}}

\newcommand\cR{{\mathcal R}}

\newcommand\fB{{\mathfrak B}}

\newcommand\fR{{\mathfrak R}}
\newcommand\fS{{\mathfrak S}}

\usepackage{dsfont}				

\usepackage{cleveref}

\newtheorem{theorem}{Theorem}[section]
\crefname{theorem}{Theorem}{Theorems}

\crefname{proof}{Proof}{Proofs}

\newtheorem{proposition}{Proposition}[section]
\crefname{proposition}{Proposition}{Propositions}

\newtheorem{lemma}{Lemma}[section]
\crefname{lemma}{Lemma}{Lemmas}

\crefname{equation}{equation}{equations}

\crefname{section}{Section}{Sections}

\newtheorem{remark}{Remark}[section]
\crefname{remark}{Remark}{Remarks}

\newtheorem{assumption}{Assumption}[section]
\crefname{assumption}{Assumption}{Assumptions}

\crefname{hypothesis}{Hypothesis}{Hypothesis}

\crefname{note}{Note}{Notes}

\renewcommand{\abstractname}

\begin{document}

\begin{center}
\textbf{A LAW OF LARGE NUMBERS IN THE SUPREMUM NORM FOR A MULTISCALE STOCHASTIC SPATIAL GENE NETWORK}\\
\end{center}

\begin{center}
\par\vspace{0.5cm}
\renewcommand*{\thefootnote}{$\dagger$}
\textsc{By Arnaud Debussche\footnotemark[1] and Mac Jugal Nguepedja Nankep\footnotemark[1] \footnotetext[1]{IRMAR, ENS Rennes, CNRS, Campus de Ker Lann, 37170 Bruz.}} 
\par\vspace{0.5cm}
\end{center}



\begin{abstract}
We study the asymptotic behavior of multiscale stochastic spatial gene networks. Multiscaling takes into account the difference of abundance between molecu\-les, and captures the dynamic of rare species at a mesoscopic level. We introduce an assumption of spatial correlations for reactions involving rare species and a new law of large numbers is obtained. According to the scales, the whole system splits into two parts with different but coupled dynamics. The high scale component converges to the usual spatial model which is the solution of a partial differential equation, whereas, the low scale component converges to the usual homogeneous model which is the solution of an ordinary differential equation. Comparisons are made in the supremum norm.
\end{abstract}

\section{Introduction}

Modern molecular biology emphasizes the important role of the gene regulatory networks in the functioning of living organisms (see \cite{Davidson2005}, \cite{Stathopoulos2005} or \cite{Lesley2011}). As many other dynamical systems in a great variety of fields (e.g. biology, chemistry, epidemic theory or physics), gene regulatory networks belong to the large family of chemical reaction systems, which has been highly investigated in mathematics, for modeling. We refer to \cite{Kurtz1970,Kurtz1971}, \cite{Arnold1980}, \cite{Arnold1980bis}, \cite{Kurtz1986}, \cite{P_Erdi1989}, \cite{Crudu2007}, \cite{Arnaud2009}, \cite{Arnaud2012}, and many other references therein. Following Arnold in \cite{Arnold1980bis}, there are two main criteria according to which reactions in a spatial domain are modeled:\vspace{0.2cm}\par 

\noindent \textbf{(\hypertarget{C1}{C1})} global description (spatially homogeneous, "well-stirred" case excluding diffusion) versus local description (spatially inhomogeneous, "spatial model" including diffusion);\vspace{0.1cm}\par 

\noindent \textbf{(\hypertarget{C2}{C2})} deterministic description (macroscopic, phenomenological, in terms of concentrations) versus stochastic description (mesoscopic, taking into account fluctuations, at the level of particles or molecules).\vspace{0.2cm}\par 

The combination of these two criteria gives rise to four mathematical models:\vspace{0.2cm}\par 

\noindent \hspace{0.2cm} $\bullet$ Deterministic-Homogeneous: the concentration is solution to an Ordinary Differential Equation (ODE).\par 
\noindent \hspace{0.2cm} $\bullet$ Stochastic-Homogeneous: the number of particles is a (time) jump Markov process.\par 
\noindent \hspace{0.2cm} $\bullet$ Deterministic-Spatial: the concentration is solution to a Partial Differential Equation (PDE).\par 
\noindent \hspace{0.2cm} $\bullet$ Stochastic-Spatial: the number of particles is a (space-time) jump Markov process. \vspace{0.2cm}\par 

There is a restrictive assumption for deterministic and many of the stochastic models, in the literature. Namely, the different species (or reactants) have the same order of population size, which is large, meaning that the whole system has a fast dynamic. There are many biological situations for which this is not true. In some situations (gene regulatory networks are typical), some molecules are present in much greater quantity than others, and reactions rate constants can vary over several orders of magnitude (see \cite{Arnaud2012} or \cite{Kurtz2006}). A third criterion of modeling then shows up: \vspace{0.2cm}\par 

\noindent \textbf{(\hypertarget{C3}{C3})} one scale description with a unique and large population size scale, with only high reaction rates and fast dynamics, versus multiscale description with at least two population size scales, with high and low reaction rates, and, fast and slow dynamics.\vspace{0.2cm}\par 

Now, with the three criteria (\hyperlink{C1}{C1}), (\hyperlink{C2}{C2}), and (\hyperlink{C3}{C3}) in hands, we can derive: the \vspace{0.2cm}\par 

\hspace{2cm}\textbf{(\hypertarget{M1}{M1})}\hspace{1cm} Deterministic Homogeneous Model,\par 
\hspace{2cm}\textbf{(\hypertarget{M2}{M2})}\hspace{1cm} Deterministic Spatial Model,\par 
\hspace{2cm}\textbf{(\hypertarget{M3}{M3})}\hspace{1cm} Stochastic Homogeneous Model,\par 
\hspace{2cm}\textbf{(\hypertarget{M4}{M4})}\hspace{1cm} Stochastic Spatial Model,\par 
\hspace{2cm}\textbf{(\hypertarget{M5}{M5})}\hspace{1cm} Multiscale Stochastic Homogeneous Model,\par 
\hspace{2cm}\textbf{(\hypertarget{M6}{M6})}\hspace{1cm} Multiscale Stochastic Spatial Model. \vspace{0.2cm}\par 

\noindent It should be emphasized that criterion (\hyperlink{C3}{C3}) has been introduced at the mesoscopic level and thus suggests to preferentially consider a stochastic description for multiscaling. This stresses the importance of, and our focus on, the multiscale stochastic modeling.

All these models are not independent one from another. The relation between (\hyperlink{M1}{M1}) and (\hyperlink{M3}{M3}) has been thoroughly investigated by Kurtz (among others, see \cite{Kurtz1970, Kurtz1971} or \cite{Kurtz1986} with Ethier). A law of large numbers (LLN) and the corresponding central limit theorem (CTL) have been proved, showing convergence of (\hyperlink{M3}{M3}) to (\hyperlink{M1}{M1}). The consistency of (\hyperlink{M1}{M1}) and (\hyperlink{M2}{M2}) as well as of (\hyperlink{M3}{M3}) and (\hyperlink{M4}{M4}) have been proved in \cite{Arnold1980bis}. In \cite{Arnold1980}, Arnold and Theodosopulu compared (\hyperlink{M2}{M2}) and (\hyperlink{M4}{M4}) in the $L^2$ norm, through a LLN. Blount after Kotelenez did the same comparison, but much more extensively. Under more and more relaxed assumptions, they proved different LLNs and the associated CLTs in spaces of distributions. Blount showed a LLN in the supremum norm. We refer to \cite{Kotelenez1986bis, Kotelenez1986, Kotelenez1987, Kotelenez1988}, \cite{Blount1987, Blount1992, Blount1993, Blount1994} among others.

In \cite{Crudu2007}, \cite{Arnaud2009} and \cite{Arnaud2012}, Crudu, Debussche, Muller and Radulescu studied (\hyperlink{M5}{M5}). In the latter paper, they proved that the model (weakly) converges to a finite-dimensional Piecewise Deterministic Markov Process (PDMP). In finite dimension, PDMPs are hybrid processes which follow between consecutive jumps, the flows of ODEs whose parameters can jump. They have been well formalised and studied by Davis in \cite{Davis1993}. Depending on the interactions and the scaling, Crudu \textit{et al.} distinguished various types of limiting PDMPs. Stochasticity does not disappear in that hybrid simplification, and multiscaling appears to be a tool for estimating asymptotically at the first order, the noise lost in the LLN. Furthermore, PDMPs increase computational efficiency. This is very useful, since direct simulation of Markov processes is extremely time consuming for models we are dealing with (see \cite{Arnaud2009} or \cite{Gillespie1976}).

In this article, we start the study of (\hyperlink{M6}{M6}) and generalize Crudu \textit{et al.} to a spatially dependent situation. Depending on their abundance, we distinguish abundant ("continuous") and rare ("discrete") species. Also, there are different speeds of reactions, according to their rate scale and the species involved. When only abundant (resp. rare) species are involved, the reaction is fast (resp. slow). Otherwise, the speed depends on the specification of the rate scale. We assume that only abundant reactants diffuse. Also, a spatial correlation for slow reactions is considered. Its effects are inversely proportional to the distance to the location where the reaction occurred. In this context, we prove a new law of large numbers, showing the convergence of (\hyperlink{M6}{M6}) to (\hyperlink{M2}{M2}) coupled with (\hyperlink{M1}{M1}), in the supremum norm.

The structure for the rest of the article is as follows. In Section 2, we present local models, starting with the existing ones. We emphasize the notions of infinitesimal generator, debit functions, scaling and density dependance. Then we present our model of interest (\hyperlink{M6}{M6}). Afterwards, we identify its limit through the asymptotics of its debit function and of its generator. We state and prove our main result in Section 3. Finally, Section 4 is dedicated to the proofs of some intermediate results.

\hspace{1cm}\\

\noindent \underline{\textbf{\textit{Some general notations.}}} Let $(Z,\Vert\cdot\Vert_Z)$ and $(\tilde{Z},\Vert\cdot\Vert_{\tilde{Z}}$ be Banach spaces. The product space $Z\times\tilde{Z}$ is equipped with the norm $\Vert\cdot\Vert_Z+\Vert\cdot\Vert_{\tilde{Z}}$. We introduce:
\begin{description}
\item[$\bullet$] $\mathcal{L}(Z,\tilde{Z})$: the space of continuous linear maps from $Z$ to $\tilde{Z}$. If $Z=\tilde{Z}$, one simply writes $\mathcal{L}(Z)$.
The operator norm is denoted $\Vert\cdot\Vert_{Z\rightarrow\tilde Z}$ and when there is no risk of confusion, we denote it $\Vert\cdot\Vert$.

\item[$\bullet$] $Z'$: the space of continuous linear forms on $Z$. It is the topological dual space of $Z$.

\item[$\bullet$] $\fB(Z)$ (resp. $\fB_b(Z)$): the space of Borel-measurable (resp. bounded Borel-measurable) real valued functions on $Z$. The space $\fB_b(Z)$ is endowed with the supremum norm \[ \Vert f\Vert_{\fB_b(Z)}=\sup_{x\in Z}|f(x)|=\Vert f\Vert_\infty. \]

\item[$\bullet$] $\disp C_b^k(Z)$, $k\in \N$: the space of real valued functions of class $C^k$, i.e. $k$-continuously Fr\'echet differentiable,  on $Z$ which are bounded and have uniformly bounded succesive differentials. It is equipped with the norm \[ \Vert f\Vert_{C_b^k(Z)}=\sum_{i=0}^k\big\Vert D^{i}f\big\Vert_\infty, \] where $D^{i}f$ is the $i$-th differential of $f\in C_b^k(Z)$, and $\disp C_b^0(Z)=C_b(Z)$ is the set of bounded continuous real valued functions on $Z$. \par 

\item[$\bullet$] $C^{l,k}(Z\times\tilde{Z})$, $l,k\in \N$: the set of real valued functions $\varphi$ of class $C^l$ w.r.t. the first variable and of class $C^k$ w.r.t. the second. In particular, $C^{0,0}(Z\times\tilde{Z)})=C(Z\times\tilde{Z)})$. \par 

For $(z,\tilde{z})\in Z\times\tilde{Z}$, we donote by $D^{l,k}\varphi(z,\tilde{z})$ the (Fr\'echet) differential of $\varphi$, of order $l$ w.r.t. $z$ and of order $k$ w.r.t. $\tilde{z}$, computed at $(z,\tilde{z})$.\par  

Also, a subscript $b$ can be added - to obtain $C_b^{l,k}(Z\times \tilde{Z})$ - in order to specify that the functions and their succesive differentials are uniformly bounded. 

\item[$\bullet$] $\disp C_p(I)$: the set of piecewise continuous real valued functions defined on $I=[0,1]$. It is equipped with the supremum norm.\par 


\item[$\bullet$] $\disp C\big(\R_+,Z\big)$: the set of continuous processes defined on $\R_+$ with values in $Z$.

\item[$\bullet$] $\disp D\big(\R_+,Z\big)$: the set of right-continuous, left-limited (or \textit{c\`adl\`ag}\footnote{From French \textit{continu \`a droite et admettant une limite \`a gauche}.}) processes defined on $\R_+$ and taking values in $Z$. It is endowed with the Skorohod topologie. 
\end{description}

\section{Modeling and Asymptotics}

\subsection{One scale spatial models}

Since we are not concerned with the homogeneous models, we only review the mathematical details of the local models as given by Arnold and Theodosopulu in \cite{Arnold1980}.

\subsubsection{Deterministic Spatial Model}\label{model_M2}

One scale chemical reactions with diffusion in a spatial domain $I\subset\mathds{R}^d$ are modeled deterministically, at a macroscopic level, by the reaction-diffusion equation
\begin{equation}{\label{Eqn_Reaction_Diffusion}}
\displaystyle \frac{\partial v}{\partial t}=D\Delta v+R(v),
\end{equation}
\noindent where $v=v(t,x)$ is an $M$-dimensional real vector which gives the concentrations for the $M\geq 1$ reactants involved in the interactions, $x\in I$ is the spatial coordinate, $t\geq 0$ is the time. The domain $I$ can be bounded or unbounded, $v$ is subjected to boundary and initial conditions, $\Delta$ is the Laplace operator, $D$ is a diagonal matrix of size $M$, with non-negative coefficients on the diagonal, called the diffusion matrix. As it is often the case, we consider for simplicity $D=I_d$, the idendity matrix (in the sequel, $D$ is the differential operator). Reactions are represented by $R$, a polynomial vector field in $\R^M$, defined by \[ R(y)=(R_1(y),\cdots,R_M(y))\hspace{0.5cm}\text{with}\hspace*{0.5cm} R_i(y):=\sum_{|\alpha|\leq n}a_\alpha^i y^\alpha, \] where $n\in\mathds{N}^*$, $\alpha=(\alpha_1,\cdots,\alpha_M)\in\mathds{N}^M$ is a multi-index, $|\alpha|:=\alpha_1+\cdots+\alpha_M$, $a_\alpha^i\in\R$, $y=(y_1,\cdots,y_M)\in\mathds{R}^M$ and $\displaystyle y^\alpha:=y_1^{\alpha_1}\cdots y_M^{\alpha_M}$. \\

\noindent \textbf{\textit{Well-posedness, generator and debit function.}} For simplicity, we work on a one dimensional domain, $d=1$, that we project on the unit interval $I=[0,1]$ with periodic boundary conditions: $v(t,0)=v(t,1),$ $\forall t\geq0$. Since concentrations are positive quantities, a positive initial data is considered: $\forall x\in I$, $v(0,x)\geq0$ in the sense $v_i(0,x)\geq 0$, $\forall 1\leq i\leq M$. In order to have some control on the concentrations, it is assumed that $v(0)=v(0,\cdot)$ is bounded: $v(0,x)<\rho_1<\infty$, $\forall x\in I$, for some $\rho_1>0$, and that:
\begin{assumption}\label{well-posedness_assumpt_on_the_deb_funct_of_react_for_the_det_loc_model}
\hspace{1cm}\par  
(i) For all $y=(y_1,\cdots,y_M)\in\R^M$, $R_i(y)\geq 0$ when $y_i=0$ for some $1\leq i\leq M$.\par 
(ii) There exists $\rho_2>0$ such that $\langle R(y),y\rangle<0$ for all $y\in\R^M$ satisfying $|y|>\rho_2$. 
\end{assumption}
\noindent We have denoted by $|\cdot|$ and $\langle\cdot,\cdot\rangle$ the vector norm and the inner product of $\R^M$ respectively. \Cref{well-posedness_assumpt_on_the_deb_funct_of_react_for_the_det_loc_model} is natural, and ensures the consistency of the model through \textit{a priori} estimates for \eqref{Eqn_Reaction_Diffusion}. In fact, (i) yieds positiveness: $v_i(t,x)\geq0$, $\forall t,x$, $\forall i$, whereas (ii) yields boundedness: $v(t,x)\leq \rho$, $\forall t>0,\forall x$ with 
$\rho:=\max(\rho_1,\rho_2)+1$, thanks to the maximum principle.

As species all have the same dynamic in one scale models, a unique reactant is often considered, i.e. $M=1$. The function $R$ is written in the form $R(y)=b(y)-d(y)$, where $b$ and $d$ refer to variations due to births and deaths respectively. They are real valued polynomials with non-negative coefficients, such that $d(0)=0$ and $a_n<0$. We consider this case in this section.

We now consider an initial condition $v(0)=v_0\in C^3(I)$. In \cite{Kotelenez1986bis}, Kotelenez showed that, under the preceding conditions, the Cauchy problem associated with $(\ref{Eqn_Reaction_Diffusion})$ has a unique mild solution $v\in C\left(\mathds{R}_+;C^3(I)\right)$ satisfying \[ v(t)=T(t)v_0+\int_0^tT(t-s)R(v(s))ds, \] and $0\leq v(t) <\rho$ for all $t\geq0$. Here, $T(t):=\text{e}^{\Delta t}$ is the semigroup associated with the Laplace operator $\Delta$. We notice that $v_0\in C^3(I)$ is not necessary if the purpose is to solve the equation. Well-posedness still holds for much more general $v_0$. For more details, see \cite{Cazenave1998}, Section 5, or \cite{Arnold1980}, Lemma 2 (based on \cite{Kuiper1977}).\par 

Informally, the (weak) infinitesimal generator of the PDE $(\ref{Eqn_Reaction_Diffusion})$ reads
\begin{equation}\label{Eq_Reaction_Diffusion_Generator}
\displaystyle \mathcal{A}\varphi(u)=\big(D_u\varphi(u),\Delta u+R(u)\big)=\left\langle D_u\varphi(u),\Delta u+R(u)\right\rangle_2,
\end{equation}
\noindent for test functions $\varphi\in C_b^1\left(C^1(I)\right)$. The bracket $\langle\cdot,\cdot\rangle_2$ is the inner product of $L^2:=L^2(I)$, the space of square integrable real valued functions on $I$, and $(\cdot,\cdot)$ denotes the duality pairing. We give some probabilistic heuristics for understanding $(\ref{Eq_Reaction_Diffusion_Generator})$. Let $X(\cdot)$ be a solution to $(\ref{Eqn_Reaction_Diffusion})$, and denote by $\Phi(\cdot,u)$ the solution starting at $u\in C^3(I)$. Viewing $\disp \big(X(t)\big)_{t\geq 0}$ as a $C^3(I)$-valued deterministic - thus Markov - process on some filtered probability space $\disp \big(\Omega,\cF,(\cF_t)_{t\geq0},\P\big)$, its semigroup $\disp (P_t)_{t\geq0}$ on $\disp \fB\big(C^3(I)\big)$ reads \[ P_t\varphi(u):=\E\left[\varphi\big(X(t)\big)\big|X(0)=u\right]=\varphi\big(\Phi(t,u)\big). \] Then as usual, we interpret the generator $\cA$ as the derivative of the semigroup w.r.t. $t$, at $t=0$. Thus, for every test function $\varphi\in C_b^1\left(C^1(I)\right)$, the chain rule leads to \[ \mathcal{A}\varphi(u):=\left.\frac{\partial}{\partial t}\varphi\left(\Phi(t,u)\right)\right|_{t=0}=\left(D\varphi(u),\Delta u+R(u)\right). \] Finally, $D\varphi(u)$ being the differential of $\varphi$ at $u$, it belongs to $\left(C^1(I)\right)'$. We identify it with the (abstract) gradient of $\varphi$ at $u$ since $C^1(I)\subset L^2$, and write \[ \left(D\varphi(u),\Delta u+R(u)\right)=\left\langle D\varphi(u),\Delta u+R(u)\right\rangle_2, \] thanks to the Riesz representation theorem. Relation $(\ref{Eq_Reaction_Diffusion_Generator})$ then follows. $\square$\par 

After all, we define the corresponding debit function as \[ \psi(u)=\Delta u+R(u). \] It is the vector field, in $C^1(I)$ here, determining the flow of the equation. More precisely, $u\mapsto\Delta u$ (resp. $u\mapsto R(u)$) is the debit function related to diffusions (resp. reactions).

\subsubsection{Stochastic Spatial Model}

Following \cite{Arnold1980}, we devide the unit interval into $N$ smaller intervals of equal length $\displaystyle N^{-1}$: $\displaystyle I_j=\big((j-1)N^{-1},jN^{-1}\big]$, for $j=1,\cdots,N$, called \tt \textbf{sites}\normalfont. Molecules are produced - birth - or removed - death - on each site according to chemical reactions. They also diffuse between sites by simple random walks, at rates proportional to $N^2$. This couples the site reactants and extends Kurtz's model (\hyperlink{M3}{M3}) to the spatially inhomogeneous case (\hyperlink{M4}{M4}). Define: \vspace{0.1cm}\par 
 
$\bullet$ $\displaystyle X_{j}^N$: the number of molecules on site $j$, \hspace{0.2cm} $j=1,\cdots,N$, \vspace{0.1cm}\par 
$\bullet$ $\displaystyle X^N=\left(X_j^N\right)_{1\leq j\leq N}\in\mathds{N}^{N}$: the molecular composition of the whole system,\vspace{0.1cm}\par 
$\bullet$ $\mathfrak{R}$: the set of possible onsite reactions (we assume that $\mathfrak{R}$ is finite).\vspace{0.1cm}\par 

\noindent Consider the natural completed filtration $\displaystyle \left\{\mathcal{F}_t^{N}, t\geq 0\right\}$, where $\displaystyle \mathcal{F}_t^{N}$ is the completion of the $\sigma$-algebra $\displaystyle \sigma\left(\left\{X^N(s):s\leq t\right\}\right)$ with null probability measure sets. That will be the filtration considered by default, for each process we will define. From \cite{Kurtz1970,Kurtz1971} or Chapter4, Section2, p.162-164 of \cite{Kurtz1986}, we know that $\displaystyle \left\{X^N(t),t\geq0\right\}$ is a time homogeneous $\mathds{N}^{N}-$valued jump Markov process, with the following transitions on any $j$:\vspace{0.1cm}\par 

\noindent \underline{Onsite reactions.} A chemical reaction $r\in\mathfrak{R}$ occurs on the site $j$: \normalfont \[ X_{j}^N\longrightarrow X_{j}^N+\gamma_{j,r},\hspace{0.3cm} \gamma_{j,r}\in\mathds{Z},\hspace{0.3cm} \text{at rate}\hspace{0.3cm} \lambda_r\left(X_j^N\right). \]
\noindent \underline{Diffusions.} A molecule moves from the site $j$ to the site $j-1$ or  to the site $j+1$: \normalfont \[ \left\{
\begin{array}{l}
\displaystyle \left(X_{j-1}^N,X_{j}^N\right)\longrightarrow \left(X_{j-1}^N+1,X_{j}^N-1\right),\hspace{0.3cm} \text{at rate}\hspace{0.2cm} N^2X_{j}^N\vspace{0.1cm},\\

\displaystyle \left(X_{j}^N,X_{j+1}^N\right)\longrightarrow \left(X_{j}^N-1,X_{j+1}^N+1\right),\hspace{0.3cm} \text{at rate}\hspace{0.2cm} N^2X_{j}^N.
\end{array}
\right. \]
%
At the level of the whole vector, we get for $1\le j \le N$
\[
\left\{
\begin{array}{l}
\displaystyle X^N\longrightarrow X^N+\gamma_{j,r}e_j\hspace{0.3cm} \text{at rate}\hspace{0.2cm} \lambda_r\left(X_j^N\right),\vspace{0.1cm}\\
\displaystyle X^N\longrightarrow X^N+e_{j-1}-e_j\hspace{0.3cm} \text{at rate}\hspace{0.2cm} N^2X_{j}^N\vspace{0.1cm}\\
\displaystyle X^N\longrightarrow X^N+e_{j+1}-e_j\hspace{0.3cm} \text{at rate}\hspace{0.2cm} N^2X_{j}^N,\vspace{0.1cm}\\
\end{array}
\right.
\]

\noindent where $\left\lbrace e_j, j=1,\cdots,N\right\rbrace$ is the canonical basis of $\mathds{R}^N$. 

In all the article, we assume that the rates of onsite reactions are polynomial functions and that there are a finite number of reaction in $\mathfrak{R}$.\\

Two types of events are clearly distinguishable: onsite or chemical reactions, and diffusions. We sometimes say reaction for each of them, but we avoid any possibility of confusion. The vectors $\gamma_{j,r}e_j$, $e_{j-1}-e_j$, $\cdots$ appearing in the transitions are the stocheometric coefficients of the corresponding event. They point out the different possible directions that the jump due to an event can take. We also denote them by jump height or simply jump. 

We insist on the fact that the directions of jumps are actually characterized by two parameters: a "true" event - reaction or diffusion - and a location - site - where that event occurs. Furthermore, they are independently chosen and the total jump rate, denoted by $\lambda^N$, is the sum of the rates over all the possible directions of jump:
\[ \lambda^N\left(X^N\right)=\sum_{j=1}^N\left[\sum_{r\in\mathfrak{R}}\lambda_r\left(X_j^N\right)+ N^2\left(X_{j-1}^N+2X_j^N+X_{j+1}^N\right)\right]. \]

\noindent In general, the jump height of an onsite reaction may depend both on that reaction and on the site where it occurs. In the present model, onsite reactions are assumed to be spatially homogeneous, that is \[ \gamma_{j,r}=\gamma_r\hspace{0.2cm} \text{ for all }\hspace{0.2cm} 1\leq j\leq N. \]
%

\noindent \textbf{\textit{Well-posedness, generator and debit function.}}\label{debit_function_as_generator_derivative_1} Such a Markov process can be constructed, relying on a sequence $(\tau_i)_{i\geq1}$ of independent random waiting times with exponential distributions. We set $T_0=0$ and $T_k=\tau_1+\cdots+\tau_k$ for $k\geq1$. $\displaystyle \left\lbrace T_k, k\geq 1\right\rbrace$ is a family of $\displaystyle \left\{\mathcal{F}_t^{N}\right\}$-stopping times such that $X^N$ is constant on each $[T_{k-1},T_k)$, and has a jump at $T_k$. The parameter of $\tau_k$ is the total jump rate $\disp \lambda^N\left(X^N(T_{k-1})\right)$, at time $T_{k-1}$. For $t>0$, \[ \mathds{P}\left\lbrace \tau_k>t\left|T_{k-1}\right.\right\rbrace=\text{exp}\left\lbrace-t\lambda^N\left(X^N(T_{k-1})\right)\right\rbrace. \] At time $T_k$, a given event (or direction of jump) is choosen with probability \[ \frac{"\text{rate of the direction of jump at time }T_{k-1}"}{"\text{total rate of jump at time }T_{k-1}"}. \]

Periodic boundary conditions are considered on the spatial domain, and $X_j^N$ is viewed as a sequence satisfying $X_{j+N}^N=X_j^N$ for all $j$ (as in \cite{Blount1987}, see p.10). The infinitesimal generator of our jump Markov process then reads
\begin{equation}\label{Infinitesimal_generator_zero}
\begin{array}{l}
\displaystyle \mathcal{A}^N\varphi(X):=\left.\frac{d}{dt}\mathds{E}_X\left[\varphi\left(X^N(t)\right)\right]\right|_{t=0}\hspace{0.2cm}:=\hspace{0.2cm}\lim_{t\rightarrow0}\frac{1}{t}\mathds{E}_X
\left[\varphi\left(X^N(t)\right)-\varphi(X)\right]\\


\displaystyle \hspace{1.7cm} =\sum_{j=1}^N\sum_{r\in\mathfrak{R}}\left[\varphi
\left(X+\gamma_re_j\right)-\varphi(X)\right]\lambda_r\left(X_j\right)\\

\displaystyle \hspace{2cm}+\sum_{j=1}^N \left[\varphi\left(X+e_{j-1}-e_j\right)+\varphi\left(X+e_{j+1}-e_j\right)-2\varphi(X)\right] N^2X_j.
\end{array}
\end{equation}

We need not specify the complete domain of the generator $\mathcal{A}^N$. The limit above holds, for test functions $\displaystyle \varphi\in C_b\left(\mathds{R}^N\right)$ (see \cite{Kurtz1986} p.162-164).\\

In this stochastic description, the global debit is the function which sums the jump heights weighted with the corresponding rates, under the possible directions of jump. It is a vector field in the state space. If we denote it by $\Psi^N$, then  
\begin{equation}\label{debit_of_the_one_scale_stochastic_spatial_model}
\begin{array}{l}
\disp\hspace{0.2cm} \Psi^N(X):=\sum_{j=1}^N\sum_{r\in\mathfrak{R}}\gamma_r\lambda_r(X_j)e_j+ N^2(e_{j-1}-2e_j+e_{j+1})X_j\\

\disp\hspace{1.7cm} =\sum_{j=1}^N\left(\sum_{r\in\mathfrak{R}}\gamma_r\lambda_r(X_j)
+N^2\left(X_{j-1}-2X_j+X_{j+1}\right)\right)e_j
\end{array}
\end{equation}
for $X\in\mathds{N}^N$, thanks to periodicity and a change of subscript. We clearly identify the debit of onsite reactions, and that of diffusions. The latter, the linear part in $X$, defines a discrete Laplace operator that we introduce later.\\

\noindent\textbf{\textit{Scaling and density dependence:}} The following are parts of the main assumptions in the framework of one scale models.
\begin{itemize}
\item \noindent On each site, the reactant has a population size of order $\mu$. Typically, $\mu$ is large and represents the initial average number of particles on each site.

\item Density dependence holds for the reaction rates. That is: for all $r\in\fR$, there exists $\tilde{\lambda}_r$ satisfying \[ \lambda_r\left(X_j^{N}\right)=\mu\tilde{\lambda}_r\left(\frac{X_j^{N}}{\mu}\right)\hspace{0.4cm} \text{for all}\hspace{0.3cm} j=1,\cdots,N. \]
\end{itemize}

\noindent Density dependence is a natural assumption when dealing with chemical reaction systems (see \cite{Kurtz1970}, or Chapter11 of \cite{Kurtz1986}). It clearly holds for linear rates for instance, and can be understood as the fact that reducing the scale is somehow increasing the density, and results in speeding up the dynamic. 

In general $\mu$ depends on $N$ and we omit to mention the dependance on $\mu$.  

We rescale our process, setting \[ U_j^{N}:=U_j^{N,\mu}:=\frac{X_j^{N}}{\mu},\hspace{0.6cm} \text{and}\hspace{0.4cm}  U^N:=\left(U_j^{N}\right)_{1\leq j\leq N}, \] and do not distinguish $\lambda_r$ from $\tilde{\lambda}_r$ in the sequel, for all $r\in\fR$. The new scaled process $U^N$ has values in $\displaystyle \mathds{R}^N$ and gives the proportions of the reactant on each site. With some abuse, we often say concentration instead of proportion.\\

In order to achieve a pointwise modeling over the whole spatial domain, a space-time jump Markov process is constructed, which is the stochastic counterpart of the solution to $(\ref{Eqn_Reaction_Diffusion})$. As in \cite{Arnold1980}, \cite{Kotelenez1986} or \cite{Blount1987,Blount1992}, we introduce the step function 
\begin{equation}\label{Modeling_space}
u^N(t,x)=\sum_{j=1}^NU_j^N(t)\mathds{1}_j(x), \hspace{0.3cm} t\geq 0,\hspace{0.2cm} x\in I_j,
\end{equation}
\noindent where $\displaystyle \mathds{1}_j(\cdot):=\mathds{1}_{I_j}(\cdot)$ is the indicator function of the $j-$th site $\displaystyle I_j$. Note that for all $t\geq 0$, $N\geq1$, the function $u^N(t,\cdot)$ can be identified with the vector $U^N(t)$ of $\R^N$. Furthermore, since we consider a periodic framework we have $U_{j+N}^N(t)=U_j^N$ and $u^N(t,\cdot)$ is $1-$periodic.
Below, we use the standard identification $\displaystyle u^N(t):=u^N(t,\cdot)$.

Let $\displaystyle \mathds{H}^N$ denote the subspace of $\displaystyle L^2$ which consists of real-valued and 1-periodic step functions that are constant on the intervals $I_j$, $1\leq j\leq N$. Therefore, $\displaystyle u^N:=\left\lbrace u^N(t),t\geq 0\right\rbrace$ is a c\`adl\`ag (time-space) Markov process, with values in $\displaystyle \mathds{H}^N$ and jumps given by: 
\[ \left\{
\begin{array}{l}
\displaystyle u^{N}\longrightarrow u^{N}+\frac{\gamma_r}{\mu}\mathds{1}_j,\hspace{0.3cm} \text{at rate}\hspace{0.2cm} \mu\lambda_r\left(u_j^{N}\right),\hspace{0.2cm}\text{for}\hspace{0.2cm} r\in\mathfrak{R},\vspace{0.1cm}\\

\displaystyle u^N\longrightarrow u^N+\frac{\mathds{1}_{j-1}-\mathds{1}_j}{\mu},\hspace{0.3cm} \text{at rate} \hspace{0.2cm} \mu N^2u_j^{N},\hspace{0.2cm}\text{for}\hspace{0.2cm} 1\le j\le N,\vspace{0.1cm}\\

\displaystyle u^N\longrightarrow u^N+\frac{\mathds{1}_{j+1}-\mathds{1}_j}{\mu},\hspace{0.3cm} \text{at rate} \hspace{0.2cm} \mu N^2u_j^{N}, \hspace{0.2cm}\text{for}\hspace{0.2cm} 1\le j\le N.
\end{array}
\right. \]
Here, $u_j^N=U_j^N$ is the $j-$th "coordinate" of $u^N$, obtained by the canonical projection 
\begin{equation}\label{Canonical_Projection}
\begin{tabular}{lll}
$\displaystyle P_N$ & $:$ & $L^2 \longrightarrow \mathds{H}^N$\vspace{0cm}\\
 & & $\displaystyle u\longmapsto P_Nu=\sum_{j=1}^Nu_j\mathds{1}_j,\hspace{0.4cm}\text{where}\hspace{0.4cm} u_j:=N\int_{I_j}u(x)dx$.
\end{tabular}
\end{equation}

\noindent The infinitesimal generator $\mathcal{A}^{N}:=\displaystyle \mathcal{A}^{N,\mu}$ of $u^N$ is given by 
\begin{equation}\label{Stochastic_Generator}
\begin{array}{l}
\displaystyle \mathcal{A}^{N}\varphi(u)=\sum_{j=1}^N\sum_{r\in\mathfrak{R}}\left[\varphi\left(u+
\frac{\gamma_r}{\mu}\mathds{1}_j\right)-\varphi(u)\right]\mu\lambda_r\left(u_j\right)\\

\displaystyle \hspace{1.8cm}+\sum_{j=1}^N\left[\varphi\left(u+\frac{\mathds{1}_{j-1}-\mathds{1}_j}{\mu}\right)+\varphi\left(u+\frac{\mathds{1}_{j+1}-\mathds{1}_j}{\mu}\right)-2\varphi(u)\right]\mu N^2u_j,\\
%
%
\end{array}
\end{equation} 
\noindent on the domain $C_b\left(\mathds{H}^N\right)$. It can be extended to $C_b\left(L^2\right)$ by \[ \bar{\mathcal{A}}^{N}\varphi(u):=\mathcal{A}^{N}\varphi(P_Nu). \] Such a Markov process does exist and is unique (see \cite{Kotelenez1988bis} whose result is based on \cite{Kurtz1986}) until a possible blow-up time. In addition, under natural assumptions on the reaction rates, we have $u^N(t)\geq0$ for all $t$, as soon as $u^N(0)\geq0$. 

As already mentioned, the convergence of $u^N$ to $u$ solution of \eqref{Eqn_Reaction_Diffusion} has been the object of several articles. This implies in particular that for $N$ large enough there is no blow-up.\\

For $f,g\in\mathds{H}^N$, the $L^2$ inner product reads $\left\langle f,g\right\rangle_2=N^{-1}\sum_{j=1}^Nf_jg_j$. We denote by $\Vert\cdot\Vert_2$ the $L^2$ norm. Also, the supremum norm is given by $\Vert f\Vert_\infty=\sup_{1\leq j\leq N}|f_j|$. Setting $\mathds{H}:=\cup_{N\geq1}\mathds{H}^N$, the following obviously holds.

%
%

\begin{proposition}\label{Piecewise_Approximation_via_Projection} \hspace{1cm}\vspace{0.1cm}\par 
\textbf{(i)} $\displaystyle \left(\mathds{H}^N,\langle\cdot,\cdot\rangle_2\right)$ is a finite dimensional Hilbert space with $\displaystyle \left\{\sqrt{N}\mathds{1}_j, 1\leq j\leq N\right\}$ as an orthonormal basis.\par 
\textbf{(ii)} $\displaystyle \lim_{N\rightarrow\infty}\Vert P_Nu-u\Vert_2\longrightarrow 0$ for $u\in L^2$, and $\left(\mathds{H},\Vert\cdot\Vert_2\right)$ is dense in $\big(L^2,\Vert\cdot\Vert_2\big)$. \par 
\textbf{(iii)} $P_N$ is a contracting linear continuous operator on $\big(C_p(I),\Vert\cdot\Vert_\infty\big)$ such that $\displaystyle \lim_{N\rightarrow\infty}\Vert P_Nu-u\Vert_\infty\longrightarrow 0$ for $u\in C_p(I)$, and $\left(\mathds{H},\Vert\cdot\Vert_\infty\right)$ is dense in $\big(C_p(I),\Vert\cdot\Vert_\infty\big)$.
\end{proposition}

\subsection{Multiscale Stochastic Spatial Model}

Now we wish to generalize \cite{Arnaud2012} to a spatially dependent context. We consider two orders of population sizes and repeat the previous procedure. Note that at least two types of reactants are needed, i.e. $M\geq2$ is necessary. For simplicity, we consider exactly two reactants, i.e. we take $M=2$, that we denote by $C$ and $D$. The former is abundant and the latter is rare. Henceforth, the super/subscript $C$ (resp. $D$) refers to the species $C$ (resp. $D$). 

The set $\mathfrak{R}$ of possible onsite reactions is divided in three disjoint subsets \[ \mathfrak{R}=\mathfrak{R}_C\cup\mathfrak{R}_{DC}\cup\mathfrak{R}_D. \] Reactions in $\fR_C$ (resp. $\fR_D$) involve only molecules of $C$ (resp. $D$) as reactants and products, whereas, reactions in $\mathfrak{R}_{DC}$ involve both types of reactants and/or products. 
We assume:

\begin{assumption}{\label{density_dependence}} \hspace{1cm}\vspace{0.1cm}\par 
\textbf{(i)} Density dependence holds for reactions in $\fR_C$.\par 
\textbf{(ii)} Reactions $r\in\fR_C$ are spatially homogeneous and fast ( "pure fast reactions"), while reactions $r\in\fR_D$ are slow ("pure slow reactions").\par 
\textbf{(iii)} The molecules of $C$ diffuse, while those of $D$ do not. 
\end{assumption}

\noindent Now, we define\vspace{0.1cm}\par 

$\bullet$ $\displaystyle \hspace{0cm} X_j^{N,C}$\hspace{0.1cm} $\displaystyle \left(\hspace{0cm} \text{resp.} \hspace{0.1cm}  X_j^{N,D}\hspace{0.1cm}\right)$ as the number of molecules of $C$ (resp. $D$) on the site $j$, \vspace{0.1cm}\par 

$\bullet$ $\displaystyle \hspace{0cm} X_C^N:=\left(X_j^{N,C}\right)_{1\leq j\leq N}$\hspace{0cm},\hspace{0.2cm} $\displaystyle X_D^N:=\left(X_j^{N,D}\right)_{1\leq j\leq N}$\hspace{0.2cm} and \hspace{0.2cm}$\displaystyle X^N:=\left(X_C^N,X_D^N\right)\in\N^{2N}$.\vspace{0.1cm}\par 

\noindent Scales are set such that on every site, the initial average number of molecules for $C$ is of order $\mu$, while it is of order $\kappa$ for $D$, with $\mu\gg \kappa$. Namely, if the total initial population is $\disp \mathfrak{M}=\mathfrak{M}_C+\mathfrak{M}_D$, then $\displaystyle \mathfrak{M}_C\approx N\times\mu$ while $\disp \mathfrak{M}_D\approx N\times \kappa$. We set \\

$\displaystyle \hspace{3cm} U_j^{N,C}:=U_j^{N,\mu,C}:=\frac{X_j^{N,C}}{\mu},\hspace{1cm}U_j^{N,D}:=U_j^{N,\mu,D}:=X_j^{N,D}$,\par

$\displaystyle \hspace{4cm} U_C^N:=\left(U_j^{N,C}\right)_{1\leq j\leq N},\hspace{1cm}U_D^N:=\left(U_j^{N,D}\right)_{1\leq j\leq N}$,\par 

\noindent and finally \par 

$\displaystyle\hspace{5cm} U^{N,\mu} =:U^N:= \left(U_C^{N},U_D^{N}\right)$. \\

\noindent The new scaled process $U^N$ has values in $\displaystyle \mathds{R}^N\times\mathds{N}^N$. Its generator has the form \\

\noindent $\displaystyle \mathcal{A}^{N,\mu}\varphi\left(U_C,U_D\right)$\par \vspace{0.1cm}

$\displaystyle \hspace{0cm}=\sum_{j=1}^N\left\lbrace\sum_{r\in\mathfrak{R}_C}\left[\varphi\left(U_C+
\frac{\gamma_{j,r}^C}{\mu}e_j,U_D\right)-\varphi(U_C,U_D)\right]\mu\lambda_r\left(U_j^C\right)\right.$\vspace{0.1cm}\par 

$\displaystyle \hspace{1.5cm}\left.+\sum_{r\in\mathfrak{R}_{DC}}\left[\varphi\left(U_C+\frac{\gamma_{j,r}^C}{\mu}e_j,U_D+
\gamma_{j,r}^{D}e_j\right)-\varphi(U_C,U_D)\right]\lambda_r\left(U_j^C,U_j^D\right)\right.$\par 

$\displaystyle \hspace{1.5cm}\left.+\sum_{r\in\mathfrak{R}_D}\left[\varphi\left(U_C,U_D+\gamma_{j,r}^{D}e_j\right)-\varphi(U_C,U_D)\right]\lambda_r\left(U_j^D\right)\right\rbrace$\par 

$\displaystyle \hspace{0.1cm}+\sum_{j=1}^N\left\{\left[\varphi\left(U_C+\frac{e_{j-1}-e_j}{\mu},U_D\right)+\varphi\left(U_C+\frac{e_{j+1}-e_j}{\mu},U_D\right)-2\varphi(U_C,U_D)\right] \mu N^2U_j^C \right.$\par 



\noindent for test functions $\displaystyle \varphi\in C_b\left(\mathds{R}^{2N}\right)$ and $U=(U_C,U_D)\in\mathds{R}^{2N}$. 

Below, we discuss the nature of the reactions in $\fR_{DC}$, and consider situations where $\gamma_j^{r,D}$ and $\lambda_r$ need to be rescaled, for $r\in\fR_D$ and some $r\in\fR_{DC}$.\\

Next, we introduce the step function 
\begin{equation}\label{time_space_function}
u^N(t,x)=\sum_{j=1}^NU_j^N(t)\mathds{1}_j(x), \hspace{0.3cm} t\geq 0,\hspace{0.2cm} x\in I_j,
\end{equation} 
\noindent which belongs to $\displaystyle \mathds{H}^N\times\mathds{H}^N$. It involves the following functions for each component \vspace{0.2cm}\par 

$\disp \hspace{1.3cm} u_C^N(t,x)=\sum_{j=1}^Nu_j^{N,C}(t)\mathds{1}_j(x)\hspace{0.5cm}\text{and}\hspace{0.4cm}u_D^N(t,x)=\sum_{j=1}^Nu_j^{N,D}(t)\mathds{1}_j(x)$,\par  

\noindent where\par 

$\disp\hspace{2cm}u_j^{N,C}(t):=P_Nu_C^N(t)=N\int_{I_j}u_C^N(t,x)dx = U_j^{N,C}(t)$,\\

\noindent and a similar relation holds for the reactant $D$. Therefore, $\displaystyle u^N:=\left\{u^N(t),t\geq 0\right\}$ is a $\displaystyle \mathds{H}^N\times\mathds{H}^N$-valued c\`adl\`ag jump Markov process, with the transitions:
\begin{equation}\label{transitions_for_the_spatial_model}
\left\{
\begin{tabular}{l}
$\displaystyle \left(u_C^N,u_D^N\right)\longrightarrow\left(u_C^N+\frac{\gamma_{j,r}^C}{\mu}\mathds{1}_j,u_D^N\right)$,\hspace{0cm} at rate $\displaystyle \mu\lambda_r\left(u_j^{N,C}\right)$, \hspace{0cm} for $r\in\mathfrak{R}_C$,\vspace{0.1cm}\\

$\displaystyle \left(u_C^N,u_D^N\right)\longrightarrow\left(u_C^N+\frac{\gamma_{j,r}^C}{\mu}\mathds{1}_j,u_D^N+\gamma_{j,r}^{D}\mathds{1}_j\right)$,\hspace{0cm} at  rate $\displaystyle \lambda_r\left(u_j^{N,C},u_j^{N,D}\right)$, $r\in\mathfrak{R}_{DC}$,\vspace{0.1cm}\\

$\displaystyle \left(u_C^N,u_D^N\right)\longrightarrow\left(u_C^N,u_D^N+\gamma_{j,r}^{D}\mathds{1}_j\right)$,\hspace{0cm} at rate $\displaystyle \lambda_r\left(u_j^{N,D}\right)$,\hspace{0cm} for $r\in\mathfrak{R}_{D}$,\vspace{0.1cm}\\

$\displaystyle \left(u_C^N,u_D^N\right)\longrightarrow \left(u_C^N+\frac{\mathds{1}_{j-1}-\mathds{1}_{j}}{\mu},u_D^N\right)$,\hspace{0cm} at rate $\displaystyle N^2u_j^{N,C}$, for $1\le j\le N$,\vspace{0.1cm}\\

$\displaystyle \left(u_C^N,u_D^N\right)\longrightarrow\left(u_C^N+\frac{\mathds{1}_{j+1}-\mathds{1}_{j}}{\mu},u_D^N\right)$,\hspace{0cm} at rate $\displaystyle N^2u_j^{N,C}$, for $1\le j\le N$,\\

%
\end{tabular} \hspace{0cm}
\right.
\end{equation}
We endow $\displaystyle \mathds{H}^N\times\mathds{H}^N$ with the norm \[ \Vert (f_1,f_2)\Vert_{\infty,\infty}:=\Vert f_1\Vert_\infty+\Vert f_2\Vert_\infty,\hspace{0.4cm}\text{for}\hspace{0.3cm}f_1,f_2\in\mathds{H}^N. \] 
\indent As before, we assume that onsite reactions in $\mathfrak{R}_C$ are spatially homogeneous, that is \[ \gamma_{j,r}^C=\gamma_r^C\hspace{0.2cm} \text{ for all }\hspace{0.2cm} 1\leq j\leq N,\; r\in \mathfrak{R}_C. \]

Let us specify the description of the mixed reactions.
\begin{assumption}\label{Partitioning_Mixed_Reactions_subset}
In some $\displaystyle S_1\subset\mathfrak{R}_{DC}$, reactions are spatially homogeneous, fast and do not affect the discrete species:
\begin{itemize}
%
%
\item $\gamma_{j,r}^D=0$ and the rate is $\mu \tilde\lambda_r\left(u_j^{N,C},u_j^{N,D}\right)$ for $
r\in S_1\subset\mathfrak{R}_{DC}$

\end{itemize}
\end{assumption}

Again we omit the tilde below.

The limit $N\to \infty$ of the above system creates mathematical difficulties (see \cref{r3.3} and \cref{r4.1} below). We introduce some 
spatial correlations for reactions in $\disp \left(\mathfrak{R}_{DC}\backslash S_1\right)\cup\mathfrak{R}_D$.

All reactions in $\disp \left(\mathfrak{R}_{DC}\backslash S_1\right)\cup\mathfrak{R}_D$ are slow. We assume that when such a reaction occurs on a site, it affects the neighboring sites. A natural way to do this is to assume that when the reaction $r\in \disp \left(\mathfrak{R}_{DC}\backslash S_1\right)\cup\mathfrak{R}_D$ occurs at  a site $j\in \{1,\dots,N\}$, each site $i\in \{1,\dots,N\}$, is affected in the following way:
\begin{equation}\label{e2.10}
(u^{N,C}_j,u^{N,D}_j)\longrightarrow (u^{N,C}_j+\frac{\gamma_r^C}{\mu} \gamma_{ij}^N,u^{N,D}_j+\gamma_r^D \gamma_{ij}^N).
\end{equation}
A first possibility is to take $\gamma_{ij}^N=\frac1A$ for $|i-j|< A$, describing the case when a reaction at site $i$ triggers simultaneously the same reaction at sites at distance less than
$\frac{A}{N}$, while keeping the global effect of order $1$. 
The above model corresponds to $A=1$ but below we assume that $A$ has the same order as $N$.  
More generally,  
we choose a $1$-periodic $C_p(I)$ function $a$ and choose:
$$
\gamma_{ij}^N=\int_{I_i}a\left(x-\frac{j}{N}\right)dx.
$$
It is natural to assume that $a$ is maximum at zero, and is even decaying with respect to $|x|$. The preceding case is recovered with $a$ being an indicator function. 
Note that we should require that $\gamma_r^D \gamma_{ij}^N$ is an integer, but it is complicated to write general conditions to ensure this. 

This modeling has a problem. Indeed, the positivity of the proportion of molecules is ensured by the fact that the rate vanishes when the proportion vanishes, but now, since we have a nonlocal effect this is not true. It would not be realistic to consider a rate which depends on the proportion at all neighboring sites. 
Instead, we simply assume that when a reaction occurs at site $j$, it has an effect at a site $i$ only if the reaction is possible.
Thus we would like to replace \eqref{e2.10} by 
\begin{equation}\label{e2.11}
(u^{N,C}_j,u^{N,D}_j)\longrightarrow \left(u^{N,C}_j+\frac{\gamma_r^C}{\mu} \gamma_{ij}^N\theta_{ij}^r,u^{N,D}_j+\gamma_r^D \gamma_{ij}^N\theta_{ij}^r\right)
\end{equation}
with $\theta_{ij}^r=\theta_{ij}^r\big(u_j^{N,C},u_j^{N,D}\big)=\1_{u^{N,C}_j+\frac{\gamma_r^C}{\mu} \gamma_{ij}^N\ge 0}\1_{u^{N,D}_j+\gamma_r^D \gamma_{ij}^N\ge 0}$, where $\gamma_r^C=0$ for $r\in\fR_D$.
Again, this creates mathematical difficulties, since the indicator function is not smooth. We replace it by a smooth function $\theta$ equal to $0$ on $(-\infty,0]$, to $1$ on $[1,\infty)$ and smoothly increasing on $[0,1]$. We are now able to write the multiscale model we consider.

\begin{assumption}\label{spatial_correlation_for_discrete_variable}
When a slow reaction $r\in\left(\mathfrak{R}_{DC}\backslash S_1\right)\cup\mathfrak{R}_D$ occurs on a site $j$, it affects $C$ and $D$ on a site $1\leq i\leq N$ in the following way: 
\[ \gamma_{j,r}^{C}\equiv\gamma_{j,r}^{N,C}\big(u_C^N,u_D^N\big):=\gamma_r^C\sum_{i=1}^N\gamma_{ij}^N\mathds{1}_i\theta_{ij}^r\big(u_i^{N,C},u_i^{N,D}\big) \]
and
\[ \gamma_{j,r}^{D}\equiv\gamma_{j,r}^{N,D}\big(u_C^N,u_D^N\big):=\gamma_r^{D}\sum_{i=1}^N\gamma_{ij}^N\mathds{1}_i\theta_{ij}^r\big(u_i^{N,C},u_i^{N,D}\big), \] 
where $\disp \gamma_{ij}^N:=\int_{I_i}a\left(x-\frac{j}{N}\right)dx$, and $\theta_{ij}^r\big(u_i^{N,C},u_i^{N,D}\big)=\theta\big(u^{N,C}_i+
\frac{\gamma_r^C}{\mu} \gamma_{ij}^N\big) \theta\big(u^{N,D}_i+\gamma_r^D \gamma_{ij}^N\big)$.
\end{assumption}

\begin{remark}\label{bound_of_jumps}
\Cref{spatial_correlation_for_discrete_variable} points out some spatial correlation for slow reactions. These latter start on a source site $j$, and then, their effect spreads inside the system. This is a natural assumption in molecular biologie. In the cell, the synthesis of some rare proteins is the launch for a sequence of chemical reactions in chain. This kind of phenomenon preceeds the activation of a gene on DNA (DeoxyriboNucleic Acid) for instance. 
\end{remark}


%
Since $\1_i(x)$ is nonzero only on $I_i$  and $\gamma_{ij}^N\le \frac{a(0)}{N} $, we clearly have:
\begin{equation}\label{maximum_jump_amplitude_for_slow_reactions}
\disp \left\Vert\gamma_{j,r}^{N,C}(u_C^N,u_D^N)\right\Vert_\infty\leq\frac{\big|\gamma_r^C| a(0)}{ N}\hspace{0.5cm}\text{and}\hspace{0.5cm} \disp \left\Vert\gamma_{j,r}^{N,D}(u_C^N,u_D^N)\right\Vert_\infty\leq\frac{\big|\gamma_r^D| a(0)}N
\end{equation}


The infinitesimal generator of $u^N$ on $\displaystyle C_b\left(\mathds{H}^N\times\mathds{H}^N\right)$ then reads
\small
\begin{equation}\label{Hybrid_Stochastic_Generator}
\begin{array}{l}
\displaystyle \mathcal{A}^{N}\varphi(u_C,u_D)\hspace{1cm}\\


\displaystyle \hspace{0.4cm}=\sum_{j=1}^N\left\lbrace\sum_{r\in\mathfrak{R}_C}\left[\varphi\left(u_C+
\frac{\gamma_r^C}{\mu}\mathds{1}_j,u_D\right)-\varphi(u_C,u_D)\right]\mu\lambda_r\left(u_j^C\right)\right.\\


\displaystyle \left.\hspace{1.8cm}+\sum_{r\in S_1}\left[\varphi\left(u_C+
\frac{\gamma_r^C}{\mu}\mathds{1}_j,u_D\right)-
\varphi(u_C,u_D)\right]\mu\lambda_r\left(u_j^C,u_j^D\right)\right.\\


\displaystyle \hspace{1.8cm}\left.+\sum_{r\in\mathfrak{R}_{DC}\backslash S_1}\left[\varphi\left(u_C+\frac{\gamma_{j,r}^{N,C}(u_C,u_D)}{\mu},u_D+\gamma_{j,r}^{N,D}(u_C,u_D)\right)-\varphi(u_C,u_D)\right]\right.\\
\disp \hspace{13cm}\left.\times\lambda_r\left(u_j^C,u_j^D\right)\right.\\


\displaystyle \hspace{1.8cm}\left.+\sum_{r\in\mathfrak{R}_D}\left[\varphi\left(u_C,u_D+
\gamma_{j,r}^{N,D}(u_C,u_D)\right)-\varphi(u_C,u_D)\right]\lambda_r\left(u_j^D\right)\right\rbrace\\


\displaystyle \hspace{0.5cm}+\sum_{j=1}^N\left[\varphi\left(u_C+\frac{\mathds{1}_{j-1}-\mathds{1}_j}{\mu},u_D\right)+\varphi\left(u_C+\frac{\mathds{1}_{j+1}-\mathds{1}_j}{\mu},u_D\right)-2\varphi(u_C,u_D)\right]\mu N^2u_j^C
\end{array}
\end{equation}
\normalsize

\noindent It can be extended to $\displaystyle C_b\left(L^2\times L^2\right)$ by \[ \mathcal{A}^{N}\varphi(u_C,u_D):=\mathcal{A}^{N,\mu}\varphi(P_Nu_C,P_Nu_D). \] 
For $u:=(u_C,u_D)\in L^2\times L^2$, we may often use the notation $\displaystyle \tilde{P}_Nu:=\left(P_Nu_C,P_Nu_D\right)$.

Again, such a process exists until a possible blow-up time. If the process blows up, we say that from this time it takes a cemetary value $\bar \Delta$ whose distance
to any point is $1$.

\subsection{Convergence of (M6)}\label{approximations}

We discuss the asymptotic behavior of the sequence $u^N$ of Markov processes representing (\hyperlink{M6}{M6}), as $N$ and $\mu$ go to the infinity, for polynomial (onsite) reaction rates.

\subsubsection{Identification of the limit}

\noindent \textbf{\textit{Some preliminaries: a discrete Laplace.}} For $f\in C_p(I)$ and for $x\in[0,1]$, we set \[ \nabla_N^+ f(x):=N\left[f\left(x+\frac{1}{N}\right)-f(x)\right]\hspace{0.5cm}\text{and}\hspace{0.5cm}\nabla_N^- f(x):=N\left[f(x)-f\left(x-\frac{1}{N}\right)\right]. \] Then, we define the discrete Laplace \\

$\displaystyle \hspace{1cm} \Delta_Nf(x):=\nabla_N^+\nabla_N^-f(x)$\par 
$\displaystyle \hspace{2.7cm}=\nabla_N^-\nabla_N^+f(x)=N^2\left[f\left(x-\frac{1}{N}\right)-2f(x)+f\left(x+\frac{1}{N}\right)\right]$.\\

\noindent If $f\in\H^N$ in particular, then \[ \Delta_Nf(x)=\sum_{j=1}^N\big[N^2(f_{j-1}-2f_j+f_{j+1})\big]\1_j(x). \] 

From the spectral analysis of $\Delta_N$, it is well known that, if $N$ is an odd integer, letting $0\leq m\leq N-1$ with $m$ even, letting $\varphi_{0,N}\equiv 1$, 
$\displaystyle \varphi_{m,N}(x)=\sqrt{2}cos\big(\pi mjN^{-1}\big)$ and $\displaystyle \psi_{m,N}(x)=\sqrt{2}sin\big(\pi mjN^{-1}\big)$ for $\displaystyle x\in I_j$, then, $\displaystyle \big\{\varphi_{m,N},\psi_{m,N}\big\}$ are eigenfunctions of $\Delta_N$ with eigenvalues given by $\displaystyle -\beta_{m,N}=-2N^2\left(1-cos\big(\pi mN^{-1}\big)\right)\leq 0$. If $N$ is even, we need the additional eigenfunction $\displaystyle \varphi_{N,N}=cos(\pi j)$ for $x\in I_j$. The following (classical) properties are derived from \cite{Blount1987}, Lemma 2.12 p.12, \cite{Blount1992}, Lemma 4.2 for the parts (i)-(vi), and from \cite{Kato1966}, chapter 9, Section 3 for the part (vi).

\begin{proposition}\label{Discrete_Laplacian_properties}\textbf{\texttt{(Some properties of the discrete Laplace)}}\vspace{0.1cm}\par 

\noindent\hspace{0.2cm} \textbf{(i)} The family $\{\varphi_{m,N},\psi_{m,N}\}$ forms an orthonormal basis of $\displaystyle \left(\mathds{H}^N,\langle\cdot,\cdot\rangle_2\right)$.\vspace{0.2cm}\par 

\noindent Consider $f,g\in\mathds{H}^N$ and let $\displaystyle T_N(t)=e^{\Delta_Nt}$ denote the semigroup on $\mathds{H}^N$ generated by $\Delta_N$.\vspace{0.1cm}\par 

\noindent\hspace{0.2cm} \textbf{(ii)} $\disp T_N(t)f=\sum_{m}e^{-\beta_{m,N}t}\left(\langle f,\varphi_{m,N}\rangle_2\varphi_{m,N}+\langle f,\psi_{m,N}\rangle_2\psi_{m,N}\right)$.\vspace{0cm}\par 

\noindent\hspace{0.2cm} \textbf{(iii)} $\langle \nabla_N^+f,g\rangle_2=\langle f,\nabla_N^-g\rangle_2$\hspace{0.1cm} and \hspace{0.1cm}$T_N(t)\Delta_Nf=\Delta_NT_N(t)f$. \vspace{0.1cm}\par 

\noindent\hspace{0.2cm} \textbf{(iv)} $\Delta_N$ and $T_N(t)$ are self-adjoint on $\displaystyle \left(\mathds{H}^N,\langle\cdot,\cdot\rangle_2\right)$.\vspace{0.1cm}\par 

\noindent\hspace{0.2cm} \textbf{(v)} $T_N(t)$ is a positive contraction semigroup on both $\displaystyle \left(\mathds{H}^N,\langle\cdot,\cdot\rangle_2\right)$ and $\displaystyle \left(\mathds{H}^N,\Vert\cdot\Vert_\infty\right)$.\vspace{0.1cm}\par 

\noindent\hspace{0.2cm} \textbf{(vi)} The projection $P_N$ commutes with $\Delta_N$, and $\forall T>0$, $f\in C^3(I)$,
\[ \left\{
\begin{array}{l}
\disp \Vert \Delta_Nf-\Delta f\Vert_\infty\longrightarrow0\hspace{0.3cm}\text{as}\hspace{0.2cm}N\rightarrow\infty,\vspace{0.1cm}\\
\disp \Vert \Delta_NP_Nf-\Delta f\Vert_\infty\longrightarrow0\hspace{0.3cm}\text{as}\hspace{0.2cm}N\rightarrow\infty,\vspace{0.1cm}\\
\disp \sup_{t\in[0,T]}\Vert T_N(t)P_Nf-T(t)f\Vert_\infty\longrightarrow0\hspace{0.3cm}\text{as}\hspace{0.2cm}N\rightarrow\infty.
\end{array} 
\right. \]

\noindent We introduce the operator 
$\displaystyle 
\tilde{\Delta}_N:=
\begin{pmatrix}
\Delta_N & 0 \\ 0 & 0
\end{pmatrix} 
$ on $\displaystyle \mathds{H}^N\times \mathds{H}^N$, and denote by $\displaystyle \tilde{T}_N(t):=e^{\tilde{\Delta}_Nt}$ the associated semigroup.\vspace{0.1cm}\par 

\noindent\hspace{0.2cm} \textbf{(vii)} $\tilde{\Delta}_N$ is a bounded linear operator on $\displaystyle \left(\mathds{H}^N\times \mathds{H}^N,\Vert\cdot\Vert_{\infty,\infty}\right)$.\vspace{0cm}\par 

\noindent\hspace{0.2cm} \textbf{(viii)} $\displaystyle \tilde{T}_N(t)=
\begin{pmatrix}
T_N(t) & 0 \\ 0 & I_d
\end{pmatrix},
$ and defines a bounded positive contraction semigroup on $\displaystyle \left(\mathds{H}^N\times \mathds{H}^N,\Vert\cdot\Vert_{\infty,\infty}\right)$. We have denoted by $I_d$ is the identity operator on $\mathds{H}^N$.

\end{proposition}

\noindent \textbf{\textit{Debit functions and formal limit of the generator.}}\label{debit_function_as_generator_derivative_2} We consider the generator $\cA^N$ of (\hyperlink{M6}{M6}) given by $(\ref{Hybrid_Stochastic_Generator})$, applied to a  $C_b^2\big(\H^N\times \H^N\big)$ test function $\varphi$. We use a Taylor expansion of order 2 of $\varphi$. From \Cref{debit_function_as_generator_derivative_1}, we know that the generator is (informally) obtained as the image of the corresponding debit by the differential of the test function. 

In \Cref{formal_limit_of_the_generator}, we give heuristics showing that the terms related to the second order of the Taylor expansion vanish, under the very strong condition $\displaystyle \mu^{-1}N^2\rightarrow0$ as $N,\mu\rightarrow\infty$. In our rigorous proof below, we only assume $\mu^{-1}\log N\rightarrow 0$.

Concerning the first order terms, we identify the differential operator with the gradient as usual, and get\\ 

\small
\noindent $\disp\hspace{0.2cm} \big(D\varphi(u_C,u_D),\Psi^N(u_C,u_D)\big)$\par 
$\disp \hspace{0.3cm}=\left\langle D^{1,0}\varphi(u_C,u_D),\sum_{j=1}^N\left(\sum_{r\in\mathfrak{R}_C}\gamma_r^C\lambda_r\big(u_j^C\big)+\sum_{r\in S_1}\gamma_r^C\lambda_r\big(u_j^C,u_j^D\big)\right)\mathds{1}_j\right\rangle_2$\par
$\disp \hspace{1.5cm}+\left\langle D^{1,0}\varphi(u_C,u_D),\sum_{j=1}^N\left(\sum_{r\in\mathfrak{R}_{DC}\backslash S_1}\frac{\gamma_{j,r}^{N,C}(u_C,u_D)}{\mu}\lambda_r\big(u_j^C,u_j^D\big)\right)\1_j\right\rangle_2$\par  
$\disp \hspace{1.5cm}+\left\langle D^{1,0}\varphi(u_C,u_D), \sum_{j=1}^NN^2\left(\mathds{1}_{j-1}-2\mathds{1}_{j}+\mathds{1}_{j+1}\right)u_j^C\right\rangle_2$\par 

\noindent $\disp \hspace{0.1cm}+\left\langle D^{0,1}\varphi(u_C,u_D),\sum_{j=1}^N\left(\sum_{r\in\mathfrak{R}_{DC}\backslash S_1}\gamma_{j,r}^{N,D}(u_C,u_D)\lambda_r\big(u_j^C,u_j^D\big)+\sum_{r\in \fR_D}\gamma_{j,r}^{N,D}(u_C,u_D)\lambda_r\big(u_j^D\big)\right)\right\rangle_2$,\\

\normalsize
\noindent where $D^{1,0}$ (resp. $D^{0,1}$) denotes the differential with respect to $u_C$ (resp. $u_D$). 
 
On the one hand, we derive the debit function $\displaystyle \Psi_C^N:\H^N\times \H^N\longrightarrow \H^N$, related to the reactant $C$ and defined by 
\begin{equation}\label{debit_function_related_to_C}
\disp \Psi_C^N(u_C,u_D):=\Delta_Nu_C+F(u_C,u_D)+F_1^N(u_C,u_D).
\end{equation}
The first term on the r.h.s. is the debit of the diffusions of molecules of $C$. As in $(\ref{debit_of_the_one_scale_stochastic_spatial_model})$, it is obtained using the periodicity and a change of index: 
\begin{equation}\label{Deterministic_discretized_diffusion_debit_function} 
 \sum_{j=1}^N N^2 \left(\mathds{1}_{j-1}-2\mathds{1}_{j}+\mathds{1}_{j+1}\right)u_j^C
 = \sum_{j=1}^NN^2\left(u_{j+1}^C-2u_{j}^C+u_{j-1}^C\right)\mathds{1}_{j}=\Delta_Nu_C.
\end{equation}
\noindent The second term is associated with fast onsite reactions (they all influence $C$). It is given by
\begin{equation}\label{Deterministic_continuous_debit_function}
\begin{tabular}{l}
$\displaystyle F\left(u_C,u_D\right)=\sum_{j=1}^N\left(\sum_{r\in\mathfrak{R}_C}
\gamma_r^C\lambda_r\left(u_j^C\right)
+\sum_{r\in S_1}\gamma_r^C\lambda_r\left(u_j^C,u_j^D\right)\right)\mathds{1}_j$\vspace{0.1cm}\\

$\displaystyle\hspace{2cm}=\sum_{r\in\mathfrak{R}_C}\gamma_r^C\lambda_r\left(u_C\right)
+\sum_{r\in S_1}\gamma_r^C\lambda_r\left(u_C,u_D\right)$.
\end{tabular}
\end{equation}
The third term corresponds to slow onsite reactions influencing $C$ and vanishes at the limit. It is defined by
\begin{equation}\label{debit_function_retated_to_slow_onsite_involving_C}
\begin{array}{ll}
\disp F_1^N(u_C,u_D) & \disp =\frac{1}{\mu}\sum_{j=1}^N \sum_{r\in\mathfrak{R}_{DC}\backslash S_1}\gamma_r^C\sum_{i=1}^N\gamma_{ij}^N\1_i\theta_{ij}^r(u_i^C,u_i^D)\lambda_r(u_j^C,u_j^D)
\end{array}
\end{equation} 

\noindent On the other hand, we have $\displaystyle \Psi_D^N:\H^N\times \H^N \longrightarrow \H^N$, the debit related to $D$. It corresponds to slow onsite reactions - only them has an effect on $D$ - and is given by
\begin{equation}\label{debit_function_related_to_D}
\disp \Psi_D^N(u_C,u_D)=\sum_{j=1}^N\sum_{i=1}^N\gamma_{ij}^N\1_ig_{ij}(u_C,u_D)=:G^N(u_C,u_D),
\end{equation}
where 
\[ g_{ij}(u_C,u_D) =\sum_{r\in\mathfrak{R}_{DC}\backslash S_1}\gamma_r^D\theta_{ij}^r\big(u_i^C,u_i^D\big)\lambda_r\big(u_j^C,u_j^D\big)+\sum_{r\in\mathfrak{R}_D}\gamma_r^D\theta_{ij}^r\big(u_i^C,u_i^D\big)\lambda_r\big(u_j^D\big). \]

\noindent The global debit on $\H^N\times\H^N$ then reads \[ \Psi^N(u_C,u_D)=\left(\Psi_C^N(u_C,u_D),\Psi_D^N(u_C,u_D)\right). \] 

We finally define the extra debit type function $\displaystyle G:C_p(I)\times C_p(I) \longrightarrow C_p(I)$, by
\begin{equation}\label{Deterministic_debit_function_related_to_D}
\displaystyle G(u_C,u_D)(x):=\theta(u_C(x))\theta(u_D(x))\int_0^1g(u_C,u_D)(y)a(x-y)dy
\end{equation}
for all $x\in I$.
 It is the limit of $G^N$ (see \cref{proof_of_debit_function_regularity_and_convergence}). Here, \[ g(u_C,u_D) =\sum_{r\in\mathfrak{R}_{DC}\backslash S_1}\gamma_r^D\lambda_r(u_C,u_D)+\sum_{r\in\mathfrak{R}_D}\gamma_r^D\lambda_r(u_D). \]
\begin{remark}

We notice that $\Delta_N$, $F$, $F_1^N$ and $G^N$ map $C_p(I)\times C_p(I)$ on $C_p(I)$. This follows from the definition of the discrete Laplace, and from we are considering polynomial reaction rates. 
Here above, we have considered these three functions on $\H^N\times\H^N$, as they were introduced to define the debit function, on $\H^N\times\H^N$. 

Besides, the debit can  be extended to $C_p(I)\times C_p(I)$, not by taking $\Delta_N$, $F$, $F_1^N$ and $G^N$ as maps on $C_p(I)\times C_p(I)$, but using the projection $P_N$ introduced earlier, in a similar way as for extending the generator. Namely, the global debit extends to \[ (u_C,u_D)\mapsto\left(\Psi_C^N(P_Nu_C,P_Nu_D),\Psi_D^N(P_Nu_C,P_Nu_D)\right). \]
\indent Furthermore, it should be emphasized that each of $F$, $F_1^N$ and $\lambda_r$, $r\in\mathfrak{R}$ is a function of real variables, actually. Still, we use the notation $F(u)$ for $\displaystyle u\in C_p(I)\times C_p(I)$. That is, for a given $u$, we define a function $f(u): I\rightarrow \mathds{R}$ by \[ f(u)(x)=F(u(x)),\hspace{0.5cm}\forall x\in [0,1]. \] If $F$ is (locally) Lipschitz, then $f$ is (locally) Lipschitz too, accordingly, with the corresponding norms. We do not distinguish between $F$ and $f$ in our notations, and do the same for the functions $F_1^N$ and $\lambda_r$, $r\in\mathfrak{R}$.\par 
\end{remark}

Fix $u=(u_C,u_D)\in C_p(I)\times C_p(I)$ and let $N,\mu\rightarrow\infty$. Clearly, $F_1^N\big(\tilde{P}_Nu\big)\rightarrow0$. From \cref{Piecewise_Approximation_via_Projection} and the continuity of $F$, $F\big(\tilde{P}_Nu\big)\rightarrow F(u)$. Therefore, if $u_C\in C^3(I)$, $\Psi_C^N\big(\tilde{P}_Nu\big)\rightarrow \Delta u_C+F(u)$ thanks to \cref{Discrete_Laplacian_properties} (vi). We also prove below that $\Psi_D^N\big(\tilde{P}_Nu\big)=G^N\big(\tilde{P}_Nu\big)\rightarrow G(u)$. Therefore, under the strong assumption $\displaystyle \mu^{-1}N^2\rightarrow 0$ and $u_C\in C^3(I)$, we obtain the limit generator
\begin{equation}\label{limit_generator}
\begin{tabular}{l}
$\displaystyle \mathcal{A}^\infty\varphi\left(u_C,u_D\right)=
\left\langle D^{1,0}\varphi(u_C,u_D),\Delta u_C+F(u_C,u_D)\right\rangle_2$\vspace{0.2cm}\\
$\displaystyle \hspace{6cm}+\left\langle D^{0,1}\varphi(u_C,u_D),G(u_C,u_D)\right\rangle_2$.
\end{tabular}
\end{equation} 

The condition $\displaystyle \mu^{-1}N^2\rightarrow 0$ was used By Arnold and Theodosopulu in \cite{Arnold1980}, to prove a LLN for (\hyperlink{M4}{M4}) in the $L^2$ norm. However, that condition is not optimal. Indeed, in \cite{Blount1992}, Blount proved a LLN for (\hyperlink{M4}{M4}) in the supremum norm, requiring only $\displaystyle \mu^{-1}\log N\rightarrow0$. Our convergence result for (\hyperlink{M6}{M6}) falls within the latter framework.\\

We now look for a process admitting $\mathcal{A}^\infty$ defined by $(\ref{limit_generator})$ as infinitesimal generator.

\subsubsection{The limiting problem}\label{the_limiting_problem_1}
\noindent \textbf{\textit{Well-posedness.}} Proceding as in \cref{model_M2}, we see that to $\mathcal{A}^\infty$, corresponds \[ \frac{\partial v_C}{\partial t}=\Delta v_C+F(v_C,v_D)\hspace{0.5cm} \text{(PDE)\hspace{0.2cm} coupled to\hspace{0.2cm} (ODE)}\hspace{0.5cm}\frac{dv_D}{dt}=G(v_C,v_D). \] Hence,  we shall consider the following system of differential equations, for $x\in I$, $t\geq 0$, with periodic boundary conditions and initial data $v_0=\big(v_0^C,v_0^D\big)$:
\[
\left\{
\begin{array}{l}
\displaystyle \frac{\partial}{\partial t}v_C(t,x)=\Delta v_C(t,x)+F(v_C(t,x),v_D(t,x))\vspace{0.1cm}\\
\displaystyle \frac{\partial}{\partial t}v_D(t,x)=G(v_C(t,x),v_D(t,x))\vspace{0.1cm}\\
\displaystyle v_C(t,0)=v_C(t,1)\hspace{0.5cm}\text{and}\hspace{0.5cm}v_D(t,0)=v_D(t,1)\hspace{0.5cm}\forall t\geq0\vspace{0.1cm}\\
\disp v_C(0,x)=v_0^C(x)\geq0\hspace{0.5cm}\text{and}\hspace{0.5cm}v_D(0,x)=v_0^D(x)\geq0\hspace{0.5cm}\forall x.
\end{array}
\right.
\]
\noindent A compact form of the system is:
\begin{equation}\label{Limiting_problem_compact_form}
\left\{
\begin{array}{l}
\displaystyle \frac{d}{dt}v(t)=\tilde{\Delta}v(t)+R(v(t))\vspace{0.1cm}\\
\displaystyle v(t,0)=v(t,1)\hspace{0.5cm}\forall t\geq0\vspace{0.1cm}\\
\disp v(0)=v_0\geq0,
\end{array}
\right.
\end{equation}
where 
$\tilde{\Delta}=
\begin{pmatrix}
\Delta & 0 \\ 0 & 0
\end{pmatrix}:C_p(I)\times C_p(I)\rightarrow C_p(I)$ is a linear operator we define on the domain $C^2(I)\times C_p(I)$, and whose associated semigroup $\tilde{T}(t):=e^{\tilde{\Delta}t}$ is of contraction. For $u\in C_p(I)\times C_p(I)$,  
$\displaystyle R(u)=
\begin{pmatrix}
F (u)\\ G(u)
\end{pmatrix}$, with $G(u)=\theta(u_C)\theta(u_D)a*g(u)$. 

Before, we give the counterpart of $\cref{well-posedness_assumpt_on_the_deb_funct_of_react_for_the_det_loc_model}$ for $(\ref{Limiting_problem_compact_form})$, let us introduce the $\H^N$-valued function $|g|$ defined by $|g|(u):=\sum_{j=1}^N|g|_j(u)\1_j$  for $u=(u_C,u_D)\in C_p(I)\times C_p(I)$, where \[ |g|_j(u_C,u_D):=\sum_{r\in\cR_{DC}\backslash S_1}|\gamma_r^D|\lambda_r(u_j^C,u_j^D)+\sum_{r\in\cR_{D}}|\gamma_r^D|\lambda_r(u_j^D), \] so that 
\[ |g|(u)=\sum_{r\in\cR_{DC}\backslash S_1}|\gamma_r^D|\lambda_r(u_C,u_D)+\sum_{r\in\cR_{D}}|\gamma_r^D|\lambda_r(u_D). \] 
Below, we introduce different quantities and functions such as $|g|$, related to debit type functions, and explain the intuition behind them. We call $|g|$ the "amplitude function" associated with the debit type function $g$.
Then as usual, we may view $|g|$ as a real valued function defined on $\R^2$. Now, we make the following
\begin{assumption}\label{well-posedness_assumpt_on_the_deb_funct_of_react_for_the_det_loc_multiscale_model} \hspace{1cm} \par 
\textbf{(C1)} $F$ is locally Lipschitz and $F(y)\geq 0$ for $y=(y_1,y_2)\in\R^2$ such that $y_1=0$.\par
\textbf{(C2)} There exists $0<\rho_C<\infty$ such that $|y|>\rho_C$ yields $F(y)<0$ for all $y\in\R^m$.\vspace{0.1cm}\par 

\textbf{(D1)} $\theta(y)= 0$ for $y\le 0$. \par 
\textbf{(D2)} $g$ is locally Lipschitz, and $|g|$ has at most linear growth w.r.t. its second variable: for all $\bar{c}>0 $ there exists $M_1(\bar{c})\geq0$ such that, for $|y_1|\le \bar{c}$, $|g|(y_1,y_2)\leq M_1(\bar{c})(|y_2|+1)$.
\end{assumption}

Concerning the reactant $C$, \cref{well-posedness_assumpt_on_the_deb_funct_of_react_for_the_det_loc_multiscale_model} (C1) and (C2) correspond to \cref{well-posedness_assumpt_on_the_deb_funct_of_react_for_the_det_loc_model}. These conditions are met in particular, as soon as the rate of each fast reaction satisfies (C1) and (C2) in the place of $F$. Yet, there is difference with $D$. First of all, in order to ensure positivity for its concentration, we have introduced the function $\theta$ in \cref{spatial_correlation_for_discrete_variable}. Moreover, the debit associated with $D$ involves a convolution product
at the limit. It seems difficult to write an assumption similar to (C2) and we assume linear growth which is less general.

%


The upcoming result states that $\mathcal{A}^\infty$ is the generator of a unique process. It is a straightforward adaptation of the result of Kotelenez in \cite{Kotelenez1986bis}, about well-posedness for the one scale deterministic spatial model. 

\begin{proposition}\label{Limit_system_Well_posedness}Assume that:\vspace{0.1cm}\par 

\noindent \hspace{0.2cm}\textbf{(i)} the rates of onsite reactions are polynomial and such that \cref{well-posedness_assumpt_on_the_deb_funct_of_react_for_the_det_loc_multiscale_model} holds,\vspace{0cm}\par 

\noindent \hspace{0.2cm}\textbf{(ii)} $\displaystyle v_0=\left(v_0^C,v_0^D\right)\in C^3(I)\times C^3(I)$,\vspace{0cm}\par 

\noindent \hspace{0.2cm}\textbf{(iii)} $\rho_C,\rho_D\ge 1$ such that $\Vert v_0^C\Vert_{\infty}< \rho_C<\infty$ and $\Vert v_0^D\Vert_{\infty}< \rho_D<\infty$.\vspace{0.1cm}\par 

\noindent Then, $(\ref{Limiting_problem_compact_form})$ has a unique global mild solution $v:=v(t,v_0)$ satisfying 
\begin{equation}\label{Limit_regularity}
\displaystyle v=(v_C,v_D)\in C\big(\R_+;C^3(I)\times C^3(I)\big)
\end{equation}
and 
\begin{equation}\label{Limiting_problem_compact_integral_form}
v(t)=\tilde{T}(t)v_0+\int_0^t\tilde{T}(t-s)R(v(s))ds\hspace{0.5cm}\forall t\geq0.
\end{equation}
Moreover, 
\begin{equation}\label{limit_boundary}
\left\{
\begin{array}{l}
\disp v_C(t)\geq0,\quad v_D(t)\geq0\hspace{2.1cm}\forall t\geq0, \vspace{0.1cm}\\
\displaystyle \Vert v_C(t)\Vert_{\infty}\leq\rho_C\hspace{0.5cm}\forall t\geq0,\vspace{0.1cm}\\
\displaystyle \Vert v_D(t)\Vert_{\infty}\leq(\rho_D+1)e^{a(0)M_1(\rho_C)t}\hspace{0.5cm}\forall t\geq0.
\end{array}
\right.
\end{equation}\vspace{0cm}

\end{proposition}

\vspace{0.2cm}
\noindent \textbf{\textit{A discretization of the limit.}} We define a discrete version $v^N=(v_C^N,v_D^N)$ of the limiting problem. Henceforth, $v=(v_C,v_D)$ denotes the solution to $(\ref{Limiting_problem_compact_form})$ given by $\cref{Limit_system_Well_posedness}$. For $x\in I$, $t\geq0$, we consider the system: 
\[ \left\{
\begin{array}{l}
\displaystyle \frac{\partial}{\partial t}v_C^N(t,x)=\Delta_N v_C^N(t,x)+F\left(v_C^N(t,x),v_D^N(t,x)\right)\vspace{0.1cm}\\
\displaystyle \frac{\partial}{\partial t}v_D^N(t,x)=G\left(v_C^N(t,x),v_D^N(t,x)\right),\vspace{0.2cm}\\
\displaystyle v_C^N(t,0)=v_C^N(t,1)\hspace{0.1cm}\text{ and }\hspace{0.1cm} v_D^N(t,0)=v_D^N(t,1), \vspace{0.2cm}\\
\displaystyle v_C^N(0,x)=P_Nv_C(0)\hspace{0.1cm}\text{ and }\hspace{0.1cm} v_D^N(0):=P_Nv_D(0)\hspace{0.2cm}.
\end{array}
\right. \]

\noindent Note that $\displaystyle v^N(0)=\tilde{P}_Nv(0)$,  where we recall $\displaystyle \tilde{P}_Nv:=\left(P_Nv_C,P_Nv_D\right)$. Using the operators introduced in $\cref{Discrete_Laplacian_properties}$, a compact version for this system of ODEs reads
\begin{equation}\label{discrete_limiting_problem_compact_form}
\left\{
\begin{array}{l}
\displaystyle \frac{d}{dt}v^N(t)=\tilde{\Delta}_Nv^N(t)+R\left(v^N(t),v^N(t)\right)\vspace{0.1cm}\\
\displaystyle v^N(t,0)=v^N(t,1) \vspace{0.1cm}\\
\displaystyle \displaystyle v^N(0)=\tilde{P}_Nv(0).
\end{array}
\right.
\end{equation}

The next result gives a relation between the limiting problem and its discretization. Its proof is reported to Appendix $\ref{Discrete_Continuous_relation_proof}$.

\begin{theorem}\label{Discrete_Continuous_relation} If the rates of onsite reactions are polynomial and $\cref{well-posedness_assumpt_on_the_deb_funct_of_react_for_the_det_loc_multiscale_model}$ holds, then
$(\ref{discrete_limiting_problem_compact_form})$ has a unique global mild solution $\displaystyle v^N=(v_C^N,v_D^N)\in C\big(\R_+;\H^N \times \H^N\big)$ satisfying
\begin{equation}\label{discrete_limiting_problem_compact_mild_form}
\displaystyle v^N(t)=\tilde{T}_N(t)\tilde{P}_Nv(0)+\int_0^t\tilde{T}_N(t-s)
\begin{pmatrix}
F\left(v^N(s)\right) \vspace{0.1cm}\\ G\left(v^N(s)\right)
\end{pmatrix}ds
\end{equation} 
for all $t\geq0$. Furthermore, for any fixed $T>0$, 
\begin{equation}\label{discretization_boundedness}
\disp \left\{
\begin{array}{l}
\disp \Vert v_C^N(t)\Vert_\infty\leq (\rho+1)/2\hspace{1.5cm}\forall t\geq 0,\vspace{0.1cm}\\
\disp \Vert v_D^N(t)\Vert_{\infty}\leq \bar{c}_D=\bar{c}_D(M_1,T,a(0)),\hspace{0.5cm}\forall 0\leq t\leq T,
\end{array}
\right.
\end{equation}
and
\begin{equation}\label{discretization_convergence}
\disp \sup_{t\in[0,T]}\left\Vert v^N(t)-v(t)\right\Vert_{\infty,\infty}\longrightarrow 0.
\end{equation}

\end{theorem}


\section{The Law of Large Numbers}

Now, we state and prove our main result.

\begin{theorem}\label{Hybrid_law_of_large_numbers}
Consider a sequence $\displaystyle u^N=\left(u_C^N,u_D^N\right)$ of Markov processes starting at $\displaystyle u^N(0)=\left(u_C^N(0),u_D^N(0)\right)$, with infinitesimal generators $\displaystyle \mathcal{A}^N$ given by $(\ref{Stochastic_Generator})$. Assume that the assumptions of \cref{Limit_system_Well_posedness} 
hold and:\vspace{0.1cm}\par 

\textbf{(i)} $N,\mu\longrightarrow\infty$ in such a way that $\displaystyle \mu^{-1}\log N\longrightarrow0$,\par 
\textbf{(ii)} $\disp u^N(0)\geq0$ \hspace{0.1cm} and \hspace{0.1cm} $\displaystyle \left\Vert u^N(0)-v(0)\right\Vert_{\infty,\infty}\longrightarrow 0$ in probability.\vspace{0.1cm}\par 

\noindent Then, for all $T>0$,
\begin{equation}\label{LLN}
\displaystyle \sup_{t\in [0,T]}\left\Vert u^N(t)-v(t)\right\Vert_{\infty,\infty}\longrightarrow 0\hspace{0.5cm}\text{in\hspace{0.1cm} probability}.
\end{equation}

\end{theorem}

\noindent \textbf{\textit{Proof.}} 
Let $T>0$ be fixed. Our goal is to show that,\vspace{0.2cm}\par 
$\disp\hspace{3cm} \forall\epsilon>0,\hspace{0.5cm}\mathds{P}\left\{\sup_{[0,T]}\left\Vert u^N(t)-v(t)\right\Vert_{\infty,\infty}>\epsilon\right\}\longrightarrow 0$.\\ 

\noindent It follows from \cref{Discrete_Continuous_relation} that it is sufficient to prove
\begin{equation}\label{CQFD_1}
\disp \forall\epsilon>0,\hspace{0.5cm} \mathds{P}\left\{\sup_{[0,T]}\left\Vert u^N(t)-v^N(t)\right\Vert_{\infty,\infty}>\epsilon\right\}\longrightarrow 0.
\end{equation}

 The rest of the proof is divided in two principal steps. In the first step, we successively consider some martingales associated to our model (\hyperlink{M6}{M6}) and, an adequate truncation of our process in time, using the properties of these martingales. The aim is to work only with the truncated model in the place of the initial one. Next, we conclude using a Gronwall-Bellman argument.\par 

\subsection{Accompanying martingales}\label{accompanying_martingales}

Different types of martingales are associated to jump Markov processes such as $u^N$. Before specifying some of them, we introduce: \\

\noindent \textbf{\textit{Some useful notations.}} Let $j=1,\cdots,N$ and $u=(u_C,u_D)\in C_p(I)\times C_p(I)$ be fixed.\vspace{0.2cm}\par 

$\bullet$ $\lambda_j^C(u)$\hspace{0.1cm} $\big($resp. $\lambda_j^D(u)\big)$ is the rate for fast (resp. slow) reactions on site $j$:\par 
$\disp\hspace{1.5cm} \lambda_j^{C}(u_C,u_D) = \mu\left(\sum_{r\in\mathfrak{R}_C}\lambda_r\left(u_j^C\right)+\sum_{r\in S_1}\lambda_r\left(u_j^C,u_j^D\right)+2N^2u_j^C\right)$,\par 
\noindent and \par 
$\disp\hspace{3cm} \lambda_j^{D}(u_C,u_D)=\sum_{r\in\mathfrak{R}_{DC}\backslash S_1}\lambda_r\left(u_j^C,u_j^D\right)+\sum_{r\in\mathfrak{R}_D}
\lambda_r\left(u_j^D\right)$.\\

$\bullet$ $\lambda_j(u):=\lambda_j^{C}(u)+\lambda_j^{D}(u)$ is the total rate for reactions on site $j$.

$\bullet$ $\displaystyle \lambda^N(u):=\sum_{j=1}^N\lambda_j(u)$ is the total rate for reactions in the system: \par
 
$\displaystyle\hspace{0.7cm} \lambda^N(u_C,u_D) = \mu\sum_{j=1}^N\left(\sum_{r\in\mathfrak{R}_C}\lambda_r\left(u_j^C\right)+\sum_{r\in S_1}\lambda_r\left(u_j^C,u_j^D\right)+2N^2u_j^C\right)$\par 

$\displaystyle\hspace{3.5cm}+\sum_{j=1}^N\left(\sum_{r\in\mathfrak{R}_{DC}\backslash S_1}\lambda_r\left(u_j^C,u_j^D\right)+\sum_{r\in\mathfrak{R}_D}
\lambda_r\left(u_j^D\right)\right)$.\\

$\bullet$ If the reaction rates are bounded,  we set $\displaystyle \bar{\lambda}:=\max\left(\bar{\lambda}_{C},\bar{\lambda}_{D}\right)$ where, \[ \bar{\lambda}_{C}=\sum_{r\in\mathfrak{R}_C\cup S_1}\left\Vert \lambda^C\right\Vert_\infty
\hspace{0.5cm} \text{with}\hspace{0.5cm} \left\Vert \lambda^C\right\Vert_\infty:=\sup\left\{\Vert\lambda_r\Vert_\infty, r\in\mathfrak{R}_{C}\cup S_1\right\}. \] The notation $\displaystyle \bar{\lambda}_D$ is the counterpart of $\displaystyle \bar{\lambda}_C$ for slow reactions.\\

$\bullet$ $\disp \left|\gamma\right|:=\max\left(\left|\gamma^C\right|,\left|\gamma^D\right|\right) \hspace{0.5cm} \text{ and } \hspace{0.8cm}  \bar{\gamma}:=\max\left(\bar{\gamma}_C\hspace{0.1cm},\hspace{0.1cm}\bar{\gamma}_D\right)$. 

$\displaystyle\hspace{1cm} \bar{\gamma}_D:=\max\left(\sum_{r\in\mathfrak{R}_{DC}\backslash S_1}\left|\gamma^D\right|,\sum_{r\in\mathfrak{R}_{D}}\left|\gamma^D\right|\right),\hspace{0.3cm}\text{with}\hspace{0.2cm} \left|\gamma^D\right|:=\sup_{r\in(\mathfrak{R}_{DC}\backslash S_1)\cup \mathfrak{R}_{D}}\left|\gamma_r^D\right|$.\\

\noindent The notations $\disp \bar{\gamma}_C$ and $\displaystyle \left|\gamma^C\right|$ are the counterparts for fast reactions.\\

Recall that the debit of a process is defined as the sum of its jumps weighted by the corresponding rates, over all its possible jump directions. Let \[ \delta u_j^N(t):=u_j^N(t)-u_j^N(t^-) \] be the jump of $u_j^N$ at time $t$, and denote by $|\delta u_j^N(t)|$ the amplitude of that jump, where  $|(y_1,y_2)|=|y_1|+|y_2|$. Firstly, we define the "amplitude" $|\Psi_C^N|_j$ (resp. $|\Psi_D^N|_j$) of the debit $\Psi_j^{N,C}$ (resp. $\Psi_j^{N,D}$), as the debit function of the process \[ \left(\sum_{s\leq t}\big|\delta u_j^{N,C}(t)\big|\right)_{t\geq0}\hspace{1cm}\left[\text{resp. }\hspace{0.2cm} \left(\sum_{s\leq t}\big|\delta u_j^{N,D}(t)\big|\right)_{t\geq0}\right]. \] 
As a fact, that notion of "amplitude" is perfectly adapted to any of the specific debit functions $\Delta_N$, $F$, $F_1^N$, or $G^N$. Especially, we have for the part related to $C$:\vspace{0.2cm}\par 

$\bullet$ $\displaystyle \left|\Delta_{N}\right|_j(u_C)=N^2\left(u_{j-1}^C+2u_j^C+u_{j+1}^C\right)$,\vspace{0.1cm}\par 

$\bullet$ $\displaystyle |F|_j(u_C,u_D)=\sum_{r\in\mathfrak{R}_C}|\gamma_r^C|\lambda_r\left(u_j^C\right)
+\sum_{r\in S_1}|\gamma_r^C|\lambda_r\left(u_j^C,u_j^D\right)$,\par 

$\bullet$ $\disp |F_1^N|_j(u_C,u_D)=\frac{1}{\mu}\sum_{i=1}^N\sum_{r\in\fR_{DC}\backslash S_1}|\gamma_{ij}^N|\theta_{ij}^r(u_i^{C},u_i^D)\lambda_r(u_j^C,u_j^D)$,\vspace{0.1cm}\par 

$\bullet$ $\displaystyle \left|\Psi_C^{N}\right|_j(u_C,u_D)=\left|\Delta_{N}\right|_j(u)+|F|_j(u)+\left|F_1^N\right|_j(u_C,u_D)$.\\ 

\noindent For the part related to the reactant $D$, we have\vspace{0.2cm}\par 

$\bullet$ $\disp |g|_{ij}(u_C,u_D)=\sum_{r\in\mathfrak{R}_{DC}\backslash S_1}\big|\gamma_r^D\big|\theta_{ij}^r(u_i^C,u_i^D)\lambda_r(u_j^C,u_j^D)+\sum_{r\in\mathfrak{R}_{D}}\big|\gamma_r^D\big|\theta_{ij}^r(u_i^C,u_i^D)\lambda_r(u_j^D)$,\\ 
for all $1\leq i\leq N$,\par 

$\bullet$ $\displaystyle \big|\Psi_D^{N}\big|_j(u_C,u_D)=\sum_{i=1}^N|\gamma_{ij}^N||g|_{ij}(u_C,u_D)=:\big|G^N\big|_j(u_C,u_D)$.\vspace{0cm}\par 

\noindent Then we set $\displaystyle \left|\Psi^{N}\right|_j(u):=\left|\Psi_C^{N}\right|_j(u)+\left|\Psi_D^{N}\right|_j(u)$ and $\disp \left|\Psi^N\right|(u):=\sum_{j=1}^N\left|\Psi^N\right|_j(u)$.\\

Next, we analogously define the "square amplitude" $|\Psi_C^N|_j^2$ (resp. $|\Psi_D^N|_j^2$) of $\Psi_j^{N,C}$ (resp. $\Psi_j^{N,D}$), as the debit function of the process \[ \left(\sum_{s\leq t}\big|\delta u_j^{N,C}(t)\big|^2\right)_{t\geq0}\hspace{1cm}\left[\text{resp. }\hspace{0.2cm} \left(\sum_{s\leq t}\big|\delta u_j^{N,D}(t)\big|^2,\right)_{t\geq0}\right]. \] 
Then,\vspace{0cm}\par 

$\bullet$ $\displaystyle \left|\Delta_{N}\right|_j^2(u_C):=N^2\left(u_{j-1}^C+2u_j^C+u_{j+1}^C\right)$,\vspace{0.1cm}\par 

$\bullet$ $\displaystyle |F|_j^2(u_C,u_D)=\sum_{r\in\mathfrak{R}_C}|\gamma_r^C|^2\lambda_r\left(u_j^C\right)+\sum_{r\in S_1}|\gamma_r^C|^2\lambda_r\left(u_j^C,u_j^D\right)$,\par 

$\bullet$ $\disp |F_1^N|_j^2(u_C,u_D)=\frac{1}{\mu^2}\sum_{i=1}^N\sum_{r\in\fR_{DC}\backslash S_1}|\gamma_{ij}^N|^2|\theta_{ij}^r(u_i^{C},u_i^D)|^2\lambda_r(u_j^C,u_j^D)$,\vspace{0.1cm}\par 

$\bullet$ $\displaystyle \big|\Psi_C^{N}\big|_j^2(u)=\left|\Delta_{N}\right|_j^2(u)+|F|_j^2(u)+\left|F_1^N\right|_j^2(u_C,u_D)$.\vspace{0.1cm}\par 
\small
$\bullet$ $\disp |g|_{ij}^2(u_C,u_D)=\sum_{r\in\mathfrak{R}_{DC}\backslash S_1}\big|\gamma_r^D\big|^2|\theta_{ij}^r(u_i^C,u_i^D)|^2\lambda_r\left(u_j^C,u_j^D\right)+\sum_{r\in\mathfrak{R}_{D}}\big|\gamma_r^D\big|^2|\theta_{ij}^r(u_i^C,u_i^D)|^2\lambda_r\left(u_j^D\right)$,\\
\normalsize  for all $1\leq i\leq N$,\par 

$\bullet$ $\displaystyle \big|\Psi_D^{N}\big|_j^2(u_C,u_D)=\sum_{i=1}^N|\gamma_{ij}^N|^2|g|_{ij}^2(u_C,u_D)=:\big|G^N\big|_j^2(u_C,u_D)$.\vspace{0cm}\par 

\noindent And we set $\displaystyle \left|\Psi^{N}\right|_j^2(u):=\left|\Psi_C^{N}\right|_j^2(u)+\left|\Psi_D^{N}\right|_j^2(u)$ and $\disp \left|\Psi^N\right|^2(u):=\sum_{j=1}^N\left|\Psi^{N}\right|_j^2(u)$.\\

In addition, we introduce at last, the very useful "square amplitude function" associated with the debit $\Psi^N$. It is defined by \[ \left|\Psi^N\right|^2(u)=\sum_{j=1}^N\left|\Psi^{N}\right|_j^2(u)\1_j. \] We also denote it by $|\Psi^N|^2$, as the "square amplitude". Both functions are defined on $\H^N\times\H^N$. But still, if the former is $\H^N-$ valued, the latter has non-negative real values. In order to avoid any confusion in the sequel, the notation $|\Psi^N|^2$ refers to the "square amplitude function", unless another precision is made.

Moreover, as for the "square amplitude" introduced first, the "square amplitude function" can be easily derived for each of the specific debit functions we have introduced above. For instance, the "square amplitude function" corresponding to $F$ is just \[ |F|^2(u):=\sum_{j=1}^N|F|_j^2(u)\1_j. \]

\noindent \textbf{\textit{The so-called accompanying martingales.}} It is well known from \cite{Kurtz1970} or Proposition 2.1 of \cite{Kurtz1971} that, if the total reaction rate $\lambda^N$ and the global "amplitude" - that of the global debit $\Psi^N$ - of $u^N$ are bounded, or in other words if
\begin{equation}\label{bounded_total_rate_martingale_condition}
\displaystyle \sup_{(u)\in\mathds{H}^{N}\times\mathds{H}^{N}}\lambda^N(u)<\infty \hspace{0.5cm} \text{and}
\hspace{0.5cm}\sup_{(u)\in\mathds{H}^{N}\times\mathds{H}^{N}}\left\vert\Psi^{N}(u)\right\vert<\infty
\end{equation}

\noindent respectively, then the processes defined by
\begin{equation}\label{Zero_level_accompanying_martingales}
\left\{
\begin{tabular}{l}
$\displaystyle Z^N(t)=u^N(t)-u^N(0)-\int_0^t\hspace{0.1cm} \Psi^N\left(u^N(s)\right) \hspace{0.1cm} ds $\vspace{0.2cm}\\

$\displaystyle Z_C^N(t):=u_C^N(t)-u_C^N(0)-\int_0^t\hspace{0.1cm}\Psi_C^N\left(u^N(s)\right)\hspace{0.1cm} ds $\vspace{0.2cm}\\

$\displaystyle Z_D^N(t):=u_D^N(t)-u_D^N(0)-\int_0^t \Psi_D^N\left(u^N(s)\right)$\\
\end{tabular}
\right.
\end{equation}
are  $\displaystyle \mathcal{F}_t^N=\mathcal{F}_t^{N,\mu}-$martingales, with $\displaystyle \displaystyle Z^N(t)=\left(Z_C^N(t),Z_D^N(t)\right)$ taking values in $\mathds{H}^N\times\mathds{H}^N$, for all $t\geq 0$. Note that $(\ref{bounded_total_rate_martingale_condition})$ holds when the reaction rates are bounded for instance. Furthermore, if $\tau$ is a $\displaystyle \mathcal{F}_t^N$-stopping time such that
\begin{equation}\label{truncation_stopping_time_condition}
\displaystyle \sup_{[0,T]} \left\Vert u^N(t\wedge\tau)\right\Vert_{\infty,\infty}\leq C(T,N,\mu)<\infty,
\end{equation}
then $u^N(t\wedge\tau)$ has a bounded total jump rate, and $Z^N(t\wedge\tau)$, $Z_C^N(t\wedge\tau)$ and $Z_D^N(t\wedge\tau)$ are martingales. In all these, $t\wedge\tau$ denotes the infimum between $t$ and $\tau$.\par 

The following result is similar to Lemma 2.2, p.8. in \cite{Blount1987}. It presents a first type of martingales related to $u^N$. We give its proof in \Cref{proof_about_accompanying_martingales_1}.

\begin{lemma}\label{First_level_acccompanying_martingales} Assume there exists $\tau$ satisfying $(\ref{truncation_stopping_time_condition})$. Then, for all $j\in \lbrace 1,\cdots,N\rbrace$ and all $t\geq 0$, the following define mean zero $\mathcal{F}_t^N$-martingales:
\begin{description}
\item[(Mg1)]
\begin{tabular}{l}
$\displaystyle\hspace{2cm} \sum_{s\leq t\wedge\tau}\left[\delta u_j^{N,C}(s)\right]^2-\frac{1}{\mu}\int_0^{t\wedge\tau}\left(|\Delta|_j^2+|F|_j^2+\mu\big|F_1^{N}\big|^2\right)\left(u^N(s)\right)ds$,
\end{tabular}

\item[(Mg2)-(Mg3)]
\begin{tabular}{l}
$\displaystyle\hspace{0.5cm} \sum_{s\leq t\wedge\tau}\left[\delta u_{j}^{N,C}(s)\right]\left[\delta u_{j\pm 1}^{N,C}(s)\right]+\frac{1}{\mu}\int_0^{t\wedge\tau}N^2\left(u_j^{N,C}(s)+u_{j\pm 1}^{N,C}(s)\right)ds$,
\end{tabular}

\item[(Mg4)]
\begin{tabular}{l}
$\displaystyle\hspace{2cm} \sum_{s\leq t\wedge\tau}\left[\delta u_j^{N,D}(s)\right]^2-\int_0^{t\wedge\tau}\left|\Psi_j^{N,D}\right|^2\left(u^N(s)\right)ds$.
\end{tabular}
\end{description}

\end{lemma}

We move on to the second martingale type. The result is a variant of Lemma 2.16, p.19 in \cite{Blount1987}, and Lemma 1.1 in \cite{Kotelenez1986}. Find the proof in \Cref{proof_about_accompanying_martingales_2}.

\begin{lemma}\label{Second_level_acccompanying_martingales} Consider $\displaystyle \varphi\in \mathds{H}^N$ and $\tau$ satisfying $(\ref{truncation_stopping_time_condition})$. Then for $t\geq 0$, the following are mean 0 martigales:
\begin{description}
\item[(\hypertarget{Mg5}{Mg5})] $\displaystyle \sum_{s\leq t\wedge\tau}\left[\delta\left\langle Z_C^N(s),\varphi\right\rangle_2\right]^2-\frac{1}{N\mu}\int_0^{t\wedge\tau}\left[\left\langle u_C^N(s),\left(\nabla_N^+\varphi\right)^2+\left(\nabla_N^-\varphi\right)
^2\right\rangle_2\right.$\par 
$\disp\hspace{6.5cm}\left.+\left\langle\left(|F|^2+\mu\big|F_1^N\big|^2\right)\left(u^N(s)\right),\varphi^2\right\rangle_2\right]ds$,

\item[(\hypertarget{Mg6}{Mg6})] $\displaystyle \sum_{s\leq t\wedge\tau}\left[\delta\left\langle Z_D^N(s),\varphi(s)\right\rangle_2\right]^2-\frac{1}{N}\int_0^{t\wedge\tau}\left\langle\left|\Psi_D^{N}\right|^2\left(u^N(s)\right),\varphi^2\right\rangle_2\hspace{0.1cm}ds$.
\end{description}

\end{lemma}

Now we look for a stopping time with the same properties as at $(\ref{truncation_stopping_time_condition})$.

\subsection{Truncation}\label{truncation}

Set \[ \tau:=\tau(N,\epsilon_0)=\inf\left\lbrace t\geq 0:\left\Vert u^N(t)-v^N(t)\right\Vert_{\infty,\infty}>\epsilon_0 \right\rbrace, \]
for fixed $\displaystyle \epsilon_0\in \left] 0,1\right[$. By definition of $\tau$, the process $u^N$ exists and takes finite 
value on $[0,\tau]$ . Now, define $\disp \bar{u}^N=\left(\bar{u}_C^N,\bar{u}_D^N\right)$ by
\begin{equation}\label{Truncated_process_definition}
\left\{
\begin{array}{l}
\displaystyle \bar{u}^N(t)=u^N(t\wedge\tau)\hspace{6.5cm}\text{for}\hspace{0.2cm}0\leq t\leq\tau\leq\infty,\vspace{0.2cm}\\

\displaystyle \bar{u}^N(t)=u^N(t\wedge\tau)+\int_{t\wedge\tau}^t\left(\tilde{\Delta}_N\bar{u}^N(s)+R^N\left(\bar{u}^N(s)\right) 
\right)ds\hspace{0.8cm}\text{for}\hspace{0.2cm}\tau<t<\infty,
\end{array}
\right.
\end{equation}
where $\displaystyle R^N(u):=
\begin{pmatrix}
F(u)\vspace{0.1cm}\\ G^N(u)
\end{pmatrix}$ for $u\in C_p(I)\times C_p(I)$. The process $\bar{u}^N$ is obtained by running $u^N$ until time $\tau$ and then running it deterministically afterwards, if $\tau<\infty$.\par 

Note that $\tau$ is a stopping time such that\\

\noindent $\displaystyle\hspace{0.2cm} \mathds{P}\left\{\sup_{[0,T]}\left\Vert u^N(t)-v^N(t)\right\Vert_{\infty,\infty}>\epsilon_0\right\}\leq \mathds{P}\left\{\sup_{[0,T]}\left\Vert u^N(t\wedge\tau)-v^N(t\wedge\tau)\right\Vert_{\infty,\infty}\geq\epsilon_0\right\}$\par 

$\displaystyle \hspace{6cm}\leq \mathds{P}\left\{\sup_{[0,T]}\left\Vert \bar{u}^N(t)-v^N(t)\right\Vert_{\infty,\infty}\geq\epsilon_0\right\}$.\\

\noindent Therefore, it follows from $(\ref{CQFD_1})$ that showing
\begin{equation}\label{CQFD_2}
\mathds{P}\left\{\sup_{[0,T]}\left\Vert \bar{u}^N(t)-v^N(t)\right\Vert_{\infty,\infty}>\epsilon_0\right\}\longrightarrow 0
\end{equation} 
is sufficient for proving the main result.\\


\noindent \textbf{\textit{Boundedness.}} From \cref{Limit_system_Well_posedness}, we know that $\Vert v(t)\Vert_{\infty,\infty}<\rho_{T}$ with 
$\rho_{T}=\rho_C+(\rho_D+1)e^{a(0)M_1(\rho_C)T}$ for all $t\in[0,T]$. Since we are assuming $\displaystyle \left\Vert u^N(0)-v(0)\right\Vert_{\infty,\infty} \rightarrow 0$ in probability, we may, by conditioning on $\displaystyle \left\Vert u^N(0)\right\Vert_{\infty,\infty}<\rho_{T}+1$ if necessary, assume without loss of generality that 
\begin{equation}\label{Markov_sequence_initial_value_boundedness}
\displaystyle \left\Vert u^N(0)\right\Vert_{\infty,\infty}<\rho_{T}+1 \hspace{0.4cm}\text{for all}\hspace{0.3cm}N.
\end{equation}

\noindent From \cref{Discrete_Continuous_relation}, we also have
\[ \sup_{t\in[0,T]}\big\Vert v^N_C(t)\big\Vert_{\infty}\leq \rho_{C},\quad \sup_{t\in[0,T]}\big\Vert v^N_D(t)\big\Vert_{\infty}\leq (\rho_D+1)e^{a(0)M_1(\rho_C)T}. \]
Therefore, by definition of $\tau$, for $\epsilon_0\le 1$,
\begin{equation}\label{Troncated_Markov_sequence_values_boundedness_before}
\displaystyle \left\Vert u^N_D(t\wedge\tau)\right\Vert_{\infty,\infty}<\rho_{C}+1,\quad \left\Vert u^N_D(t\wedge\tau)\right\Vert_{\infty,\infty}< 
(\rho_D+1)e^{a(0)M_1(\rho_C)T}+1\hspace{0.7cm}\text{for}\hspace{0.4cm}0\leq t\leq T.
\end{equation}
Thanks to \cref{well-posedness_assumpt_on_the_deb_funct_of_react_for_the_det_loc_multiscale_model}, we obtain 
\begin{equation}\label{Troncated_Markov_sequence_values_boundedness}
\displaystyle \left\Vert\bar{u}^N(t)\right\Vert_{\infty,\infty}\le \tilde\rho_{T}=\rho_C+1+(\rho_D+2)e^{a(0)M_1(\rho_C+1)T}\hspace{0.5cm}\text{for}\hspace{0.2cm}0\leq t\leq  T.
\end{equation}
\noindent Indeed, given $(\ref{Troncated_Markov_sequence_values_boundedness_before})$, we prove $(\ref{Troncated_Markov_sequence_values_boundedness})$ if for $\tau<t\leq T$, we show the inequalities
\[ \left\{
\begin{array}{l}
\disp \left\Vert\bar{u}_C^N(t)\right\Vert_\infty\le \rho_C+1\hspace{1.5cm}(C\star)\vspace{0.1cm}\\
\disp \left\Vert\bar{u}_D^N(t)\right\Vert_\infty\leq(\rho_D+2)e^{2a(0)M_1(\rho_C+1)T}\hspace{2.5cm}(D\star).
\end{array}
\right. \] 
Suppose that $\disp \bar{u}_{j}^{N,C}(t)>\rho_C+1$ for some $t\in[\tau,T]$, $1\leq j\leq N$. Since $\bar{u}_C^N(t)\in\H^N$, there is $1\leq i\leq N$ and $t_0$
satisfying 
\[ \bar{u}_{i}^{N,C}(t_0)=\sup_{s\in [\tau,T]}\left\Vert \bar{u}_{C}^N(s)\right\Vert_\infty\geq \bar{u}_{j}^{N,C}(t)>\rho_C+1. \] 

Therefore, \cref{well-posedness_assumpt_on_the_deb_funct_of_react_for_the_det_loc_multiscale_model} (i) yields 
\[ F\big(\bar{u}_{i}^{N}(t_0)\big)=F\big(\bar{u}_{i}^{N,C}(t_0),\bar{u}_{i}^{N,D}(t_0)\big)<0. \] 
Moreover 
\[ \Delta_N\bar{u}_{i}^{N,C}(t_0)=\bar{u}_{i+1}^{N,C}(t_0)-2\bar{u}_{i}^{N,C}(t_0)+\bar{u}_{i-1}^{N,C}(t_0)<0. \] 
Combining these two arguments with $(\ref{Truncated_process_definition})$, we get 
\[ \frac{d}{dt}\bar{u}_{i}^{N,C}(t_0)=\Delta_N\bar{u}_{i}^{N,C}(t_0)+F\big(\bar{u}_{i}^{N}(t_0)\big)<0.\] 
By \eqref{Troncated_Markov_sequence_values_boundedness_before} $t_0>\tau$ and necessarily $\frac{d}{dt}\bar{u}_{i}^{N,C}(t_0)\ge 0$ in contradiction with the above.
Thus, ($C\star$) follows. 

Now, let $\tau<t\leq T$ and $1\leq i\leq N$ be fixed. Let $\tau\leq s\leq t$. Observing that $|g_{ij}(u^N(s)|\leq \big\Vert |g|(u^N(s))\Vert_\infty$, \cref{well-posedness_assumpt_on_the_deb_funct_of_react_for_the_det_loc_multiscale_model} (ii) and ($C\star$) yield
\begin{align*}
\left|G^N\big(\bar{u}_i^{N,D}(s)\big)\right| & \leq \sum_{j=1}^N|\gamma_{ij}^N||g_{ij}(u^N(s))| \leq \sum_{j=1}^N\frac{a(0)}{N}\big\Vert |g|(u^N(s))\Vert_\infty\\
& \leq a(0)M_1(\rho_C+1)\big(\Vert \bar{u}_D^{N}(s)\Vert_\infty+1\big).
\end{align*}

\noindent Thus, $\disp \left\Vert G^N\big(\bar{u}_D^{N}(s)\big)\right\Vert_\infty \leq a(0)M_1(\rho_C+1)(\left\Vert \bar{u}_D^{N}(s)\right\Vert_\infty+1)$, and, $(\ref{Truncated_process_definition})$ with $(\ref{Troncated_Markov_sequence_values_boundedness_before})$ lead to
\begin{align*}
\left\Vert \bar{u}_D^{N}(t)\right\Vert_\infty & \leq \left\Vert u_D^{N}(t\wedge\tau)\right\Vert_\infty+\int_\tau^t\left\Vert G^N\big(\bar{u}_i^{N,D}(s)\big)\right\Vert_\infty ds\\
& \leq (\rho_D+1)e^{a(0)M_1(\rho_C)T}+1+M_1(\rho_C+1)a(0)\int_\tau^t\left\Vert \bar{u}_D^{N}(s)\right\Vert_\infty +1ds
\end{align*}

\noindent The result follows from Gronwall lemma.
($D\star$) follows, and $(\ref{Troncated_Markov_sequence_values_boundedness})$ is proved. $\square$\vspace{0cm}\par 

In the sequel, we consider $\bar{c}$ as a generic constant depending on $\tilde{\rho}_T$ and $T$.

\begin{remark}\label{regularity_from_boundedness}
As we are considering the truncated process $\bar{u}^N$ in the following, we consider that reaction rates are bounded and Lipschitz as functions defined on $\R^2$. Indeed, thanks to $(\ref{Troncated_Markov_sequence_values_boundedness})$, we know that the family $\disp \left\{\bar{u}^N(t),t\geq0\right\}_N$ lies in the bounded set \[ \fS_{\tilde{\rho}_T}:=\left\{u\in \H^N\times\H^N:\Vert u\Vert_{\infty,\infty}\leq \tilde{\rho}_T\right\}. \] Thus $\disp\hspace{0cm} \left\{\bar{u}^N(t,x),t\geq0,x\in I\right\}_N \subset \bar{B}(0,\tilde{\rho}_T)$, where $\bar{B}(0,\tilde{\rho}_T)$, the closed ball of $\R^2$ of radius $\tilde{\rho}_T$, and centered at zero. Now, since reaction rates are functions of the concentrations of the truncated process, we only consider their restrictions to the compact $\bar{B}(0,\tilde{\rho}_T)$. These latter are bounded and Lipshitz, as the rates are polynomial. Their bounds depend on $\tilde{\rho}_T$ and are also denoted 
$\bar{\lambda}(\tilde{\rho}_T)$. $\square$
\end{remark}

We have the following result, proved in \Cref{proof_of_debit_function_regularity_and_convergence}.

\begin{proposition}\label{debit_function_regularity_and_convergence}\hspace{1cm}\vspace{0.1cm}\par  
(i) The debit functions $F$, $F_1^N$, $G$ and $G^N$ are Lipschitz.\par 
(ii) For all $u\in C(I)\times C(I)$, $\disp \left\Vert G^N\left(u\right)-G(u)\right\Vert_\infty\longrightarrow 0$. 
\end{proposition}

\noindent \textbf{\textit{Martingale and jumps.}} The stopping time $\tau$ satisfies $(\ref{truncation_stopping_time_condition})$. Thus, \[ Z^N(t\wedge\tau):=u^N(t\wedge\tau)-u^N(0)-\int_0^{t\wedge\tau}\left(\tilde{\Delta}_Nu^N(s)+R^N\left(u^N(s)\right)\right)ds \] is a mean $0$ martingale. From the definition of $\bar{u}^N$,
\begin{equation}\label{compact_form_for_the_truncated_process_u_N_bar}
\bar{u}^N(t)=u^N(0)+\int_0^{t}\left(\tilde{\Delta}_N\bar{u}^N(s)+R^N\left(\bar{u}^N(s)\right)\right)ds+Z^N(t\wedge\tau).
\end{equation}  
In fact, if $t\leq \tau$, then\vspace{0.1cm}\par 

$\displaystyle \bar{u}^N(t)-u^N(0)-\int_0^t\left(\tilde{\Delta}_N\bar{u}^N(s)+R^N\left(\bar{u}^N(s)\right)\right)ds$\par 

$\displaystyle \hspace{1cm}:=u^N(t\wedge\tau)-u^N(0)-\int_0^t\left(\tilde{\Delta}_N{u}^N(s\wedge\tau)+R^N\left(u^N(s\wedge\tau)\right)\right)ds$\par 

$\displaystyle \hspace{1cm}=u^N(t\wedge\tau)-u^N(0)-\int_0^{t\wedge\tau}\left(\tilde{\Delta}_N{u}^N(s)+R^N\left(u^N(s)\right)\right)ds\hspace{0.2cm}=\hspace{0.2cm}Z^N(t\wedge\tau)$.\\ 

\noindent Otherwise $t>\tau$ and, since $\bar{u}^N$ and $u^N$ are the same before time $\tau$,\vspace{0.1cm}\par 

\noindent $\displaystyle \bar{u}^N(t)-u^N(0)-\int_0^t\left(\tilde{\Delta}_N\bar{u}^N(s)+R^N\left(\bar{u}^N(s)\right)\right)ds$\par 

$\displaystyle \hspace{0.4cm}:=u^N(t\wedge\tau)+\int_{t\wedge\tau}^t\left(\tilde{\Delta}_N\bar{u}^N(s)+R^N\left(\bar{u}^N(s)\right)\right)ds$\par 
$\displaystyle \hspace{6cm}-u^N(0)-\int_0^t\left(\tilde{\Delta}_N\bar{u}^N(s)+R^N\left(\bar{u}^N(s)\right)\right)ds$\par 

$\displaystyle \hspace{0.4cm}=u^N(t\wedge\tau)-u^N(0)-\int_0^{t\wedge\tau}\left(\tilde{\Delta}_N\bar{u}^N(s)+R^N\left(\bar{u}^N(s)\right)\right)ds$\par 

$\displaystyle \hspace{0.4cm}=u^N(t\wedge\tau)-u^N(0)-\int_0^{t\wedge\tau}\left(\tilde{\Delta}_N{u}^N(s)+R^N\left(u^N(s)\right)\right)ds\hspace{0.1cm}=\hspace{0.1cm}Z^N(t\wedge\tau)$. $\square$\\

\noindent Also, setting $\disp \delta f:=\left(\delta f_j\right)_{1\leq j\leq N}$ for $f\in\H^N$, the jumps related to $u^N$ and $\bar{u}^N$ satisfy
\begin{equation}\label{Truncated_process_jump_boundedness}
\left\{
\begin{array}{l}
\displaystyle \left\Vert \delta\bar{u}_C^N(t)\right\Vert_{\infty}=\left\Vert\delta Z_C^N(t\wedge\tau)\right\Vert_{\infty}=\left\Vert\delta u_C^N(t\wedge\tau)\right\Vert_{\infty}\leq (1\vee\bar{\gamma}_C)\mu^{-1}\leq \gamma\mu^{-1}\vspace{0.2cm}\\

\displaystyle \left\Vert \delta\bar{u}_D^N(t)\right\Vert_{\infty}=\left\Vert\delta Z_D^N(t\wedge\tau)\right\Vert_{\infty}=\left\Vert\delta u_D^N(t\wedge\tau)\right\Vert_{\infty}\leq \bar{\gamma}_Da(0)N^{-1}\leq\gamma N^{-1}
\end{array}
\right.
\end{equation}
for all $t\geq0$, for some constant $\gamma>0$. We refer to $(\ref{maximum_jump_amplitude_for_slow_reactions})$ in particular for the bound related to $D$.

\subsection{A Gronwall-Bellman argument}

We want to study the difference \[ \bar{u}^N(t)-v^N(t)=\left(\bar{u}_C^N(t)-v_C^N(t)\hspace{0.1cm},\hspace{0.1cm}\bar{u}_D^N(t)-v_D^N(t)\right). \]
Using variation of constant at $(\ref{compact_form_for_the_truncated_process_u_N_bar})$, we get
\begin{equation}\label{mild_compact_form_for_u_bar_N}
\displaystyle \bar{u}^N(t)=\tilde{T}_N(t)u^N(0)+\int_0^t\tilde{T}_N(t-s)
\begin{pmatrix}
F\left(\bar{u}^N(s)\right) \vspace{0.2cm} \\ G^N\left(\bar{u}^N(s)\right)
\end{pmatrix}
ds+Y^N(t),
\end{equation}
where $\displaystyle Y^N(t)=\int_0^t\tilde{T}_N(t-s)dZ^N(s\wedge\tau)$. Note that 
$\displaystyle Z^N\left(s\wedge\tau,\frac{j}{N}\right):=
\begin{pmatrix}
Z_C^N\left(s\wedge\tau,\frac{j}{N}\right) \vspace{0.2cm} \\ Z_D^N\left(s\wedge\tau,\frac{j}{N}\right)
\end{pmatrix}
$, $1\leq j\leq N$, is of bounded variation in $s$, and $\tilde{T}_N$ may be viewed as a $2N\times 2N$ matrix-valued function. $T_N$ is indeed identified to an $N\times N$ matrix-valued function. Hence, $\displaystyle Y^N\left(t,\frac{j}{N}\right)$, $1\leq j\leq N$, is defined as a Stieltjes integral. From the mild forms of $v^N$ and $\bar{u}^N$ at $(\ref{discrete_limiting_problem_compact_mild_form})$ and $(\ref{mild_compact_form_for_u_bar_N})$ respectively,\\

$\displaystyle \bar{u}^N(t)-v^N(t)=\tilde{T}_N(t)\left(u^N(0)-v^N(0)\right)$\\

$\disp \hspace{3cm}+\int_0^t\tilde{T}_N(t-s)
\begin{pmatrix}
F\left(\bar{u}^N(s)\right)-F\left(v^N(s)\right)\vspace{0.2cm} \\ G^N\left(\bar{u}^N(s)\right)-G\left(v^N(s)\right)
\end{pmatrix} + Y^N(t)$.\\

\noindent Since the semigroup $\tilde{T}_N(t)$ is of contraction, \\

\noindent $\displaystyle \left\Vert\bar{u}^N(t)-v^N(t)\right\Vert_{\infty,\infty}$\\

$\disp \hspace{0cm} \leq \left\Vert u^N(0)-v^N(0)\right\Vert_{\infty,\infty} \hspace{0cm}+\hspace{0cm} \left\Vert Y^N(t)\right\Vert_{\infty,\infty}$\\

$\displaystyle \hspace{0.2cm}+\int_0^t\left\Vert
\begin{pmatrix}
F\left(\bar{u}^N(s)\right)-F\left(v^N(s)\right)\vspace{0.1cm} \\ G^N\left(\bar{u}^N(s)\right)-G^N\left(v^N(s)\right)
\end{pmatrix}\right\Vert_{\infty,\infty} + \left\Vert
\begin{pmatrix}
0\vspace{0.1cm} \\ G^N\left(v^N(s)\right)-G\left(v^N(s)\right)
\end{pmatrix}\right\Vert_{\infty,\infty}ds$.\\


\noindent We have \\

$\displaystyle \hspace{0.5cm} \left\Vert
\begin{pmatrix}
F\left(\bar{u}^N(s)\right)-F\left(v^N(s)\right)\vspace{0.1cm} \\ G^N\left(\bar{u}^N(s)\right)-G^N\left(v^N(s)\right)
\end{pmatrix}\right\Vert_{\infty,\infty}\leq 2L\left\Vert \bar{u}^N(s)-v^N(s)\right\Vert_{\infty,\infty}$,\\

\noindent and \\

$\disp\hspace{0.5cm} \left\Vert
\begin{pmatrix}
0\vspace{0.1cm} \\ G^N\left(v^N(s)\right)-G\left(v^N(s)\right)
\end{pmatrix}\right\Vert_{\infty,\infty}$\\

$\disp \hspace{2.5cm} =\left\Vert  G^N\left(v^N(s)\right)-G\left(v^N(s)\right)\right\Vert_{\infty}$\\

$\disp \hspace{2.5cm} \leq \left\Vert  G^N\left(v^N(s)\right)-G^N\left(v(s)\right)\right\Vert_{\infty} + \left\Vert  G^N\left(v(s)\right)-G\left(v(s)\right)\right\Vert_{\infty}$\vspace{0.2cm}\par 

$\disp \hspace{3.5cm}+ \left\Vert  G\left(v(s)\right)-G\left(v^N(s)\right)\right\Vert_{\infty}$\\

$\disp \hspace{2.5cm} \leq 2L  \left\Vert v^N(s)-v(s)\right\Vert_{\infty,\infty} +\left\Vert  G^N\left(v(s)\right)-G\left(v(s)\right)\right\Vert_{\infty}$.\\

%
%
%

\noindent where $L$ is a Lipschitz constant common to $F$, $G$ and $G^N$. Therefore, \\

$\displaystyle \left\Vert\bar{u}^N(t)-v^N(t)\right\Vert_{\infty,\infty}$\par 

$\disp \hspace{1cm} \leq\left\Vert u^N(0)-v^N(0)\right\Vert_{\infty,\infty} + \left\Vert Y^N(t)\right\Vert_{\infty,\infty}+2L\int_0^t\left\Vert \bar{u}^N(s)-v^N(s)\right\Vert_{\infty,\infty}ds$\\

$\displaystyle\hspace{2.5cm} + \int_0^t\left(2L\left\Vert v^N(s)-v(s)\right\Vert_{\infty,\infty}+\left\Vert G^N(v(s))-G(v(s))\right\Vert_{\infty}\right)ds$\\


\noindent Taking the supremum in time on $[0,T]$ and applying Gronwall lemma, we obtain\\

\noindent $\displaystyle \sup_{[0,T]}\left\Vert\bar{u}^N(t)-v^N(t)\right\Vert_{\infty,\infty}$\\
\begin{equation}\label{Gronwall_Estimation_2}
\begin{array}{l}
\displaystyle \leq\left\lbrace\left\Vert u^N(0)-v^N(0)\right\Vert_{\infty,\infty} + \sup_{[0,T]}\left\Vert Y^N(t)\right\Vert_{\infty,\infty}   \right. \vspace{0.1cm}\\

\displaystyle \hspace{0.5cm}\left. + 2LT\sup_{[0,T]}\left\Vert v^N(t)-v(t)\right\Vert_{\infty,\infty}\hspace{0cm}+\hspace{0cm}\int_0^{T}\left\Vert G^N(v(s))-G(v(s))\right\Vert_{\infty}ds\right\rbrace
\times\text{e}^{2LT}.
\end{array}
\end{equation}

\noindent Now, by assumption, \[ \left\Vert u^N(0)-v^N(0)\right\Vert_{\infty,\infty}\leq\left\Vert u^N(0)-v(0)\right\Vert_{\infty,\infty}+\left\Vert \tilde{P}_Nv(0)-v(0)\right\Vert_{\infty,\infty} \longrightarrow 0 \]

\noindent in probability.   
Also,  $\disp \hspace{0cm} \int_0^{T}\left\Vert G^N(v(s))-G(v(s))\right\Vert_{\infty}\longrightarrow 0$, thanks to \cref{debit_function_regularity_and_convergence} and dominated convergence. Moreover, \cref{Discrete_Continuous_relation} yields $\disp \sup_{[0,T]}\left\Vert v^N(t)-v(t)\right\Vert_{\infty,\infty}\longrightarrow 0$. \par 

Finally, we claim that:
\begin{equation}\label{Martingale_convergence}
\displaystyle \sup_{[0,T]}\left\Vert Y^N(t)\right\Vert_{\infty,\infty}\longrightarrow 0\hspace{0.3cm}\text{in probability}.
\end{equation}
The proof of \cref{Hybrid_law_of_large_numbers} is completed, if our claim is justified. We write
\[
Y^N(t)=\int_0^t\tilde{T}_N(t-s)d
\begin{pmatrix}
Z_C^N(s\wedge\tau)\vspace{0.2cm} \\ Z_D^N(s\wedge\tau)
\end{pmatrix} = 
\begin{pmatrix}
\int_0^tT_N(t-s)dZ_C^N(s\wedge\tau)\vspace{0.2cm} \\ Z_D^N(t\wedge\tau)
\end{pmatrix} =:
\begin{pmatrix}
Y_C^N(t)\vspace{0.2cm} \\ Y_D^N(t)
\end{pmatrix},
\]
where $\displaystyle Z_C^N(t\wedge\tau)$ and $\displaystyle Z_D^N(t\wedge\tau)$ are $\mathds{H}^N-$valued martingales. It follows that \[ \sup_{[0,T]}\left\Vert Y^N(t)\right\Vert_{\infty,\infty}\leq \sup_{[0,T]}\left\Vert Y_C^N(t)\right\Vert_{\infty}+\sup_{[0,T]}\left\Vert Z_D^N(t\wedge\tau)\right\Vert_{\infty}. \] Therefore, it is sufficient to show that each term on the r.h.s. of the inequality converges to zero in probability. For that purpose, we need the subsequent results, whose proofs are differed to \Cref{proof_related_to_the_Gronwall-Bellman_argument}.

\begin{lemma}\label{Blount_scalar_product_estimate}\textbf{\texttt{(Lemma 4.3, \cite{Blount1992})}}
Set $\displaystyle f=N\mathds{1}_j$. Then \par 
\noindent $\displaystyle \hspace{0cm} \left\langle\left(\nabla_N^+T_N(t)f\right)^2+\left(\nabla_N^-T_N(t)f\right)^2+\left(T_N(t)f\right)^2,1\right\rangle_2\leq h_N(t),\hspace{0cm}\text{ with }\hspace{0cm}\int_0^th_N(s)ds\leq KN+t$. 
\end{lemma}

\begin{lemma}\label{Blount_expectation_estimate} \textbf{\texttt{(Lemma 4.4, \cite{Blount1992})}}
Let $m(t)$ be a bounded martingale of finite variation defined on $[t_0,t_1]$, with $m(t_0)~=~0$, and satisfying:\vspace{0.1cm}\par 

\textbf{(i)} $m$ is right-continuous with left limits.\par 
\textbf{(ii)} $|\delta m(t)|\leq 1$ for $t_0\leq t\leq t_1$. \par 
\textbf{(iii)} $\displaystyle \sum_{t_0\leq s\leq t}[\delta m(s)]^2-\int_{t_0}^tg(s)ds$ is a mean $0$ martingale with $0\leq g(s)\leq h(s)$, where $h(s)$ is a bounded deterministic function and $g(s)$ is $\displaystyle \mathcal{F}_t^N$-adapted.\par 

Then, \[ \mathds{E}\left[e^{m(t_1)}\right]\leq exp\left(\frac{3}{2}\int_{t_0}^{t_1}h(s)ds\right). \]
\end{lemma}

Let us go back to the proof of $\ref{Martingale_convergence}$.  

\subsubsection{Martingale continuous component term}

We want to show that\par 
$\disp\hspace{4.5cm} \mathds{P}\left\lbrace\sup_{[0,T]}\left\Vert Y_C^N(t)\right\Vert_{\infty}>\epsilon_0\right\rbrace\longrightarrow 0$.\\

Fix $\bar{t}\in(0,T]$ and $j\in\lbrace1,\cdots,N\rbrace$. For $0\leq t\leq\bar{t}$, set \[ f:=N\mathds{1}_j,\hspace{0.6cm}\text{and}\hspace{0.4cm} \bar{m}_C(t):=\left\langle\int_0^tT_N(\bar{t}-s)dZ_C^N(s\wedge\tau),f\right\rangle_2. \] 

\noindent Then, $\displaystyle \bar{m}_C:=\left\{\bar{m}_C(t),0\leq t\leq\bar{t}\right\}$  is a mean zero martingale such that \\

$\disp \hspace{1cm}\bar{m}_C(\bar{t})=N\int_0^1\int_0^{\bar{t}}T_N(\bar{t}-s)dZ_C^N\left(s\wedge\tau,x\right)\mathds{1}_j(x)dx$\par 

$\disp\hspace{2.2cm}=\int_0^{\bar{t}}T_N(\bar{t}-s)dZ_C^N\left(s\wedge\tau,\frac{j}{N}\right)\hspace{0.5cm}=\hspace{0.5cm}Y_C^N\left(\bar{t},\frac{j}{N}\right)$.\\

We need the additional upcoming result, that we prove in \Cref{proof_related_to_the_Gronwall-Bellman_argument}.
\begin{lemma}\label{third_level_accompanying_martingale}\textbf{\texttt{(A third type of accompanying martingale)}}
\[
\begin{array}{l}
\displaystyle \sum_{s\leq t\wedge\tau}\left[\delta\bar{m}_C(s)\right]^2\\

\disp\hspace{1cm} -\frac{1}{N\mu}\int_0^{t\wedge\tau}\left[\left\langle u_C^N(s),\left(\nabla_N^+T_N(\bar{t}-s)f\right)^2+\left(\nabla_N^-T_N(\bar{t}-s)f\right)^2\right\rangle_2\right.\\

\displaystyle\hspace{4.5cm}\left.+\left\langle\left(|F|^2+\mu\big|F_1^N\big|^2\right)\left(u^N(s)\right),(T_N(\bar{t}-s)f)^2\right\rangle_2\right] ds
\end{array}
\]
defines a mean 0  \textit{c\`adl\`ag} martingale for $0\leq t\leq\bar{t}$.
\end{lemma}

\begin{remark}
In \cite{Blount1992}, Blount uses this result and claims that it is a consequence of  \cref{Second_level_acccompanying_martingales}, whose proof uses \cref{First_level_acccompanying_martingales}. We have not been able to reproduce his proof. Our proof uses \cref{First_level_acccompanying_martingales} directly.
\end{remark}

\noindent From $(\ref{Discrete_Laplacian_properties})$, we know that $T_N(t)$ is a contraction on $\displaystyle \left(\mathds{H}^N,\Vert\cdot\Vert_\infty\right)$. Thus, $(\ref{Truncated_process_jump_boundedness})$ yields\\

$\displaystyle\hspace{1cm} |\delta\bar{m}_C(t)|=\left|T_N(\bar{t}-t)\delta Z_C^N\left(t\wedge\tau,\frac{j}{N}\right)\right|\leq \left\Vert T_N(\bar{t}-t)\delta Z_C^N\left(t\wedge\tau\right)\right\Vert_{\infty}$\par 

$\displaystyle \hspace{2.5cm}\leq \left\Vert\delta Z_C^N\left(t\wedge\tau\right)\right\Vert_\infty\hspace{0.1cm}=\hspace{0.1cm}\left\Vert\delta u_C^N(t\wedge\tau)\right\Vert_\infty\leq \gamma\mu^{-1}$.\\

\noindent For $\displaystyle \bar{\theta}\in\left[0,\gamma^{-1}\right]$, the process $m_C$ defined by \[ m_C(t):=\bar{\theta}\mu\bar{m}_C(t)\hspace{0.4cm}\text{ for }\hspace{0.5cm} 0\leq t\leq \bar{t} \] is a mean 0 \textit{c\`adl\`ag} martingale such that $|\delta m_C(t)|=\bar{\theta}\mu\left|\delta\bar{m}_C(s)\right|\leq 1$. 
Thus, \[ \sum_{s\leq t\wedge\tau}[\delta m_C(s)]^2-\int_0^{t\wedge\tau} g_C^N(s)ds \] 
\noindent is a mean 0 \textit{c\`adl\`ag} martingale thanks to \cref{third_level_accompanying_martingale}, where for $0\leq s\leq t\wedge\tau$,\\

$\disp \hspace{1cm} g_C^N(s)=\frac{\bar{\theta}^2\mu}{N}\left[\left\langle u_C^N(s),\left(\nabla_N^+T_N(\bar{t}-s)f\right)^2+\left(\nabla_N^-T_N(\bar{t}-s)f\right)^2\right\rangle_2\right.$\par 

$\disp\hspace{5cm}\left.+\left\langle\left(|F|^2+\mu\big|F_1^N\big|^2\right)\left(u^N(s)\right),\left(T_N(\bar{t}-s)f\right)^2\right\rangle_2\right]$.\\

\noindent Since $u_C^N$ is positive and the reaction rates are bounded, $(\ref{Troncated_Markov_sequence_values_boundedness_before})$ and \cref{Blount_scalar_product_estimate} yield \\

$\displaystyle g_C^N(s)\leq \frac{\bar{\theta}^2\mu}{N}\bar{c}\left\langle 1,\left(\nabla_N^+T_N(\bar{t}-s)f\right)^2+\left(\nabla_N^-T_N(\bar{t}-s)f\right)^2+\left(T_N(\bar{t}-s)f\right)^2\right\rangle_2$\par 

$\displaystyle \hspace{1cm}\leq \frac{\bar{\theta}^2\mu}{N}\bar{c}h_1^{N,C}\left(\bar{t}- s\right)$,\\

\noindent where we know from $(\ref{Blount_inner_product_estimate_constant})$ that we can take $\displaystyle h_1^{N,C}(t)=1+4\sum_{m>0}e^{-2\beta_{m,N}\left(\bar{t}-t\right)}\left(\beta_{m,N}+1\right)$, and it satisfies $\displaystyle \int_0^th_1^{N,C}(\bar{t}-s)ds\leq KN+t$, for $0\leq t\leq \bar{t}$. Hence, $g_C^N$ is an $\mathcal{F}_t^N$-adapted process such that  \[ 0\leq g_C^N(s)\leq h_C^N(s), \] where $\displaystyle h_C^N(s)=\bar{c}\frac{\bar{\theta}^2\mu}{N}h_1^{N,C}(\bar{t}-s)$ is a bounded deterministic function on $[0,\bar{t}]$. Since $N\rightarrow\infty$ and $\bar{t}\leq T<\infty$, we may assume $\displaystyle \bar{t}N^{-1}\leq 1$ and get \[ \int_0^{\bar{t}}h_C^N(s)ds\leq \bar{c}\bar{\theta}^2\mu\left(K+\bar{t}N^{-1}\right)\leq \bar{c}\bar{\theta}^2\mu. \] Therefore, \cref{Blount_expectation_estimate} implies $\disp \mathds{E}\left[e^{
m_C(\bar{t})}\right]\leq \text{exp}\left(\bar{c}\bar{\theta}^2\mu\right)$. By Markov's inequality,\\

\noindent $\displaystyle \mathds{P}\left\{Y_C^N\left(\bar{t},\frac{j}{N}\right)>\epsilon_0\right\}=\mathds{P}\left\{\bar{m}_C\left(\bar{t}\right)
>\epsilon_0\right\}
=\mathds{P}\left\{m_C\left(\bar{t}\right)>\bar{\theta}\mu\epsilon_0\right\}= \mathds{P}\left\{e^{m_C(\bar{t})}>e^{\bar{\theta}\mu\epsilon_0}\right\}$\par 

$\displaystyle \hspace{3.5cm}\leq e^{-\bar{\theta}\mu\epsilon_0}\mathds{E}\left[e^{m_C(t)}\right]\leq \text{exp}\left[\bar{\theta}\mu\left(\bar{c}\bar{\theta}-\epsilon_0\right)\right]$.\\

\noindent Thus we can choose $\bar{\theta}$ such that \[ \mathds{P}\left\{Y_C^N\left(\bar{t},\frac{j}{N}\right)>\epsilon_0\right\}\leq e^{-\alpha\epsilon_0^2\mu},\hspace{0.3cm} \text{for}\hspace{0.2cm} \alpha=\bar{c}(\tilde{\rho}_T,T)>0, \] 
independently of $N$, $\mu$, $j$ and $\bar{t}$ \big(one may solve $\bar{c}\bar{\theta}^2-\epsilon_0\bar{\theta}+\alpha\epsilon_0^2\leq 0$ in $\bar{\theta}$\big). 
The same relation holds for $\displaystyle \mathds{P}\left\{-Y_C^N\left(\bar{t},\frac{j}{N}\right)>\epsilon_0\right\}$, by replacing the processes $\bar{m}_C$ and $Y_C^N$ by their opposites $-\bar{m}_C$ and $-Y_C^N$ respectively, and repeating the argument.
Therefore, \[ \mathds{P}\left\{\left|Y_C^N\left(t,\frac{j}{N}\right)\right|>\epsilon_0\right\}\leq 2e^{-\alpha\epsilon_0^2\mu},\hspace{0.5cm}\text{ for}\hspace{0.3cm} 0<t\leq T\hspace{0.2cm}\text{ and }\hspace{0.2cm} j\in\{1,\cdots,N\}.  \]
Since $\displaystyle \left\Vert Y_C^N(t)\right\Vert_\infty=\sup_{1\leq j\leq N}\left|Y_C^N\left(t,\frac{j}{N}\right)\right|$ and $Y_C^N(0)=0$,
\begin{equation}\label{Convergence_of_Y_C_for_fixed_time}
\begin{array}{l}
\displaystyle \mathds{P}\left\{\left\Vert Y_C^N(t)\right\Vert_\infty>\epsilon_0\right\}=\mathds{P}\left\{\exists j\in\{1,\cdots,N\}:\left|Y_C^N\left(t,\frac{j}{N}\right)\right|>\epsilon_0\right\}\vspace{0.2cm}\\

\displaystyle \hspace{3.2cm}\leq \hspace{0.3cm}\sum_{1\leq j\leq N}\mathds{P}\left\{\left|Y_C^N\left(t,\frac{j}{N}\right)\right|>\epsilon_0\right\}\hspace{0.7cm}\leq\hspace{0.7cm} 2Ne^{-\alpha\epsilon_0^2\mu}
\end{array}
\end{equation}
for $0\leq t\leq T$, $\alpha=\alpha(\rho)>0$. We now show that $(\ref{Convergence_of_Y_C_for_fixed_time})$ holds with $\displaystyle \left\Vert Y_C^N(t)\right\Vert_\infty$ replaced by $\sup_{[0,T]}\left\Vert Y_C^N(t)\right\Vert_\infty$ and $N$ - on the r.h.s. - replaced by $N^3$. From Duhamel formula, $\disp Y_C^N(t)=\int_0^tT_N(t-s)dZ_C^N(s\wedge\tau)$ satisfies the Stochastic Differential Equation (SDE) \[ dY_C^N(t)=\Delta_NY_C^N(t)+ dZ_C^N(t\wedge\tau). \] whose integral form is
\begin{equation}\label{Y_representation_by_variation_of_constant}
\displaystyle Y_C^N(t)=Z_C^N(t\wedge\tau)+\int_0^t\Delta_NY_C^N(s)ds.
\end{equation}
\noindent We subdivide $[0,T]$ into $N^2$ subintervals $I_n(T):=\left[\frac{nT}{N^2},\frac{(n+1)T}{N^2}\right]$, $0\leq n\leq N^2-1$. 
Then taking $t=\frac{nT}{N^2}$ in $(\ref{Y_representation_by_variation_of_constant})$ yields \par 

$\displaystyle \hspace{3.7cm} Y_C^N\left(\frac{nT}{N^2}\right)=Z_C^N\left(\frac{nT}{N^2}\wedge\tau\right)+\int_0^{\frac{nT}{N^2}}\Delta_NY_C^N(s)ds$.\\

\noindent One writes \par 

$\displaystyle\hspace{0.8cm} Y_C^N(t)=Z_C^N(t\wedge\tau)+\int_0^{\frac{nT}{N^2}}\Delta_NY_C^N(s)ds+\int_{\frac{nT}{N^2}}^t\Delta_NY_C^N(s)ds$\\

$\displaystyle \hspace{2cm}=Z_C^N(t\wedge\tau)+\left[Y_C^N\left(\frac{nT}{N^2}\right)-Z_C^N\left(\frac{nT}{N^2}\wedge\tau\right)\right]+\int_{\frac{nT}{N^2}}^t\Delta_NY_C^N(s)ds$\\

$\displaystyle \hspace{2cm}=Y_C^N\left(\frac{nT}{N^2}\right)+\int_{\frac{nT}{N^2}}^t\Delta_NY_C^N(s)ds+\tilde{m}_C(t)$,\par 

\noindent where \[ \tilde{m}_C(t)=Z_C^N(t\wedge\tau)-Z_C^N\left(\frac{nT}{N^2}\wedge\tau\right) \] defines a mean 0 martingale, for $t\in I_n(T)$. 
Thus, \[ \left\Vert Y_C^N(t)\right\Vert_\infty\leq\left\Vert Y_C^N\left(\frac{nT}{N^2}\right)\right\Vert_\infty+4N^2\int_{\frac{nT}{N^2}}^t\left\Vert Y_C^N(s)\right\Vert_\infty ds+\left\Vert\tilde{m}_C(t)\right\Vert_\infty \]
since $\Vert\Delta_N\Vert_\infty\leq 4N^2$. 
With the notation $\disp\sup_{t\in I_n(T)}\equiv\sup_{I_n(T)}$, a Gronwall argument leads to
\begin{equation}\label{Gronwall_lemma_for_Y}
\displaystyle \sup_{I_n(T)}\left\Vert Y_C^N(t)\right\Vert_\infty\leq\left(\left\Vert Y_C^N\left(\frac{nT}{N^2}\right)\right\Vert_\infty+\sup_{I_n(T)}\left\Vert\tilde{m}_C(t)\right\Vert_\infty\right)e^{4T}.\hspace{0.5cm}
\end{equation}\par\vspace{0.1cm}

\noindent We now want to proceed as previously, by finding a suitable martingale associated with $\tilde{m}_C$. Fix $j\in\{1,\cdots,N\}$, $\displaystyle \bar{\theta}\in\left[0,\left|\gamma\right|^{-1}\right]$ and set \[ m_C(t)=\bar{\theta}\mu\tilde{m}_C\left(t,\frac{j}{N}\right),\hspace{0.4cm}\text{for}\hspace{0.3cm}t\in I_n(T)=\left[\frac{nT}{N^2},\frac{(n+1)T}{N^2}\right]. \]
For $s\in I_n(T)$, $\displaystyle s\geq \frac{nT}{N^2}$ yields $\displaystyle \left|\delta m_C(s)\right|=\bar{\theta}\mu\left|\delta\tilde{m}_C(s)\right|=\bar{\theta}\mu\left|\delta Z_{j}^{N,C}(s\wedge\tau)\right|\leq 1$. 
By \cref{Second_level_acccompanying_martingales}, \[ \sum_{\frac{nT}{N^2}\wedge\tau\leq s\leq t\wedge\tau}\left[\delta m_C(s)\right]^2-\int_{\frac{nT}{N^2}\wedge\tau}^{t\wedge\tau}g_C^N(s)ds \] is a mean $0$ martingale for $t\in I_n(T)$,
where \[ g_C^N(s)=\bar{\theta}^2\mu\left[N^2\left(u_{j-1}^{N,C}(s)+2u_j^{N,C}(s)+u_{j+1}^{N,C}(s)\right)+\left(|F|_j^2+\mu\big|F_1^N\big|_j^2\right)\big(u^N(s)\big)\right]. \]
Now, $\displaystyle 0\leq g_C^N(s) \leq h_C^N(s)$ by $(\ref{Troncated_Markov_sequence_values_boundedness_before})$ and the local boundedness of the reaction rates, where \par  
$\displaystyle \hspace{1cm} h_C^N(s)=\bar{\theta}^2\mu\left(4N^2\tilde{\rho}_T+\bar{\gamma}_C\bar{\lambda}(\tilde{\rho}_T\right)$\hspace{0.2cm}
satisfies\hspace{0.2cm}$\disp \int_{\frac{nT}{N^2}\wedge\tau}^{t\wedge\tau}h_C^N(s)\leq \bar{c}\bar{\theta}^2\mu T$.\\

\noindent \Cref{Blount_expectation_estimate} then implies \[ \mathds{E}\left[\text{exp}\left\{m_C\left(\frac{(n+1)T}{N^2}\right)\right\}\right]\leq\text{exp}\left(\bar{c}\bar{\theta}^2\mu T\right). \]
Applying Doob's inequalities, we have \\

$\displaystyle\hspace{0cm} \mathds{P}\left\{\sup_{I_n(T)}\tilde{m}_C\left(t,\frac{j}{N}\right)>\epsilon_0\right\}\leq e^{-\bar{\theta}\mu\epsilon_0}\mathds{E}\left[\text{exp}\left\{m_C\left(\frac{(n+1)T}{N^2}\right)\right\}\right]$\par 

$\displaystyle \hspace{5.2cm}\leq \hspace{0.2cm}\text{exp}\left[\bar{\theta}\mu(\bar{c}\bar{\theta}-\epsilon_0)\right]\hspace{0.7cm}\leq\hspace{0.7cm} e^{-\alpha\epsilon_0^2\mu}$,\\

\noindent where $\alpha=\bar{c}(\tilde{\rho}_T,T)>0$, independently of $N$, $\mu$ and $j$. We have choosen a suitable $\bar{\theta}$ as previously. Again, the same inequality holds for $\displaystyle -\tilde{m}_C\left(t,\frac{j}{N}\right)$ and consequently,
\begin{equation}\label{Proba_norme_infinie_m_tilde}
\displaystyle \mathds{P}\left\{\sup_{I_n(T)}\left\Vert\tilde{m}_C(t)\right\Vert_\infty>\epsilon_0\right\}\leq 2Ne^{-\alpha\epsilon_0^2\mu}.
\end{equation}
From $(\ref{Convergence_of_Y_C_for_fixed_time})$, $(\ref{Gronwall_lemma_for_Y})$ and $(\ref{Proba_norme_infinie_m_tilde})$,\\

\noindent $\displaystyle \mathds{P}\left\{e^{-4T}\sup_{I_n(T)}\left\Vert Y_C^N(t)\right\Vert_\infty>\epsilon_0\right\}\leq \mathds{P}\left\{\left\Vert Y_C^N\left(\frac{nT}{N^2}\right)\right\Vert_\infty+\sup_{I_n(T)}\left\Vert\tilde{m}_C(t)\right\Vert_\infty>\epsilon_0\right\}$
\begin{equation}\label{Convergence_of_norm_of_Y_C_for_fixed_time_interval}
\displaystyle \hspace{2cm}\leq \mathds{P}\left\{\left\Vert Y_C^N\left(\frac{nT}{N^2}\right)\right\Vert_\infty>\frac{\epsilon_0}{2}\right\}+\mathds{P}\left\{+\sup_{I_n(T)}\left\Vert\tilde{m}_C(t)\right\Vert_\infty>\frac{\epsilon_0}{2}\right\}
\end{equation}
\noindent $\displaystyle \hspace{3cm} \leq 2Ne^{-\alpha\frac{\epsilon_0^2}{4}\mu}+2Ne^{-\alpha\frac{\epsilon_0^2}{4}\mu}=4Ne^{-\alpha\frac{\epsilon_0^2}{4}\mu} \leq 4Ne^{-\alpha\epsilon_0^2\mu}$.\par 

\noindent Hence, \[ \mathds{P}\left\{e^{-4T}\sup_{[0,T]}\left\Vert Y_C^N(t)\right\Vert_\infty>\epsilon_0\right\}\leq\sum_{n=0}^{N^2-1}\mathds{P}
\left\{e^{-4T}\sup_{I_n(T)}\left\Vert Y_C^N(t)\right\Vert_\infty>\epsilon_0\right\}\leq 4N^3e^{-\alpha\epsilon_0^2\mu}, \]

\noindent where $\alpha=\alpha(\tilde{\rho}_T,T)>0$. We deduce \[ \mathds{P}\left\{\sup_{[0,T]}\left\Vert Y_C^N(t)\right\Vert_\infty>\epsilon_0\right\}\leq 4N^3\text{exp}\left(-\alpha e^{-16T^2}\epsilon_0^2\mu\right)\leq \text{exp}(\log4+3\log N-\alpha\mu), \]
with $\alpha=\alpha(\tilde{\rho}_T,T,\epsilon_0)$. The r.h.s. vanishes, since we are assuming $\displaystyle \mu^{-1}\log N\rightarrow0$. This proves the convergence of the continuous part. $\square$

\subsubsection{Martingale discrete component term}

We want to show that \par 
$\disp\hspace{4.5cm} \mathds{P}\left\lbrace\sup_{[0,T]}\left\Vert Y_D^N(t)\right\Vert_{\infty}>\epsilon_0\right\rbrace\longrightarrow 0$. \\

\noindent We use the same procedure as for the continuous part. Fix $\bar{t}\in(0,T]$ and $j\in\lbrace1,\cdots,N\rbrace$. 
For $0\leq t\leq\bar{t}$, set \[ f:=N\mathds{1}_j,\hspace{0.6cm}\text{and}\hspace{0.4cm} \bar{m}_D(t):=\left\langle Z_D^N(t\wedge\tau),f\right\rangle_2. \] 
Then, the process $\bar{m}_D:=\left\{\bar{m}_D(t),0\leq t\leq\bar{t}\right\}$ defines a mean zero martingale, such that \hspace{0.2cm} $\displaystyle \bar{m}_D(\bar{t})=Y_D^N\left(\bar{t},\frac{j}{N}\right)$. Taking $\varphi=f$ in \cref{Second_level_acccompanying_martingales} yields 
\begin{equation}\label{Martingale_bar_first_discrete}
\displaystyle \sum_{s\leq t\wedge\tau}\left[\delta\bar{m}_D(s)\right]^2-\frac{1}{N}\int_0^{t\wedge\tau}\left\langle\left|G^N\right|^2\left(u^N(s)\right),f^2\right\rangle_2 ds
\end{equation}
is a mean 0 martingale for $0\leq t\leq\bar{t}$. From $(\ref{Truncated_process_jump_boundedness})$,\\

$\displaystyle |\delta\bar{m}_D(t)|=\left|\delta Z_D^N\left(t\wedge\tau,\frac{j}{N}\right)\right|\leq\left\Vert\delta Z_D^N\left(t\wedge\tau\right)\right\Vert_\infty\hspace{0.1cm}=\hspace{0.1cm}\left\Vert\delta u_D^N(t\wedge\tau)\right\Vert_\infty\leq \gamma N^{-1}$.\\

\noindent Since $N\rightarrow\infty$, we may assume $\displaystyle N^{-1}\leq 1$. Then, For $\displaystyle \bar{\theta}\in\left[0,\gamma^{-1}\right]$, \[ m_D(t)=\bar{\theta} N\bar{m}_D(t) \] 
defines a mean $0$ \textit{c\`adl\`ag} martingale such that $\left|\delta m_D(t)\right|\leq 1$ and 
$\displaystyle \left[\delta m_D(s)\right]^2=\bar{\theta}^2N^2\left[\delta\bar{m}_D(s)\right]^2$ for $0\leq t\leq \bar{t}$.
Thus, from \cref{Second_level_acccompanying_martingales}, \[ \sum_{s\leq t}[\delta m_D(s)]^2-\int_0^{t\wedge\tau} g_D^N(s)ds \] 
\noindent is a mean $0$ \textit{c\`adl\`ag} martingale, where\\

$\displaystyle\hspace{1cm} g_D^N(s):=\bar{\theta}^2N\left\langle\left|G^N\right|\left(u^N(s)\right),f^2\right\rangle_2=\bar{\theta}^2N^2\left|G_j^N\right|^2\left(u^N(s)\right)$\par 

$\displaystyle \hspace{2cm}=\bar{\theta}^2N^2\sum_{i=1}^N|\gamma_{ij}^N|^2|g|_{ij}^2\big(u^N(s)\big)\leq \bar{\theta}^2\bar{c}N$, \\

\noindent thanks to the boundedness of the reaction rates. It follows that $\disp 0\leq g_D^N(s)\leq h_D^N(s)$, 
where $\displaystyle h_D^N(s):=\bar{\theta}^2\bar{c}N$ satisfies $\displaystyle \int_0^th_D^N(s)ds\leq \bar{\theta}^2\bar{c}N$ for $0\leq t\leq \bar{t}\leq T$. 
Therefore, \Cref{Blount_expectation_estimate} implies $\disp \mathds{E}\left[e^{m_D(\bar{t})}\right]\leq \text{exp}\left(\bar{\theta}^2\bar{c}N\right)$. 
Markov's inequality yields\\

$\displaystyle \mathds{P}\left\{Y_D^N\left(\bar{t},\frac{j}{N}\right)>\epsilon_0\right\}=\mathds{P}\left\{\bar{m}_D\left(\bar{t}\right)
>\epsilon_0\right\}
=\mathds{P}\left\{m_D\left(\bar{t}\right)>\bar{\theta} N\epsilon_0\right\}=\mathds{P}\left\{e^{m_D(\bar{t})}>e^{\bar{\theta} N\epsilon_0}\right\}$\par 

$\disp\hspace{4cm}\leq e^{-\bar{\theta} N\epsilon_0}\mathds{E}\left[e^{m_D(t)}\right]\leq \text{exp}\big[\bar{\theta} N(\bar{c}\bar{\theta}-\epsilon_0)]$.\\

We can choose $\bar{\theta}$ such that \[ \mathds{P}\left\{Y_D^N\left(\bar{t},\frac{j}{N}\right)>\epsilon_0\right\}\leq e^{-\alpha\epsilon_0^2N},\hspace{0.3cm}\text{for}\hspace{0.2cm} \alpha=\bar{c}(M_1,\rho,T)>0, \] independently of $N$, $\mu$, $j$ and $\bar{t}$. The same relation holds for $\displaystyle \mathds{P}\left\{-Y_D^N\left(\bar{t},\frac{j}{N}\right)>\epsilon_0\right\}$. 
Therefore, \[ \mathds{P}\left\{\left|Y_D^N\left(t,\frac{j}{N}\right)\right|>\epsilon_0\right\}\leq 2e^{-\alpha\epsilon_0^2N}\hspace{0.2cm}\text{for}\hspace{0.2cm} 0<t\leq T\hspace{0.2cm}\text{ and }\hspace{0.2cm} j\in\{1,\cdots,N\}. \]
Since $\displaystyle \left\Vert Y_D^N(t)\right\Vert_\infty=\sup_{1\leq j\leq N}\left|Y_D^N\left(t,\frac{j}{N}\right)\right|$ and $Y_D^N(0)=0$,
\begin{equation}\label{Convergence_of_Y_D_for_fixed_time}
\begin{array}{l}
\displaystyle \mathds{P}\left\{\left\Vert Y_D^N(t)\right\Vert_\infty>\epsilon_0\right\}=\mathds{P}\left\{\exists j\in\{1,\cdots,N\}:\left|Y_D^N\left(t,\frac{j}{N}\right)\right|>\epsilon_0\right\}\\
\displaystyle \hspace{3.2cm}\leq\hspace{0.2cm}\sum_{1\leq j\leq N}\mathds{P}\left\{\left|Y_D^N\left(t,\frac{j}{N}\right)\right|>\epsilon_0\right\}\hspace{0.7cm}\leq\hspace{0.7cm} 2Ne^{-\alpha\epsilon_0^2N}
\end{array}
\end{equation}
for $0\leq t\leq T$. 
We now show that $(\ref{Convergence_of_Y_D_for_fixed_time})$ holds with $\displaystyle \left\Vert Y_D^N(t)\right\Vert_\infty$ replaced by $\displaystyle \sup_{[0,T]}\left\Vert Y_D^N(t)\right\Vert_\infty$ and $N$ (on the right) replaced by $N^3$. 
With the introduced subdivision of $[0,T]$ into $N^2$ subintervals $I_n(T)=\left[\frac{nT}{N^2},\frac{(n+1)T}{N^2}\right]$, $0\leq n\leq N^2-1$, one can always write
\[ \displaystyle Y_D^N(t)=\tilde{m}_D(t)+Z_D^N\left(\frac{nT}{N^2}\wedge\tau\right)=\tilde{m}_D(t)+Y_D^N\left(\frac{nT}{N^2}\right), \]
where $\displaystyle \tilde{m}_D(t):=Z_D^N(t\wedge\tau)-Z_D^N\left(\frac{nT}{N^2}\wedge\tau\right)$ is a mean zero martingale for $t\in I_n(T)$. Thus, 
\begin{equation}\label{Gronwall_lemma_for_Y_D}
\displaystyle \sup_{I_n(T)}\left\Vert Y_D^N(t)\right\Vert_{\infty}\leq \sup_{I_n(T)}\left\Vert\tilde{m}_D(t)\right\Vert_{\infty}+\left\Vert Y_D^N\left(\frac{nT}{N^2}\right)\right\Vert_{\infty}.
\end{equation}

\noindent As previously, we now want to find a suitable martingale associated with $\tilde{m}_D$. Fix $j\in\{1,\cdots,N\}$, $\displaystyle \bar{\theta}\in\left[0,\gamma^{-1}\right]$ and set \[ m_D(t)=\bar{\theta} N\tilde{m}_D\left(t,\frac{j}{N}\right),\hspace{0.4cm}\text{for}\hspace{0.3cm}t\in I_n(T)=\left[\frac{nT}{N^2},\frac{(n+1)T}{N^2}\right]. \]
For $s\in I_n(T)$, $\displaystyle s\geq \frac{nT}{N^2}$ yields $\displaystyle \left|\delta m_D(s)\right|=\bar{\theta} N\left|\delta Z_{j}^{N,D}(s\wedge\tau)\right|\leq 1$. 
From \cref{Second_level_acccompanying_martingales}, \[ \sum_{\frac{nT}{N^2}\leq s\leq t\wedge\tau}\left[\delta m_D(s)\right]^2-\int_{\frac{nT}{N^2}\wedge\tau}^{t\wedge\tau}g_D^N(s)ds \] is a mean 0 martingale, with $\displaystyle g_D^N(s)=\bar{\theta}^2N^2\left|G^N\right|_j^2\left(u^N(s)\right)\leq\bar{\theta}^2\bar{c}N$. 
Since reaction rates are bounded, $0\leq g_D^N(s)\leq h_D^N(s)$, where $h_D^N(s)=\bar{\theta}^2\bar{c}N$ satisfies $\disp \int_{\frac{nT}{N^2}\wedge\tau}^{t\wedge\tau}h_D^N(s)\leq \bar{c}\bar{\theta}^2N$. 
Thus, \Cref{Blount_expectation_estimate} implies $\disp \mathds{E}\left[\text{exp}\left\{m_D\left(\frac{(n+1)T}{N^2}\right)\right\}\right]\leq \text{exp}\big(\bar{c}\bar{\theta}^2N\big)$. 
Doob's inequalities, lead to \\

$\displaystyle\hspace{0cm} \mathds{P}\left\{\sup_{I_n(T)}\tilde{m}_D\left(t,\frac{j}{N}\right)>\epsilon_0\right\}\leq e^{-\bar{\theta}N\epsilon_0}\mathds{E}\left[\text{exp}\left\{m_D\left(\frac{(n+1)T}{N^2}\right)\right\}\right]$\par 

$\displaystyle \hspace{5.2cm}\leq \hspace{0.2cm}\text{exp}\left[\bar{\theta}N(\bar{c}\bar{\theta}-\epsilon_0)\right]\hspace{0.7cm}\leq\hspace{0.7cm} e^{-\alpha\epsilon_0^2N}$,\\

\noindent where $\alpha=\bar{c}(\tilde{\rho}_T,T)>0$, independently of $N$, $\mu$ and $j$. We have choosen a suitable $\bar{\theta}$ as previously.
Also, the same holds for $\displaystyle -\tilde{m}_D\left(t,\frac{j}{N}\right)$. This shows that
\begin{equation}\label{Proba_norme_infinie_m_tilde_D}
\displaystyle \mathds{P}\left\{\sup_{I_n(T)}\left\Vert\tilde{m}_D(t)\right\Vert_\infty>\epsilon_0\right\}\leq 2Ne^{-\alpha \epsilon_0^2N}.
\end{equation}
From $(\ref{Convergence_of_Y_D_for_fixed_time})$, $(\ref{Gronwall_lemma_for_Y_D})$ and $(\ref{Proba_norme_infinie_m_tilde_D})$, we have\\

\noindent $\displaystyle\hspace{0.2cm} \mathds{P}\left\{e^{-4T}\sup_{I_n(T)}\left\Vert Y_D^N(t)\right\Vert_\infty>\epsilon_0\right\}\leq \mathds{P}\left\{\left\Vert Y_D^N\left(\frac{nT}{N^2}\right)\right\Vert_\infty+\sup_{I_n(T)}\left\Vert\tilde{m}_D(t)\right\Vert_\infty>\epsilon_0\right\}$
\begin{equation}\label{Convergence_of_norm_of_Y_D_for_fixed_time_interval}
\displaystyle \hspace{3cm}\leq \mathds{P}\left\{\left\Vert Y_D^N\left(\frac{nT}{N^2}\right)\right\Vert_\infty>\frac{\epsilon_0}{2}\right\}+\mathds{P}\left\{+\sup_{I_n(T)}\left\Vert\tilde{m}_D(t)\right\Vert_\infty>\frac{\epsilon_0}{2}\right\}
\end{equation}
$\displaystyle \hspace{3.5cm} \leq 2Ne^{-\alpha\frac{\epsilon_0^2}{2}N}+2Ne^{-\alpha\frac{\epsilon_0^2}{2}N}\leq 4Ne^{-\alpha\epsilon_0^2N}$.\par 

\noindent Hence, \[ \mathds{P}\left\{e^{-4T}\sup_{[0,T]}\left\Vert Y_D^N(t)\right\Vert_\infty>\epsilon_0\right\}\leq\sum_{n=0}^{N^2-1}\mathds{P}
\left\{e^{-4T}\sup_{I_n(T)}\left\Vert Y_D^N(t)\right\Vert_\infty>\epsilon_0\right\}\leq 4N^3e^{-\alpha\epsilon_0^2N}. \]

\noindent Setting $c:=\alpha e^{-4T}\epsilon_0>0$, it follows that

\[ \displaystyle \mathds{P}\left\{\sup_{[0,T]}\left\Vert Y_D^N(t)\right\Vert_\infty>\epsilon_0\right\}\leq 4N^3e^{-cN}\rightarrow0 \hspace{0.3cm}\text{ as } \hspace{0.4cm}N\rightarrow\infty, \] which shows the convergence for the discrete part. This completes the proof of $(\ref{Martingale_convergence})$ and ends the proof of \cref{Hybrid_law_of_large_numbers}. $\blacksquare$

\begin{remark}\label{r3.3}
In this proof we have used in a crucial way that at each point, the jump of the discrete component is of order $\frac1N$. For the original model, without spatial 
correlation of the slow reactions, this does not hold since the jump of discrete components are of order $1$. 
\end{remark}


\noindent \textbf{\textit{Concluding remarks.}} \par 

The law of large numbers we present in this work differs from the existing one, essentially with the additional presence of discrete species, which leads to multiscaling. That consideration complements classical one scale stochastic spatial models, which were already very usefull, to a more reallistic model. 

We believe the multiscale model can be extended to higher dimensions both in the space and the number of species, at the prize of very cumbersome notations. That will be one of the directions for future works.

The possibility of a spatial correlation for slow reactions, as considered, supposes that the discrete species is present enough to react sufficiently frequently and get homogenized as the sites get closer.

If the initial model \hyperlink{(M6)}{(M6)} is a typical hybrid stochastic system, its limit looks like one, but this latter has no discrete variable anymore, and jumps totally disappear. In furthcoming works, we will study the fluctuations of the model around its deterministic limit, in order to get the speeds of the convergence of the initial fast and slow dynamics, through a central limit theorem. 

Also, we will present other multiscale stochastic spatial models, with a stronger stochasticity that does not vanishes in the approximation at the first order. The limit will remain a stochastic hybrid system.


\section{Appendix}

\subsection{Formal limit of the generator}\label{formal_limit_of_the_generator}

We consider the sequence  $\cA^N=\cA^{N,\mu}$ of generators defined by $(\ref{Hybrid_Stochastic_Generator})$. We consider a 
$L^2(I)$ framework and consider test functions $\varphi\in C_b^2\big(L^2(I)\times L^2(I)\big)$ . 
Each line is expanded at order two thanks to Taylor expansion. The first order expansion has already been examined to understand the limiting generator.
For $k=1,\cdots,5$, denote by $T_k(N)$ the second order term corresponding to the $k$-th line in the expression of the generator. Let $u=(u_C,u_D)\in C_p(I)\times C_p(I)$ be fixed and let $N,\mu\rightarrow\infty$.\\
We assume that the reaction rates are bounded.\\

\noindent \textbf{At order 1.} 
We already know that the whole terms together read\\

$\disp \cA_1^N\varphi\big(\tilde{P}_Nu\big)=\left\langle D^{1,0}\varphi\big(\tilde{P}_Nu\big),\Delta_NP_Nu_C+F\big(\tilde{P}_Nu\big)+F_1^N\big(\tilde{P}_Nu\big)\right\rangle_2$\par 
$\disp \hspace{8cm}+\left\langle D^{0,1}\varphi\big(\tilde{P}_Nu\big),G^N\big(\tilde{P}_Nu\big)\right\rangle_2$.\\

\noindent For $u_C\in C^3(I)$, it converges to \[ \cA^\infty\varphi(u)=\left\langle D^{1,0}\varphi(u), \Delta u_C+F(u)\right\rangle_2 +\left\langle D^{0,1}\varphi(u), G(u)\right\rangle_2, \] requiring $G^N\big(\tilde{P}_Nu\big)\rightarrow G(u)$, which follows from \cref{Piecewise_Approximation_via_Projection}, \cref{debit_function_regularity_and_convergence} and \\

$\disp\hspace{1cm} \left\Vert G^N\big(\tilde{P}_Nu\big)-G(u)\right\Vert_{\infty}\leq \left\Vert G^N\big(\tilde{P}_Nu\big)-G^N(u)\right\Vert_{\infty}+\left\Vert G^N(u)- G(u)\right\Vert_{\infty}$\par 

$\disp \hspace{5cm}\leq L_0\Vert P_Nu-u\Vert_{\infty,\infty}+\left\Vert G^N(u)- G(u)\right\Vert_{\infty}$,\\

\noindent where we use the same notations as in \cref{proof_of_debit_function_regularity_and_convergence}.\\

\noindent \textbf{At order 2.}  We introduce $e_j:=\sqrt{N}\1_j$, for $1\leq j\leq N$. Since the functions $\disp (\1_j)_{1\leq j\leq N}$ are pairwise orthogonal in $L^2$ and of norm $\Vert \1_j\Vert_2=\frac{1}{\sqrt{N}}$, the family $\big\{e_j,1\leq j\leq N\big\}$ forms an orthonormal basis of $\disp \big(\H^N,\Vert\cdot\Vert_2\big)$. Also, $c$ denotes a generic constant and $\Vert D^{l,k}\varphi\Vert_\infty$ denotes a uniform bound of $D^{l,k}\varphi$, for $0\leq l,k\leq 2$, $l+k\leq2$. Let us treat the term appearing from the second line,\\

$\displaystyle \hspace{1cm} T_2(N)=\frac{1}{2}\sum_{r\in S_1}\sum_{j=1}^ND^{2,0}\varphi\big(\tilde{P}_Nu\big)
\cdot\left\langle\frac{\gamma_r^C}{\mu}\mathds{1}_j,\frac{\gamma_r^C}{\mu}\mathds{1}_j\right\rangle_2\mu\lambda_r\left(u_j^C,u_j^D\right)$\par

$\disp \hspace{2.2cm}\leq \frac{\bar{\gamma}^2\bar{\lambda}}{2}\left\Vert D^{2,0}\varphi\right\Vert_\infty\frac{1}{N\mu}\sum_{j=1}^N\big\Vert e_j\big\Vert_2^2\hspace{0.2cm}\leq \hspace{0.2cm} \frac{c}{\mu}\longrightarrow0$.\\ 

\noindent The first term is similar and  also vanishes at the limit, using the same argument. Next, the third term reads\\

$\displaystyle\hspace{0.5cm} T_3(N)=\frac{1}{2}\sum_{r\in\mathfrak{R}_{DC}\backslash S_1}\sum_{j=1}^ND^{2,0}\varphi\big(\tilde{P}_Nu\big) \cdot\left\langle\frac{\gamma_r^C}{\mu}\1_j,\frac{\gamma_r^C}{\mu}\1_j)\right\rangle_2\lambda_r\left(u_j\right)$\par 

$\disp \hspace{2.5cm}+\frac{1}{2}\sum_{r\in\mathfrak{R}_{DC}\backslash S_1}\sum_{j=1}^ND^{1,1}\varphi\big(\tilde{P}_Nu\big)\cdot\left\langle\frac{\gamma_r^C}{\mu}\1_j,\gamma_{j,r}^{N,D}(\tilde{P}_Nu)\right\rangle_2\lambda_r\left(u_j\right)$\par 

$\displaystyle \hspace{2.5cm}+\frac{1}{2}\sum_{r\in\mathfrak{R}_{DC}\backslash S_1}\sum_{j=1}^ND^{0,2}\varphi\big(\tilde{P}_Nu\big) \cdot\left\langle\gamma_{j,r}^{N,D}(P_Nu_D),\gamma_{j,r}^{N,D}(\tilde{P}_Nu)\right\rangle_2\lambda_r\left(u_j\right)$.\\

\noindent We have\par 

$\disp\hspace{1cm} T_{31}(N):=\frac{1}{2N}\sum_{r\in\mathfrak{R}_{DC}\backslash S_1}\sum_{j=1}^ND^{2,0}\varphi\big(\tilde{P}_Nu\big) \cdot\left\langle\frac{\gamma_r^C}{\mu}e_j,\frac{\gamma_r^C}{\mu}e_j\right\rangle_2\lambda_r\left(u_j\right)$\par 

$\disp \hspace{2.4cm}\leq c\left\Vert D^{2,0}\varphi\right\Vert_\infty\frac{1}{N\mu^2}\sum_{j=1}^N\big\Vert e_j\big\Vert_2^2\hspace{0.2cm}\leq\hspace{0.2cm}\frac{c}{\mu^2}\longrightarrow0$.\\

\noindent Observing that \[ \left\Vert \gamma_{j,r}^{N,D}(\tilde{P}_Nu)\right\Vert_2\leq \left|\gamma_r^D\right|\sum_{i=1}^N|\gamma_{ij}(N)|\Vert\1_i\Vert_2\leq\left|\gamma_r^D\right|\frac{a(0)}{\sqrt{N}} \] for all $1\leq j\leq N$, we deduce from the inequality of Cauchy-Schwarz\\ 

$\disp T_{32}(N):=\frac{1}{2\sqrt{N}}\sum_{r\in\mathfrak{R}_{DC}\backslash S_1}\sum_{j=1}^ND^{1,1}\varphi\big(\tilde{P}_Nu\big)
\cdot\left\langle\frac{\gamma_r^C}{\mu}e_j,\gamma_{j,r}^{N,D}(\tilde{P}_Nu)\right\rangle_2\lambda_r\left(u_j\right)$\par 

$\disp \hspace{1.4cm}\leq c\Vert D^{1,1}\varphi\Vert_\infty\frac{1}{\sqrt{N}\mu}\sum_{j=1}^N\big\Vert e_j\big\Vert_2\left\Vert\gamma_{j,r}^{N,D}(\tilde{P}_Nu)\right\Vert_2\hspace{0.2cm}\leq\hspace{0.2cm}\frac{c}{\sqrt{N}\mu}\longrightarrow0$,\\

\noindent on the one hand, and, on the other hand\\

$\disp T_{33}(N):=\frac{1}{2}\sum_{r\in\mathfrak{R}_{DC}\backslash S_1}\sum_{j=1}^ND^{0,2}\varphi\big(\tilde{P}_Nu\big) \cdot\left\langle\gamma_{j,r}^{N,D}(\tilde{P}_Nu),\gamma_{j,r}^{N,D}(\tilde{P}_Nu)\right\rangle_2\lambda_r\left(u_j\right)$\par 


\noindent Recall that, setting $\tilde P_Nu=u^N$,
\[\hspace{1cm} \gamma_{j,r}^{N,D}\big(u_C^N,u_D^N\big):=\gamma_r^{D}\sum_{i=1}^N\gamma_{ij}^N\mathds{1}_i\theta_{ij}^r(u_i^{N,C},u_i^{N,D}), \] 
where $\disp \gamma_{ij}^N:=\int_{I_i}a\left(x-\frac{j}{N}\right)dx$, and $\theta_{ij}^r(u_i^{N,C},u_i^{N,D})=\theta(u^{N,C}_i+
\frac{\gamma_r^C}{\mu} \gamma_{ij}^N) \theta(u^{N,D}_i+\gamma_r^D \gamma_{ij}^N)$. 
Define the operator $A_r$ by:
$$
A_re_j=\sqrt{N}\sqrt{\lambda_r(u_j^N)}\gamma_{j,r}^{N,D}=\gamma_r^{D}\sum_{i=1}^N\gamma_{ij}^N \theta_{ij}^r(u_i^{N})\sqrt{\lambda_r(u_j^N)}e_i.
$$
Then, using standard property of the trace,
$$
\begin{array}{ll}
\disp T_{33}(N)&\disp =\frac{1}{2N}\sum_{r\in\mathfrak{R}_{DC}\backslash S_1}\sum_{j=1}^ND^{0,2}\varphi(u^N) \cdot\left\langle A_re_j,A_re_j\right\rangle_2\\
&\disp =\frac1{2N} \sum_{r\in\mathfrak{R}_{DC}\backslash S_1} Tr\left(A_r^*D^{0,2}\varphi(u^N)A_r\right)\\
&\disp \le \frac1{2N}  \|D^{0,2}\varphi(u^N)\|_\infty \sum_{r\in\mathfrak{R}_{DC}\backslash S_1}Tr\left(A_r^*A_r\right).
\end{array}
$$
Also
$$
Tr\left(A_r^*A_r\right)=\sum_{j=1}^N \|A_re_j\|_2^2\le \bar \lambda(\tilde \rho_T) {\bar \gamma}^2  \sum_{i,j}^N |\gamma_{ij}^N|^2
$$
and since $\disp |\gamma_{ij}^N|^2=\int_{I_i} a\left(x-\frac{j}N\right)dx = \frac1{\sqrt{N}} \left\langle a\left(\cdot-\frac{j}N\right),e_i\right\rangle_2^2$ we have by Parseval identity:
$$
 \sum_{i}^N |\gamma_{ij}^N|^2=\frac1N\left\|a\left(\cdot-\frac{j}N\right)\right\|_2^2=\frac1N\|a\|_2^2.
$$
It follows that
$$
Tr\left(A_r^*A_r\right)\le \bar \lambda {\bar \gamma}^2  \|a\|_2^2,
$$

\noindent and as a result, $T_3(N)\longrightarrow0$. The fourth term also vanishes at the limit, using the same argument. 

\begin{remark}\label{r4.1}
Again, we see why we had to consider some spatial correlation in our model. Indeed, the case without spatial correlation for the discrete spiecies
corresponds to $\gamma_{ij}^N=1$ for $i=j$ and is $0$ otherwise, obtained with $a$ being the Dirac mass. In this case $Tr\left(A_r^*A_r\right)$ is of order $N$ and
$T_3(N)$ does not converge to $0$. 
\end{remark}

Let us consider the fifth. From \[ \Vert\1_{j+1}-\1_j\Vert_2^2=\frac{2}{N}=\Vert\1_{j-1}-\1_j\Vert_2^2, \] we have\par 

$\disp \hspace{0.5cm} T_5(N)=\frac{1}{2}\sum_{j=1}^N\left[D^{2,0}\varphi\big(\tilde{P}_Nu\big)\cdot\left\langle\frac{\mathds{1}_{j+1}-\mathds{1}_{j}}{\mu},\frac{\mathds{1}_{j+1}-\mathds{1}_{j}}{\mu}\right\rangle_2\right.$\par 

$\disp \hspace{4.5cm}+\left.D^{2,0}\varphi\big(\tilde{P}_Nu\big)\cdot\left\langle\frac{\mathds{1}_{j-1}-\mathds{1}_{j}}{\mu},\frac{\mathds{1}_{j-1}-\mathds{1}_{j}}{\mu}\right\rangle_2 \right]\mu N^2u_j^C$\par 

$\disp \hspace{1.8cm}\leq \frac{\big\Vert u_C\big\Vert_\infty}{2}\left\Vert D^2\varphi\right\Vert_\infty\frac{N^2}{\mu}\sum_{j=1}^N\left(\Vert\1_{j+1}-\1_j\Vert_2^2+\Vert\1_{j-1}-\1_j\Vert_2^2\right)$\par 

$\disp \hspace{1.8cm}\leq c\mu^{-1}N^2\longrightarrow0$,\\

\noindent if $\disp\mu^{-1}N^2\rightarrow0$ as $N,\mu\rightarrow\infty$. 

 Hence, under the additional condition $\disp\mu^{-1}N^2\rightarrow0$, all the terms of the second order in the Taylor expansion vanish at the limit and hence, $\cA^\infty$ is indeed formally the limit of $\cA^N$. $\square$

\subsection{Relation between the limit and its discretization}\label{Discrete_Continuous_relation_proof}

\noindent \underline{Proof of \cref{Discrete_Continuous_relation}}. Let $T>0$ be fixed. Since the operator $\tilde{\Delta}_N$ is linear, it is Lipschitz. Next, the vector field $R=(F,G)$ are locally Lipschitz continuous. Therefore, the initial value problem $(\ref{discrete_limiting_problem_compact_form})$ has a unique local solution $v^N$ thanks to the Picard-Lindel\"of theorem. The bounds \eqref{discretization_boundedness} are proved thanks to the discrete maximum principle
and Gronwall lemma as in section \ref{truncation} and we deduce that $v^N$ is in fact a global solution.

Thanks to \eqref{discretization_boundedness} and \eqref{limit_boundary}, we may assume that $F$ and $G$ are globally Lipschitz and we choose 
$L$ such that 
$$
\begin{array}{l}
\|F(u^1_C,u^1_D)-F(u^2_C,u^2_D)\|_\infty \le L \|(u^1_C,u^1_D)-(u^2_C,u^2_D)\|_{\infty,\infty},\\
\\
\|G(u^1_C,u^1_D)-G(u^2_C,u^2_D)\|_\infty \le L \|(u^1_C,u^1_D)-(u^2_C,u^2_D)\|_{\infty,\infty},
\end{array}
$$
provided $\Vert u_C\Vert_{\infty}\leq\rho_C,\;  \Vert u_D\Vert_{\infty}\leq(\rho_D+1)e^{a(0)M_1(\rho_C)T}$.


we have\\

$\displaystyle v^N(t)-u(t)=\tilde{T}_N(t)\tilde{P}_Nu(0)-\tilde{T}(t)u(0)$\par 

$\disp \hspace{3.3cm}+\int_0^t\tilde{T}_N(t-s)
\begin{pmatrix}
F\left(v^N(s)\right) \vspace{0.1cm}\\ G\left(v^N(s)\right)
\end{pmatrix}-\tilde{T}(t-s)
\begin{pmatrix}
F(u(s)) \vspace{0.1cm}\\ G(u(s))
\end{pmatrix}ds.
$\\

\noindent for all $t\geq0$. Then, \\ 

$\displaystyle\hspace{0.5cm} \left\Vert v^N(t)-u(t)\right\Vert_{\infty,\infty}$\vspace{0.1cm}\par  

$\disp \hspace{1.5cm}\leq \left\Vert\tilde{T}_N(t)\tilde{P}_Nu(0)-\tilde{T}(t)u(0)\right\Vert_{\infty,\infty}$\par 

$\displaystyle \hspace{2cm}+\int_0^t\left\Vert\tilde{T}_N(t-s)\left(
\begin{pmatrix}
F\left(v^N(s)\right) \vspace{0.1cm}\\ G\left(v^N(s)\right)
\end{pmatrix} - \tilde{P}_N 
\begin{pmatrix}
F\left(v(s)\right) \vspace{0.1cm}\\ G\left(v(s)\right)
\end{pmatrix}\right)\right\Vert_{\infty,\infty}ds$\\

$\displaystyle \hspace{2.5cm}+\int_0^t\left\Vert\tilde{T}_N(t-s)\tilde{P}_N
\begin{pmatrix}
F\left(u(s)\right) \vspace{0.1cm}\\ G\left(u(s)\right)
\end{pmatrix}-\tilde{T}(t-s)
\begin{pmatrix}
F\left(u(s)\right) \vspace{0.1cm}\\ G\left(u(s)\right)
\end{pmatrix}\right\Vert_{\infty,\infty}ds$.\\

\noindent 
Observing that $\H^N$ is stable by $F$ and $G$, we rely on $\cref{Piecewise_Approximation_via_Projection}$ and $\cref{Discrete_Laplacian_properties}$ to find an upper bound for each term on the right hand side of the inequality above.
\[ T_2(N)\leq \int_0^t\left\Vert \tilde{P}_N
\begin{pmatrix}
F\left(v^N(s)\right)-F\left(v(s)\right) \vspace{0.1cm}\\ G\left(v^N(s)\right)-G\left(v(s)\right)
\end{pmatrix}\right\Vert_{\infty,\infty}ds \leq 2L\int_0^t\left\Vert v^N(s)-v(s)\right\Vert_{\infty,\infty}ds, \] and\par 
$\displaystyle\hspace{1.5cm} T_3(N) = \int_0^t\left\Vert T_N(t-s)P_NF(u(s))-T(t-s)F(u(s))\right\Vert_\infty ds$.\par 

\noindent Hence,
\begin{align*}
\left\Vert v^N(t)-v(t)\right\Vert_{\infty,\infty} & \leq \left\Vert\tilde{T}_N(t)\tilde{P}_Nv(0)-\tilde{T}(t)v(0)\right\Vert_{\infty,\infty}\\
& \hspace{1cm}+\int_0^t\left\Vert T_N(t-s)P_NF(v(s))-T(t-s)F(v(s))\right\Vert_\infty ds\\
& \hspace{0.2cm}+2L\int_0^t\left\Vert v^N(s)-v(s)\right\Vert_{\infty,\infty}ds.
\end{align*}

\noindent Taking the supremum in $t$ on $[0,T]$ and using Gronwall lemma leads to\\

\noindent $\displaystyle \sup_{[0,T]}\left\Vert v^N(t)-v(t)\right\Vert_{\infty,\infty}$\par 
$\displaystyle\hspace{0.3cm}\leq \left(\sup_{[0,T]}\left\Vert\tilde{T}_N(t)\tilde{P}_Nv(0)-\tilde{T}(t)v(0)\right\Vert_{\infty,\infty}\right.$\par 

$\displaystyle \hspace{1cm} \left.+\int_0^{T}\sup_{t\in[0,T]}\left(\left\Vert T_N(t-s)P_NF(v(s))-T(t-s)F(v(s))\right\Vert_\infty\mathds{1}_{(s\leq t)}\right)ds\right)\times \text{e}^{2LT}$.\\

\noindent Firstly, \vspace{0.1cm}\par 

$\displaystyle \hspace{0.4cm} \sup_{[0,T]}\left\Vert\tilde{T}_N(t)\tilde{P}_Nv(0)-\tilde{T}(t)v(0)\right\Vert_{\infty,\infty}$\par 

$\displaystyle \hspace{1.7cm}\leq \sup_{[0,T]}\left\Vert T_N(t)P_Nu_C(0)-T(t)u_C(0)\right\Vert_{\infty}+\left\Vert P_Nu_D(0)-u_D(0)\right\Vert_\infty\longrightarrow0$,\\

\noindent 
where the first term vanishes thanks to $\cref{Discrete_Laplacian_properties}$ (vi), since $u_C(0)\in C^3(I)$. Secondly, we fix $t\in [0,T]$ and $s\in [0,t]$. Since $F(u(s))\in C^3(I)$, we have \[ \left\Vert T_N(t-s)P_NF(u(s))-T(t-s)F(u(s))\right\Vert_\infty \longrightarrow0, \] thanks to $\cref{Discrete_Laplacian_properties}$ (vi) again. Moreover, the convergence is uniform with respect to $t$, on $[0,T]$. 
Thus, \[ \sup_{t\in[0,T]}\left(\left\Vert T_N(t-s)P_NF(u(s))-T(t-s)F(u(s))\right\Vert_\infty\mathds{1}_{(s\leq t)}\right)\longrightarrow0. \] 
Hence,\hspace{0.2cm} $\disp \int_0^{T}\sup_{t\in[0,T]}\left(\left\Vert T_N(t-s)P_NF(u(s))-T(t-s)F(u(s))\right\Vert_\infty\mathds{1}_{(s\leq t)}\right)ds \longrightarrow0$,\par 
\noindent thanks to the dominated convergence theorem. $\square$

\subsection{The accompanying martingales}


\noindent \underline{Proof of \cref{First_level_acccompanying_martingales}}.\label{proof_about_accompanying_martingales_1} Fix $1\leq j\leq N$ and define the process $\displaystyle Y^N(t)=\left(Y_1^N(t),\cdots,Y_{N+4}^N(t)\right)$ by 
\[ 
\left\{
\begin{array}{l}
\displaystyle Y_k^N(t):=u_k^{N}(t)\hspace{1cm}\text{for}\hspace{0.5cm}1\leq k\leq N,\vspace{0.1cm}\\

\displaystyle Y_{N+1}^N(t):=\sum_{s\leq t}\left[\delta u_{j-1}^{N,C}(s)\right]\left[\delta u_{j}^{N,C}(s)\right]\vspace{0.1cm}\\

\displaystyle Y_{N+2}^N(t):=\sum_{s\leq t}\left[\delta u_{j}^{N,C}(s)\right]^2\vspace{0.1cm}\\

\displaystyle Y_{N+3}^N(t):=\sum_{s\leq t}\left[\delta u_{j}^{N,C}(s)\right]\left[\delta u_{j+1}^{N,C}(s)\right]\vspace{0.1cm}\\

\displaystyle Y_{N+4}^N(t):=\sum_{s\leq t}\left[\delta u_{j}^{N,D}(s)\right]^2.
\end{array}
\right.
\]
$Y^N$ is a Markov process with values in $\displaystyle \mathds{R}^{2N+4}$. From $(\ref{transitions_for_the_spatial_model})$, Assumptions $\ref{Partitioning_Mixed_Reactions_subset}$ and $\ref{spatial_correlation_for_discrete_variable}$, we know that $Y^N$ enjoys the following transitions:\\

\noindent $\displaystyle Y_{N+1}^N\longrightarrow Y_{N+1}^N-\frac{1}{\mu^2}$ at rate $\displaystyle \mu N^2\left(u_{j-1}^{N,C}+u_j^{N,C}\right)$ so that the corresponding debit is \[ -\frac{N^2}{\mu}\left(u_{j-1}^{C}+u_j^{C}\right)=:\left|\Psi_{-}^{N,C}\right|_j^2(u). \]

\noindent Similarly, $\displaystyle Y_{N+3}^N\longrightarrow Y_{N+3}^N-\frac{1}{\mu^2}$ at rate $\displaystyle \mu N^2\left(u_{j}^{N,C}+u_{j+1}^{N,C}\right)$, and the associated debit is \[ -\frac{N^2}{\mu}\left(u_{j}^{C}+u_{j+1}^{C}\right)=:\left|\Psi_{+}^{N,C}\right|^2(u_j). \]
\noindent Now,\par 
$\disp\hspace{0.5cm} Y_{N+2}^N\rightarrow Y_{N+2}^N+\frac{\left|\gamma_r^C\right|^2}{\mu^2} \hspace{0.1cm}\text{ at rate }\hspace{0.1cm} \left\{
\begin{array}{l}
\disp \mu\lambda_r\left(u_j^{N,C}\right)\hspace{1.5cm} \text{for } r\in\mathfrak{R}_C\vspace{0.1cm}\\
\disp \mu\lambda_r\left(u_j^{N,C},u_j^{N,D}\right)\hspace{0.5cm}\text{for } r\in S_1,
\end{array} \right. $\par 

$\disp\hspace{0.5cm} Y_{N+2}^N\rightarrow Y_{N+2}^N+\frac{\left|\gamma_r^C\right|^2}{\mu^2} \hspace{0.1cm}\text{ at rate }\hspace{0.1cm} \lambda_r\left(u_j^{N,C},u_j^{N,D}\right)$ for $ r\in \fR_{DC}\backslash S_1$\vspace{0.1cm}\par 

$\disp\hspace{0.5cm} Y_{N+2}^N\rightarrow Y_{N+2}^N+\frac{1}{\mu^2} \hspace{0.1cm}\text{ at rate }\hspace{0.1cm}\mu N^2\left(u_{j-1}^{N,C}+2u_{j}^{N,C}+u_{j+1}^{N,C}\right) \text{ for a diffusion}$.\vspace{0.2cm}\par 

\noindent The joined debit is\\

$\disp\hspace{0cm} \frac{1}{\mu}\left[N^2\left(u_{j-1}^{C}+2u_j^{C}+u_{j+1}^{C}\right)+\sum_{r\in\mathfrak{R}_C}\left|\gamma_r^C\right|^2\lambda_r\left(u_j^{C}\right)+\sum_{r\in S_1}\left|\gamma_r^C\right|^2\lambda_r\left(u_j^{C},u_j^{D}\right)\right.$\par 
$\disp \hspace{2.5cm}+\frac{1}{\mu}\sum_{i=1}^N\left.\sum_{r\in\fR_{DC}\backslash S_1}|\gamma_{ji}^N|^2|\theta_{ji}^r(u_j^{C},u_j^D)|^2\lambda_r(u_j^C,u_j^D)\right]$\vspace{0.1cm}\par 

$\displaystyle \hspace{0.5cm}=\frac{1}{\mu}\left(\left|\Delta_{N}\right|^2\left(u_j^{C}\right)+|F|_j^2(u_{C},u_{D})+\left|F_1^N\right|_j^2(u_{C},u_{D})\right)\hspace{0.2cm}=\hspace{0.2cm}\frac{1}{\mu}\left|\Psi_C^{N}\right|_j^2(u_C,u_D)$.\\

\noindent Finally, for $i=1,\cdots,N$, $Y_{N+4}^N\longrightarrow Y_{N+4}^N+\left|\gamma_{ij}(N)\right|^2\left|\gamma_r^D\right|^2|\theta_{ij}^r(u_i^C,u_i^D)|^2$ at rate 
\[ \left\{
\begin{array}{l}
\disp \lambda_r\left(u_j^{N,C},u_j^{N,D}\right)\hspace{0.5cm} \text{for } r\in\mathfrak{R}_{DC}\backslash S_1\vspace{0.1cm}\\
\disp \lambda_r\left(u_j^{N,D}\right)\hspace{1.5cm} \text{for } r\in\mathfrak{R}_D.
\end{array}
\right. \]
Theses transitions admit as debit function \[ \big|\Psi_D^{N}\big|_j^2(u_C,u_D)=\sum_{i=1}^N|\gamma_{ij}^N|^2|g|_{ij}^2(u_C,u_D)=:\big|G^N\big|_j^2(u_C,u_D). \]


 
\noindent Furthermore,\\ 

$\displaystyle \sup_{[0,T]}\sup_{1\leq k\leq N}\left|Y_k^N(t\wedge\tau)\right|=\sup_{[0,T]}\sup_{1\leq k\leq N}\left|u_k^N(t\wedge\tau)\right|$\\

$\displaystyle \hspace{4.2cm}\leq c\sup_{[0,T]}\max\left(\left\Vert u_C^{N}(t\wedge\tau)\right\Vert_{\infty},\left\Vert u_D^{N}(t\wedge\tau)\right\Vert_{\infty}\right)$\\

$\displaystyle \hspace{4.2cm}\leq c\sup_{[0,T]}\left\Vert u^{N}(t\wedge\tau)\right\Vert_{\infty,\infty}\hspace{1cm}\leq\hspace{1cm}C(T_N,\mu)$\\

\noindent since $u^N$ satisfies $(\ref{truncation_stopping_time_condition})$. Therefore, $Y^N$ satisfies $(\ref{truncation_stopping_time_condition})$ too, and has a bounded total jump rates when stopped at $\tau$. Applying Proposition 2.1 of \cite{Kurtz1971} to $Y^N$ ends the proof. $\square$\\

\noindent \underline{Proof of \cref{Second_level_acccompanying_martingales}}.\label{proof_about_accompanying_martingales_2} Set $\displaystyle \varphi_j:=\varphi\left(\frac{j}{N}\right)$ for  $1\leq j\leq N$. Since $Z_C^N$ has the same jumps as $u_C^N$,\par 

$\disp\hspace{1cm} \delta\left\langle Z_C^N(s),\varphi\right\rangle_2=\delta\left[\frac{1}{N}\sum_{j=1}^NZ_C^N\left(s,\frac{j}{N}\right)\varphi\left(\frac{j}{N}\right)\right]$\par 

$\disp\hspace{4.3cm}=\hspace{0.5cm}\frac{1}{N}\sum_{j=1}^N\delta Z_j^{N,C}(s)\varphi_j\hspace{0.5cm}=\hspace{0.5cm}\frac{1}{N}\sum_{j=1}^N\delta u_j^{N,C}(s)\varphi_j$.\\

\noindent Recall that $\displaystyle \left[\delta u_j^{N,C}(s)\right]\left[\delta u_k^{N,C}(s)\right]=0$ for $k\notin\{j-1,j,j+1\}$. Thus, \\

\noindent $\displaystyle\hspace{0.2cm} \left[\delta \left\langle Z_C^{N}(s),\varphi\right\rangle_2\right]^2=\frac{1}{N^2}\sum_{j=1}^N\left(\left[\delta u_{j-1}^{N,C}(s)\right]\left[\delta u_{j}^{N,C}(s)\right]\varphi_{j-1}\varphi_j+\left[\delta u_{j}^{N,C}(s)\right]^2\varphi_j^2\right.$\par 

$\displaystyle \hspace{6cm}\left.+\left[\delta u_{j}^{N,C}(s)\right]\left[\delta u_{j+1}^{N,C}(s)\right]\varphi_j\varphi_{j+1}\right)$.\par 

\noindent By \cref{First_level_acccompanying_martingales},\\

\noindent $\displaystyle\hspace{0.2cm} \sum_{s\leq t\wedge\tau}\left[\delta \left\langle Z_C^{N}(s),\varphi\right\rangle_2\right]^2$\par 

$\disp \hspace{0.5cm}=\frac{1}{N^2\mu}\sum_{j=1}^N\left\{\varphi_j^2\int_0^{t\wedge\tau}\left[N^2\left(u_{j+1}^{N,C}(s)+2u_{j}^{N,C}(s)+u_{j-1}^{N,C}(s)\right)\right.\right.$\par 
$\disp\hspace{7cm}\left.+\left(|F|_j^2+\mu\big|F_1^N\big|_j^2\right)\left(u^{N}(s)\right)\right]ds$\par 

$\displaystyle \hspace{3.5cm}\left.-\varphi_j\varphi_{j-1}\int_0^{t\wedge\tau}N^2\left(u_j^{N,C}(s)+u_{j-1}^{N,C}(s)\right)ds\right.$\par 

$\displaystyle \hspace{4.5cm}\left.-\varphi_j\varphi_{j+1}\int_0^{t\wedge\tau}N^2\left(u_j^{N,C}(s)+u_{j+1}^{N,C}(s)\right)ds\right\}\hspace{0.1cm}+\hspace{0.1cm}M(t)$,\\

\noindent where $M(t)$ is a martingale. Therefore, using the $1-$periodicity of our processes with respect to the space variable and proceeding to a change of superscript in $j$, we get\\

\noindent $\displaystyle\hspace{0.2cm} \sum_{s\leq t\wedge\tau}\left[\delta \left\langle Z_j^{N,C}(s),\varphi\right\rangle_2\right]^2$\par 

$\displaystyle \hspace{0.5cm}=\frac{1}{N^2\mu}\int_0^{t\wedge\tau}\left[\sum_{j=1}^Nu_j^{N,C}(s)N^2\big(\varphi_{j-1}^2+2\varphi_j^2+\varphi_{j+1}^2-2\varphi_{j-1}\varphi_{j}-2\varphi_{j}\varphi_{j+1}\big)\right.$\par 
$\displaystyle \hspace{4cm}\left.+\sum_{j=1}^N\left(|F|_j^2+\mu\big|F_1^N\big|_j^2\right)\left(u^{N}(s)\right)\varphi_j^2\right]ds\hspace{0.1cm}+\hspace{0.1cm}M(t)$\par 

$\displaystyle \hspace{0.5cm}=\frac{1}{N\mu}\int_0^{t\wedge\tau}\left[\frac{1}{N}\sum_{j=1}^Nu_j^{N,C}(s)\left(\left[N\big(\varphi_{j-1}-\varphi_j\big)\right]^2+\left[N\big(\varphi_{j+1}-\varphi_j\big)\right]^2\right)\right.$\par 
$\displaystyle \hspace{4cm}\left.+\frac{1}{N}\sum_{j=1}^N\left(|F|_j^2+\mu\big|F_1^N\big|_j^2\right)\left(u^{N}(s)\right)\varphi_j^2\right]ds\hspace{0.1cm}+\hspace{0.1cm}M(t)$\par 

$\displaystyle \hspace{0.5cm}=\frac{1}{N\mu}\int_0^{t\wedge\tau}\left[\left\langle u_C^{N}(s),\left(\nabla_N^+\varphi\right)^2+\left(\nabla_N^-\varphi\right)^2\right\rangle_2\right.$\par 
$\disp \hspace{4cm}\left.+\left\langle\left(|F|^2+\mu\big|F_1^N\big|^2\right)\left(u^{N}(s)\right),\varphi^2\right\rangle_2\right] ds+M(t)$.\\

\noindent This shows that (\hyperlink{Mg5}{Mg5}) defines a martingale. The case of (\hyperlink{Mg6}{Mg6}) is treated similarly, replacing $|F|^2$ by $\displaystyle \left|G^N\right|^2$. There is no term with the discrete gradient. $\square$

\subsection{Regularity and convergence of the debits}\label{proof_of_debit_function_regularity_and_convergence}

\noindent \underline{Proof of \cref{debit_function_regularity_and_convergence}}. Denote by $L_\lambda$ the Lipschitz constant common to all the reactions rates. Fix $u=(u_C,u_D),v=(v_C,v_D)\in C_p(I)\times C_p(I)$.\\

We start with the debit $F$ defined by $(\ref{Deterministic_continuous_debit_function})$.\\

$\displaystyle \left|F(u(x))-F(v(x))\right|$\par 

$\disp \hspace{1cm}\leq \sum_{r\in\mathfrak{R}_C}\left|\gamma_r^C\right|\left|\lambda_r(u_C(x))-\lambda_r(v_C(x))\right|+\sum_{r\in S_1}\left|\gamma_r^C\right|\left|\lambda_{r}(u(x))-\lambda_{r}(v(x))\right|$\par 

$\displaystyle \hspace{1cm}\leq L_\lambda\left(\left|u_C(x)-v_C(x)\right|\sum_{r\in\mathfrak{R}_C}\left|\gamma_r^C\right|+\left|u(x)-v(x)\right|\sum_{r\in S_1}\left|\gamma_r^C\right|\right)$\vspace{0.1cm}\par 

$\displaystyle \hspace{1cm}\leq L_\lambda\bar{\gamma}_C\left[\left|u_C(x)-v_C(x)\right|+\left|u_D(x)-v_D(x)\right|+\left|u_D(x)-v_D(x)\right|\right]$\vspace{0.1cm}\par 

$\displaystyle \hspace{1cm}\leq L_F\left(\left\Vert u_C-v_C\right\Vert_\infty+\left\Vert u_D-v_D\right\Vert_{\infty}\right)\hspace{0.5cm}=\hspace{0.5cm}L_F\left\Vert u-v\right\Vert_{\infty,\infty}$.\\

\noindent Thus, $F$ is Lipschitz. The Lipschitz property of the debit $F_1^N$ given by $(\ref{debit_function_retated_to_slow_onsite_involving_C})$ immediately follows from that of the reaction rates, as $F_1^N$ is a linear combination of some of them. 

The function $g$ defined by \eqref{Deterministic_debit_function_related_to_D}  is Lipschitz. Indeed,\\

$\displaystyle\hspace{0cm} \left| g\big(u(x)\big)- g\big(v(x)\big)\right|$\vspace{0.1cm}\par 

$\disp\hspace{1cm} \leq \sum_{r\in\fR_{DC}\backslash S_1}\left|\gamma_r^D\right|\left|\lambda_{r}\big(u(x)\big)-\lambda_{r}\big(v(x)\big)\right|+\sum_{r\in\mathfrak{R}_{D}}\left|\gamma_r^D\right|\left|\lambda_{r}\big(u_D(x)\big)-\lambda_{r}\big(v_D(x)\big)\right|$\vspace{0.1cm}\par 

$\disp \hspace{1cm}\leq L_\lambda\bar{\gamma}_D\left[|u_C(x)-v_C(x)|+|u_D(x)-v_D(x)|+|u_D(x)-v_D(x)|\right]$\vspace{0.1cm}\par 

$\disp \hspace{1cm}\leq L_g\big(\Vert u_C-v_C\Vert_\infty+\Vert u_D-v_D\Vert_\infty\big)\hspace{0.1cm}=\hspace{0.1cm}L_g\Vert u-v\Vert_{\infty,\infty}$,\\


\noindent and $g$ is Lipschitz. As a result, it follows easily that $G$ is Lipschitz. \\

Recall \[ G^N(u)=\sum_{i,j}\gamma_{ij}^Ng_{ij}(u)\1_i,\hspace{0.5cm}\text{ where}\hspace{0.5cm} g_{ij}(u)=\sum_{r}\gamma_r^D\theta_{ij}^r(u_i)\lambda_r(u_i), \] with $j,i=1,\cdots,N$ and $r\in(\cR_{DC}\backslash S_1)\cup \cR_D$. Fix $j$, $i$, let $x\in I_i$ and let $L_\theta$ be a Lipschitz constant for the function $\theta$. We successively have 
\begin{align*}
|\theta_{ij}^r(u_i)-\theta_{ij}^r(v_i)| & \leq \left|\theta\big(u_i^C+\frac{\gamma_r^C}{\mu}\gamma_{ij}^N\big)-\theta\big(v_i^C+\frac{\gamma_r^C}{\mu}\gamma_{ij}^N\big)\right|\theta\big(u_i^C+\gamma_r^D\gamma_{ij}^N\big)\\
& \hspace{1cm}+\theta\big(v_i^C+\frac{\gamma_r^C}{\mu}\gamma_{ij}^N\big)\left|\theta\big(u_i^C+\gamma_r^D\gamma_{ij}^N\big)-\theta\big(v_i^C+\gamma_r^D\gamma_{ij}^N\big)\right|\\
& \leq L_\theta(|u_i^C-v_i^C|+|u_i^D-v_i^D|) \leq L_\theta \Vert u-u\Vert_{\infty,\infty},
\end{align*}
\begin{align*}
|g_{ij}(u)-g_{ij}(v)| & \leq \sum_r|\gamma_r^D|\left[|\theta_{ij}^r(u_i)-\theta_{ij}^r(vi)|\lambda_r(u_j)+\theta_{ij}^r(v_i)|\lambda_r(u_j)-\lambda_r(v_j)|\right]\\
& \leq L_\theta\bar{\gamma}\bar{\lambda}(\tilde{\rho}_T)\Vert u-v\Vert_{\infty,\infty}+L_\lambda\bar{\gamma}|u_j-v_j| \bar{c}\Vert u-v\Vert_{\infty,\infty},
\end{align*}
and it easily follows that $G^N$ is Lipschtiz uniformly in $N$.\\

Since $G^N$ and $G$ are Lipschitz with a uniform constant and using a density argument, we may assume that $u=(u_C,u_D)\in C^1(I)\times C^1(I)$. Also, $G$ and $G^N$ are continuous with respect to $a$ for the $L^1$ topology. Indeed, since $\theta$ is bounded by $1$ and $g$ is a bounded function
$$
\|G^N(u)\|_\infty \le \bar{\gamma}\bar{\lambda}(\tilde{\rho}_T) \max_{i=1,\dots,N} \sum_{j=1}^NN\int_{I_j}\int_{I_i}a\left(z-\frac{j}{N}\right)dzdy= \|g\|_\infty \|a\|_{L^1(I)}.
$$
The same bound holds for $G$ and by linearity, we deduce the continuity property. 
By density of periodic $C^1(I)$ function in $L^1(I)$, we may assume that 
$a\in C^1(I)$.

\noindent Write for $x\in I_i$, and $r$ varying over $(\cR_{DC}\backslash S_1)\cup\cR_D$: 
\[ G^N(u)(x)=\sum_{j=1}^N\sum_{r}\gamma_r^D\theta\big(u_i^C+\frac{\gamma_r^C}\mu \gamma_{ij}^N)\big)\theta\big(u_i^D+\gamma_r^D \gamma_{ij}^N)\big)N\int_{I_j}\lambda_r(u_j)\int_{I_i}a\left(z-\frac{j}{N}\right)dzdy, \] 
and
\[ G(u)(x)=\sum_{j=1}^N\sum_{r}\gamma_r^D\theta(u_C(x))\theta(u_D(x))N\int_{I_j}\lambda_r(u(y))\int_{I_i}a(x-y) dzdy. \] 

\noindent Thus, since $\theta$ is bounded by $1$, $a$ by $a(0)$ and is a bounded function, \\ 

$\displaystyle\hspace{0cm} \left|G^N(u)(x)-G(u)(x)\right|$\par 

$\displaystyle \hspace{1cm} \leq a(0)\sum_r|\gamma_r^D|\sum_{j=1}^N\int_{I_j}\left|g(u_j)-g\big(u(y)\big)\right|dy$\par 

$\displaystyle \hspace{2.5cm}+a(0) \bar{\gamma}\bar{\lambda}(\tilde{\rho}_T)  \left|\theta\big(u_i^C+\frac{\gamma_r^C}\mu \gamma_{ij}^N)\big)\theta\big(u_i^D+\gamma_r^D \gamma_{ij}^N)\big)
-\theta(u^C(x))\theta(u^D(x))\right|$\par 

$\displaystyle \hspace{4cm}+ \bar{\gamma}\bar{\lambda}(\tilde{\rho}_T) \sum_{j=1}^NN\int_{I_j}\int_{I_i}\left|a\left(z-\frac{j}{N}\right)-a(x-y)\right|dz\,dy$\\

%
%
%

This clearly goes to zero when $N\to\infty$ thanks to \eqref{Truncated_process_jump_boundedness} and the Lipschitz property of $\theta$, $g$, $a$ and $u$. $\square$

\subsection{On the Gronwall-Bellman argument}\label{proof_related_to_the_Gronwall-Bellman_argument}

\noindent \underline{Proof of \cref{Blount_scalar_product_estimate}}. From \cref{Discrete_Laplacian_properties} and the observations at the beginning of \cref{approximations}, \\

$\displaystyle \hspace{1cm}\left\langle\left(\nabla_N^\pm T_N(t)f\right)^2,1\right\rangle_2=\left\langle T_N(2t)f,\Delta_Nf\right\rangle_2$\par 

$\displaystyle \hspace{4.5cm}=\sum_{m}\left[\left\langle \left|f,\varphi_{m,N}\right\rangle_2\right|^2+\left\langle \left|f,\psi_{m,N}\right\rangle_2\right|^2\right]e^{-2\beta_{m,N}t}\beta_{m,N}$\par 

$\displaystyle \hspace{4.5cm}=\sum_{m}\left(\varphi_{m,N}^2\left(\frac{j}{N}\right)+\psi_{m,N}^2\left(\frac{j}{N}\right)\right)e^{-2\beta_{m,N}t}\beta_{m,N}$\par 

$\disp \hspace{4.5cm}\leq 2\sum_me^{-2\beta_{m,N}t}\beta_{m,N}$.\\

\noindent Similarly, \[ \left\langle \left(T_N(t)f\right)^2,1\right\rangle_2\leq 1+2\sum_{m>0}e^{-2\beta_{m,N}t}. \] The result then holds for 
\begin{equation}\label{Blount_inner_product_estimate_constant}
h_N(t)=1+4\sum_{m>0}e^{-2\beta_{m,N}t}\left(\beta_{m,N}+1\right),
\end{equation}
since $\beta_{0,N}=0$ and $\beta_{m,N}>cm^2$ for $m>0$ and $c>0$, where $c$ is independent of $m$ and $N$. $\square$\\

\noindent \underline{Proof of \cref{Blount_expectation_estimate}}. Let $f(x)=e^x$ and note \[ 0\leq f"(x+y)=f(x)f(y)\leq 3f(x)\hspace{0.5cm}\text{if}\hspace{0.3cm}|y|\leq 1. \] Using change of variables for functions of bounded variation, we have for $t_0\leq t\leq t_1$,\\

$\displaystyle\hspace{1cm} f(m(t))=1+\int_{t_0}^{t}f'\left(m(s^-)\right)dm(s)$\par 

$\disp\hspace{4.5cm}+\sum_{t_0\leq s\leq t}\left[f(m(s))-f\left(m(s^-)\right)-f'\left(m(s^-)\right)\delta m(s)\right]$\par 

$\displaystyle \hspace{2.5cm}\leq 1+\int_{t_0}^{t}f'\left(m(s^-)\right)dm(s)+\frac{3}{2}\sum_{t_0\leq s\leq t}f\left(m(s^-)\right)(\delta m(s))^2$,\\

\noindent by Taylor's theorem, (ii) of \cref{Blount_expectation_estimate} and our observations at the start of the proof. Note that $\displaystyle \int_{t_0}^tf'\left(m(s^-)\right)dm(s)$ has mean $0$ and after applying (iii) of \cref{Blount_expectation_estimate} and taking expectations, we have \[ \mathds{E}f(m(t))\leq 1+\frac{3}{2}\int_{t_0}^t\mathds{E}f(m(s))h(s)ds. \] The result then follows form Gronwall's inequality. $\square$\\

\noindent \underline{Proof of \cref{third_level_accompanying_martingale}}. We want to identify the martingale part of $\disp \sum_{s\leq t}\big[\delta\bar{m}_C(t)\big]^2$. Integrating $\disp \bar{m}_C(t)$ by parts leads to\\

$\disp \hspace{1cm}\bar{m}_C(t)=\left\langle T_N(\bar{t}-t)Z_C^N(t\wedge\tau)+\int_0^t\Delta_NT_N(\bar{t}-s)Z_C^N(s\wedge\tau)ds,f\right\rangle_2$\par 

$\disp \hspace{2.2cm}=\left\langle T_N(\bar{t}-t)Z_C^N(t\wedge\tau),f\right\rangle_2+\left\langle\int_0^t\Delta_NT_N(\bar{t}-s)Z_C^N(s\wedge\tau)ds,f\right\rangle_2$.\\

\noindent Since $\displaystyle t\mapsto \int_0^t\Delta_NT_N(\bar{t}-s)Z_C^N(s\wedge\tau)ds$ is continuous, the corresponding term does not jump and consequently
\[ \delta\bar{m}_C(t)=\delta\left\langle T_N(\bar{t}-t)Z_C^N(t\wedge\tau)\hspace{0.1cm},\hspace{0.1cm}f\right\rangle_2=\big[\delta T_N(\bar{t}-t)u_C^N(t\wedge\tau)\big]_j, \]
where $j$ has been fixed within $\{1,\cdots,N\}$. Observing that \[ [T_N(\bar{t}-t)u_C^N(t)]_j=\sum_{i=1}^Nu_i^{N,C}(t)[T_N(\bar{t}-t)\1_i]_j, \] we have \par 

$\disp \big[\delta T_N(\bar{t}-t)u_C^N(t)\big]_j^2=\left(\sum_{j=1}^N\left[\delta u_i^{N,C}(t)\right][T_N(\bar{t}-t)\1_i]_j\right)^2$\par 

$\disp\hspace{2cm} =\sum_{i=1}^N\left(\left[\delta u_i^{N,C}(t)\right]^2[T_N(\bar{t}-t)\1_i]_j^2\right.$\par 

$\disp \hspace{3.5cm}+\left.\left[\delta u_i^{N,C}(t)\right]\left[\delta u_{i-1}^{N,C}(t)\right][T_N(\bar{t}-t)\1_i]_j[T_N(\bar{t}-t)\1_{i-1}]_j\right.$\\

$\disp \hspace{4cm}+\left.\left[\delta u_i^{N,C}(t)\right]\left[\delta u_{i+1}^{N,C}(t)\right][T_N(\bar{t}-t)\1_i]_j[T_N(\bar{t}-t)\1_{i+1}]_j\right)$\\

\noindent We need the forthcoming.

\begin{lemma}\label{change_between_index_components_for_the_discrete_semigroup}
We have $\disp [T_N(t)\1_i]_j=[T_N(t)\1_j]_i$ for all $1\leq i\leq N$.
\end{lemma}
\noindent \textit{Proof.} Let $1\leq i\leq N$ and $t\geq0$ be fixed. We first remark that \[ \langle \1_i,\varphi\rangle_2=\varphi_j=v(j/N)\hspace{0.5cm}\forall\varphi\in\H^N. \] Then, from the spectral decomposition of $T_N(t)$ on $\disp\big(\H^N,\langle\cdot,\cdot\rangle_2\big)$ (see  \cref{Discrete_Laplacian_properties}),\\

$\disp\hspace{0.5cm} [T_N(t)\1_i]_j=\sum_m e^{-\beta_{m,N}t}\left(\big\langle \1_i,\varphi_{m,N}\big\rangle_2\varphi_{m,N}^j+\big\langle \1_i,\psi_{m,N}\big\rangle_2\psi_{m,N}^j\right)$\par 

$\disp \hspace{2.3cm}=\sum_m e^{-\beta_{m,N}t}\left(\varphi_{m,N}^i\big\langle \1_j,\varphi_{m,N}\big\rangle_2+\psi_{m,N}^i\big\langle \1_j,\psi_{m,N}\big\rangle_2\right)=[T_N(t)\1_j]_i$. $\square$\\

\noindent Now, from \cref{First_level_acccompanying_martingales}, there exists a martingale $M(t)$ such that, using periodicity and \cref{change_between_index_components_for_the_discrete_semigroup}, we get\\

\noindent $\disp\hspace{0.2cm} \sum_{s\leq t}\big[\delta T_N(\bar{t}-s)u_C^N(t)\big]_j^2$\par 

$\disp \hspace{0.3cm}=\sum_{i=1}^N\frac{1}{\mu}\int_0^{t\wedge\tau}\left\{\left[N^2\big(u_{i+1}^{N,C}(s)+2u_{i}^{N,C}(s)u_{i-1}^N(s)\big)\right.\right.$\par 
$\disp\hspace{5.5cm}\left.+\left(|F|_i^2+\mu\big|F_1^N\big|_i^2\right)\left(u^{N}(s)\right)\right][T_N(\bar{t}-s)\1_i]_j^2$\\
 
$\disp \hspace{3.2cm}-N^2\big(u_{i}^{N,C}(s)+u_{i-1}^{N,C}(s)\big)[T_N(\bar{t}-s)\1_{i}]_j[T_N(\bar{t}-s)\1_{i-1}]_j$\vspace{0.2cm}\par 

$\disp \hspace{2.5cm}\left.-N^2\big(u_{i}^{N,C}(s)+u_{i+1}^{N,C}(s)\big)[T_N(\bar{t}-s)\1_{i}]_j[T_N(\bar{t}-s)\1_{i+1}]_j\right\} ds\hspace{0.1cm}+M(t)$\\

$\disp \hspace{0.2cm}=\frac{1}{\mu}\int_0^{t\wedge\tau}\left\{\sum_{i=1}^Nu_i^{N,C}(s)N^2\text{\LARGE(}[T_N(\bar{t}-s)\1_{i-1}]_j^2+2[T_N(\bar{t}-s)\1_{i}]_j^2+[T_N(\bar{t}-s)\1_{i+1}]_j^2\right.$\par 
$\disp \hspace{5.5cm}-2[T_N(\bar{t}-s)\1_{i}]_j[T_N(\bar{t}-s)\1_{i-1}]_j$\vspace{0.2cm}\par 
 
$\disp \hspace{6.5cm}-2[T_N(\bar{t}-s)\1_{i+1}]_j[T_N(\bar{t}-s)\1_{i}]_j\text{\LARGE)}$\\

$\disp \hspace{2.3cm}\left.+\sum_{i=1}^N\left(|F|_i^2+\mu\big|F_1^N\big|_i^2\right)\left(u^{N}(s)\right)[T_N(\bar{t}-s)\1_{i}]_j^2\right\} ds\hspace{0.1cm}+ M(t)$\\

$\disp \hspace{0.2cm}=\frac{1}{\mu N}\int_0^{t\wedge\tau}\left\{\frac{1}{N}\sum_{i=1}^N\text{\LARGE(}N[T_N(\bar{t}-s)N\1_{j}]_{i+1}-[T_N(\bar{t}-s)N\1_{j}]_{i}\big)^2\right.$\par 
$\disp\hspace{5cm}+\big(N[T_N(\bar{t}-s)N\1_{j}]_{i}-[T_N(\bar{t}-s)N\1_{j}]_{i-1}\text{\LARGE)}^2$\par 

$\disp \hspace{2.8cm}\left.+\frac{1}{N}\sum_{i=1}^N\left(|F|_i^2+\mu\big|F_1^N\big|_i^2\right)\left(u^{N}(s)\right)[T_N(\bar{t}-s)N\1_{i}]_j^2\right\}ds\hspace{0.2cm}+\hspace{0.2cm} M(t)$\\

$\disp \hspace{0.2cm}=\frac{1}{\mu N}\int_0^{t\wedge\tau}\left[\left\langle u_C^N(s),\big(\nabla_N^+T_N(\bar{t}-s)f\big)^2+\big(\nabla_N^-T_N(\bar{t}-s)f\big)^2\right\rangle_2\right.$\par 

$\disp \hspace{4cm}\left.+\left\langle \left(|F|^2+\mu\big|F_1^N\big|^2\right)\left(u^{N}(s)\right), \big(T_N(\bar{t}-s)f\big)^2\right\rangle_2\right] ds +M(t)$. $\square$

\vspace{0.3cm}
\footnotesize
\bibliographystyle{alpha}
\bibliography{MSSM_Deterministic_Limit_Decoupling_Last}

\begin{thebibliography}{CDMR12}

\bibitem[Arn80]{Arnold1980bis}
Ludwig Arnold.
\newblock {\em Stochastic Systems: The Mathematics of Filtering and
  Identification and Applications}, volume~78 of {\em Nato advanced study
  institutes}, chapter 3, Mathematical models of chemical reactions.
\newblock Dordrecht, 1980.

\bibitem[AT80]{Arnold1980}
L.~Arnold and M.~Theodosopulu.
\newblock Deterministic limit of the stochastic model of chemical reactions
  with diffusion.
\newblock {\em Adv. Appl. Prob.}, 12:367--379, 1980.

\bibitem[BKPR06]{Kurtz2006}
Karen Ball, Thomas~G. Kurtz, Lea Propovic, and Greg Rempala.
\newblock Asymtotic analysis of multiscale approximations to reaction networks.
\newblock {\em The Annals of Applied Probability}, 16(4):1925--1961, 2006.

\bibitem[Blo87]{Blount1987}
Douglas~James Blount.
\newblock {\em Comparison of a stochastic model of a chemical reaction with
  diffusion and the deterministic model}.
\newblock Ph.d., The University of Wisconsin-Madison, 1987.

\bibitem[Blo92]{Blount1992}
D.~J. Blount.
\newblock Law of large numbers in the supremum norm for a chemical reaction
  with diffusion.
\newblock In {\em The Annals of Applied Probability}, volume~2, pages 131--141.
  1992.

\bibitem[Blo93]{Blount1993}
D.~J. Blount.
\newblock Limit theorems for a sequence of nonlinear reaction-diffusion
  systems.
\newblock In {\em Stochastic Processes and their Applications}, volume~45,
  pages 193--203. 1993.

\bibitem[Blo94]{Blount1994}
D.~J. Blount.
\newblock Density-dependent limits for a nonlinear reaction-diffusion model.
\newblock In {\em The Annals of Probability}, volume~22, pages 2040--2070.
  1994.

\bibitem[CDMR12]{Arnaud2012}
A.~Crudu, A.~Debussche, A.~Muller, and O.~Radulescu.
\newblock Convergence of stochastic gene networks to hybrid piecewise
  deterministic processes.
\newblock {\em The Anals of Applied Probability}, 22(5):1822--1859, 2012.

\bibitem[CDR09]{Arnaud2009}
A.~Crudu, A.~Debussche, and O.~Radulescu.
\newblock Hybrid stochastic simplifications for multiscale gene networks.
\newblock {\em BMC Systems Biology}, 3:89, Septembre 2009.

\bibitem[CH98]{Cazenave1998}
Thierry Cazenave and Alain Haraux.
\newblock {\em An Introduction to Semilinear Evolution Equations}.
\newblock Clarendon Press - Oxford, 1998.

\bibitem[Dav93]{Davis1993}
M.~H.~A. Davis.
\newblock Markov models and optimization.
\newblock In Chapman and London Hall, editors, {\em Monographs on Statistics
  and Applied Probability}, volume~49. 1993.

\bibitem[DL05]{Davidson2005}
E.~Davidson and M.~Levine.
\newblock Gene regulatory networks.
\newblock {\em Proc Natl Acad Sci USA}, 102(14), 2005.

\bibitem[EK86]{Kurtz1986}
Stewart~N. Ethier and Thomas~G. Kurtz.
\newblock {\em Markov Processes, Characterization and Convergence}.
\newblock John Wiley and Sons, Inc, 1986.

\bibitem[ET89]{P_Erdi1989}
P.~\'Erdi and J.~Th\'oth.
\newblock {\em Mathematical models of chemical reactions: theory and
  applications of deterministic and stochastic models}.
\newblock Nonlinear science, Manchester University Press, 1989.

\bibitem[Gil76]{Gillespie1976}
D.~T. Gillespie.
\newblock A general method for numerically simulating the stochstic time
  evolution of coupled chemical reactions.
\newblock {\em J. Comput. Phys.}, 22:403--434, 1976.

\bibitem[Kat66]{Kato1966}
T.~Kato.
\newblock {\em Perturbation Theory for Linear Operators}.
\newblock Springer-Verlag, Berlin, 1966.

\bibitem[Kot86a]{Kotelenez1986bis}
P.~Kotelenez.
\newblock {\em Gaussian approximation to the nonlinear reaction-diffusion
  equation.}
\newblock Report 146, Universit\"at Bremen Forschungsschwerpunkt Dynamische
  Systemes., 1986.

\bibitem[Kot86b]{Kotelenez1986}
P.~Kotelenez.
\newblock Law of large numbers and central limit theorem for linear chemical
  reactions with diffusion.
\newblock In {\em The Annals of Probability}, volume~14, pages 173--193.
  Universit\"at Bremen, 1986.

\bibitem[Kot87]{Kotelenez1987}
P.~Kotelenez.
\newblock Fluctuations near homogeneous states of chemical reaction with
  diffusion.
\newblock In {\em Adv. in Appl. Probab.}, volume~19, pages 352--370. 1987.

\bibitem[Kot88a]{Kotelenez1988}
P.~Kotelenez.
\newblock High density limit theorems for nonlinear chemical reactions with
  diffusion.
\newblock In {\em Probab. Theory Related Fields}, volume~78, pages 11--37.
  1988.

\bibitem[Kot88b]{Kotelenez1988bis}
P.~Kotelenez.
\newblock A stochastic reaction-diffusion model.
\newblock In {\em University of Ultrecht}. 1988.

\bibitem[Kui77]{Kuiper1977}
Hendrick~J. Kuiper.
\newblock Existence and comparison theorems for nonlinear diffusion systems.
\newblock {\em J. Math. Anal. and App.}, 60:166--181, 1977.

\bibitem[Kur70]{Kurtz1970}
T.~G. Kurtz.
\newblock Solutions of ordinary differential equations as limits of pure jump
  markov processes.
\newblock {\em J. Appl. Prob.}, 7:49--58, 1970.

\bibitem[Kur71]{Kurtz1971}
T.~G. Kurtz.
\newblock Limit theorems for sequences of jump markov processes approximating
  ordinary differential processes.
\newblock {\em J. Appl. Prob.}, 8:344--356, 1971.

\bibitem[MW11]{Lesley2011}
Lesley~T. MacNeil and Albertha J.~M. Walkout.
\newblock Gene regulatory networks and the role of robustness and stochasticity
  in the control of gene expression.
\newblock {\em Genome Res., Cold Spring Harbor Laboratory Press}, 21:645--657,
  2011.

\bibitem[RMC07]{Crudu2007}
O.~Radulescu, O.~Muller, and A.~Crudu.
\newblock Théorèmes limites pour des processus de markov à sauts: syntèses
  des résultats et applications en biologie moléculaire.
\newblock {\em Tech. Sci. Inform.}, 26:443--469, 2007.

\bibitem[SL05]{Stathopoulos2005}
A.~Stathopoulos and M.~Levine.
\newblock Genomic regulatory networks and animal development.
\newblock {\em Dev. Cell.}, 9(4), 2005.

\end{thebibliography}
\normalsize

\end{document}